%% file: AVR.tex
\documentclass[11pt,a4paper,leqno]{amsart}

\include{packages}

\usepackage[margin=70pt]{geometry}
\usepackage{xcolor}
\usepackage{bbm}

\renewcommand{\bar}{\overline}

\newcommand{\re}{\mathrm{Re}}

\newcommand{\II}{I\hspace{-0.1cm}I}
\newcommand{\III}{I\hspace{-0.1cm}I\hspace{-0.1cm}I}

\newcommand{\rrho}{_{\rho_1}^{\rho_2}}
\newcommand{\mrho}{_{\rho_1}^{\rho_2}(I_0^{*})}

\renewcommand{\bar}{\overline}

\newcommand{\argsinh}{{\mathrm{arsinh}}}

\newtheorem{innercustomdefi}{Definition}
\newenvironment{customdefi}[1]
{\renewcommand\theinnercustomdefi{#1}\innercustomdefi}
{\endinnercustomdefi}

\newtheorem*{definition*}{Definition}
\newtheorem*{theorem*}{Theorem}

\newcounter{notes}%

\begin{document}

	\title[]{Adapted Renormalized Volume for Hyperbolic 3-Manifolds with Compressible Boundary}

	\date{v0, \today}
	
	\author{Viola Giovannini}
	\address{Viola Giovannini:
		ETH, Department of Mathematics, 
		HG, 101 R\"amistrasse,
		8092 Z\"urich, Switzerland}
	\email{viola.giovannini@math.ethz.ch}
	
	\maketitle
	\begin{abstract}
	The renormalized volume is a smooth function associating to every convex co-compact hyperbolic $3$-manifold $M$ a real number.
	 When the boundary of $M$ is incompressible, the renormalized volume is always positive, otherwise there are sequences of convex co-compact structures on $M$ whose renormalized volumes diverge to minus infinity.
	
	We define here a new version of the renormalized volume which adapts to the compressible boundary case, satisfying  properties analogous to those of the classical one in the incompressible setting. In particular, the adapted renormalized volume is bounded from below, its differential has uniformly bounded supremum norm, and its gradient has uniformly bounded Weil-Petersson norm. Moreover, it stays at uniformly bounded distance from the convex core volume function. As a corollary, we obtain a bound on the convex core volume of handlebodies in terms of the Weil-Petersson distance from a certain subset of the Teichm\"uller space, where the convex core volume is bounded by a known constant. Furthermore, the adapted renormalized volume extends continuously, as a function on the Teichm\"uller space of $\partial\bar M$, to the strata in the boundary of its Weil-Petersson completion corresponding to compressible multicurves. We provide a geometric interpretation of the limit quantity by defining a renormalized volume, and its adapted version, for convex co-compact hyperbolic $3$-manifolds with a finite set of marked points in the boundary. 
	\end{abstract}
	
	\tableofcontents
	
		\section{Introduction and Results}
		
	The volume of closed or cusped hyperbolic $3$-manifolds has been extensively studied as a topological invariant. However, there exists a wide class of hyperbolic manifolds whose elements have infinite volume, and here is where the notion of renormalized volume comes into play. By the Tameness Theorem, every complete hyperbolic $3$-manifold with finitely generated fundamental group is homeomorphic to the interior of a compact manifold $N$, possibly with boundary. Whenever $N$ presents a boundary component of genus at least $2$, any complete hyperbolic metric on its interior has infinite volume.
	
	The main setting of this paper is the class of geometrically finite hyperbolic $3$-manifold, and, more precisely, the subset of convex co-compact ones. A hyperbolic $3$-manifold $M$ is \textit{geometrically finite} if it contains a non-empty convex subset $C$ of finite volume to which $M$ is homotopically equivalent, and it is \textit{convex co-compact} if $C$ is also compact. The smallest convex subset $C(M)$ of $M$ whose inclusion is an homotopy equivalence is called the \textit{convex core} of $M$. If $M$ is a convex cocompact strucutre on the interior of $N$, the boundary of its closure $\partial\bar M$ coincides with $\partial N$.
	
	The world of geometrically finite hyperbolic $3$-manifold is beautifully interconnected with that of Riemann surfaces and then, by Riemann Uniformization Theorem, with that of hyperbolic surfaces, or even conformal classes of Riemannian metrics.
	% This deep link is due to the fundational work of Ahlfors, Bers, Kra, Marden, Maskit, Sullivan and Thurston, among others. 
	In particular, the deformation space $CC(M)$ of convex co-compact hyperbolic structures on $M$, considered up to homotopy, is parameterized by a quotient of the space of marked hyperbolic structures, that is, the \textit{Teichm\"uller space} of $\partial\bar M$ (see \cite[Theorem $5.1.3$]{Ma2016}, or \cite[Theorem 5.27]{MT1998}). Briefly, the correspondence is realized as follows. The Riemannian metric tensor on a convex co-compact hyperbolic $3$-manifold $M$ diverges exponentially on the components of the complement of the convex core, and thus it does not extend to the boundary $\partial\bar M$. However, the manifold $M$ can be \textit{conformally compactified}, that is, the hyperbolic metric on $M$ induces a well defined conformal structure on $\partial\bar M$, which then defines a point in the Teichm\"uller space $\mathcal{T}(\partial\bar M)$. This conformal structure is referred to as the \textit{conformal boundary at infinity} of $M$, denoted $\partial_{\infty}M$, and defined as follows. 
	
	\bdefi\label{cnoformalboundary}
	Let $M$ be a convex hyperbolic $3$-manifold, and let $\Gamma<\mathbb{P}SL(2, \mathbb{C})$ be a discrete subgroup of M\"obius tranformation such that $M=\mathbb{H}^3/\Gamma$. Consider the \textit{limit set} of $M$ \[\Lambda(\Gamma)= \overline{\Gamma\cdot x}\cap \mathbb{CP}^1~,\]
where $x$ is any point of $\mathbb{H}^3$, and the Riemann sphere $\mathbb{CP}^1$ coincides with the conformal structure $\partial_{\infty}\mathbb{H}^3$ induced on the boundary of the Poincar\'e model for $\mathbb{H}^3$. The \textit{discontinuity domain} of $M$ is then defined as \[\Omega(\Gamma)=\mathbb{CP}^1\setminus \Lambda(\Gamma)~,\] and the \textit{conformal boundary at infinity} of $M$ as  \[\partial_{\infty}M=\Omega(\Gamma)/\Gamma~.\]
	\edefi
	
	We remark that the action of $\Gamma$ extends properly discontinuously to $\Omega(\Gamma)$, and that then $\partial_{\infty}M$ is a closed surface
	homeomorphic to $\partial\bar M$. Moreover, since $\Omega(\Gamma)$ lies in $\mathbb{CP}^1$ and $\Gamma$ is in particular a subgroup of M\"obius transformation, the conformal boundary at infinity is naturally equipped with a \textit{complex projective structure}, which in turn induces its underlying complex structure.

	\bthm \label{UnifCCM}(Ahlfors, Bers, Kra, Marden, Maskit, Sullivan, Thurston, \cite[Theorem 5.1.3]{Ma2016}).
	The deformation space $CC(M)$ of convex co-compact hyperbolic structures on $M$ considered up to homotopy is biholomorphically parameterized by a quotient of the space of conformal structures on the boundary at infinity, as 
	\[ CC(M)=\quotient{\mathcal T(\partial\overline M)}{T_0(D)},\]
	where $\mathcal T(\partial\overline M)$ is the product of the Teichm\"uller spaces of the connected components of $\partial \bar{M}$, and $T_0(D)$ the subgroup of the mapping class group of $\partial\bar M$ generated by all the \textit{Dehn twists} along compressible simple closed curves in $\partial\bar M$
	\ethm
	
	Within this framework, given a convex co-compact hyperbolic $3$-manifold $M$, the \textit{renormalized volume} is a real-analytic function on $CC(M)$, or, equivalently, on the Teichm\"uller space of $\partial\bar M$ \cite{S2008} \[V_R\colon \mathcal{T}(\partial\bar M)\rightarrow \mathbb{R}~.\]
	 The central idea behind its definition is to renormalize the divergent volumes associated to an exhaustion of $M$ by homotopically equivalent convex compact submanifolds. There are many possible choices for such an exhaustion. Indeed, there is one for each smooth metric in the conformal boundary at infinity \cite{epsteinSurf}. The renormalized volume of $M$ is defined using the exhaustion in convex compacts determined by the unique hyperbolic metric conformal to $\partial_{\infty}M$. For all other choices of exhaustion, or, equivalently, of conformal metric, the resulting quantity is referred to as a \textit{W-volume.} 
	
	In this paper we focus in the setting where the boundary of $M$ is \textit{compressible}, that is, there exists a non-trivial simple closed curve in $\partial\bar M$ which bounds a disk in $M$. Equivalently, the inclusion $\partial\bar M\hookrightarrow \bar M$ is not $\pi_1$-injective.
	If $\partial\bar{M}$ is incompressible, Bridgeman, Brock and Bromberg proved in \cite{BBB2018} that $V_R$ is always non-negative and, if $M$ is also (relatively) acylindrical, it attains its minimum at the unique convex co-compact metric with totally geodesic convex core boundary, which, in this setting, always exists and it is unique. In contrast, if $\partial\bar{M}$ is compressible, the works of Canary and Bridgeman in \cite{BC2017}, and of Schlenker and Witten in \cite{averages}, show that the renormalized volume of a sequence of convex co-compact manifolds associated to the \textit{pinching} of a compressible curve in the boundary diverges to $-\infty$. Pinching a simple closed curve in a Riemann surface means making the length of its geodesic representative go to zero with respect to the hyperbolic metric.
	
	The differential of the renormalized volume has a beautiful complex analytic description. Given $M=\mathbb{H}^3/\Gamma$ convex co-compact, to the natural complex projective structure on $\partial_{\infty}M$ is associated its \textit{developing map} $f$ (see \cite{dumas-survey}). In this work, we always consider the developing map comparing the complex projective structure of $\partial_{\infty}M$ with its Fuchsian one. In this way, for any connected component $X$ of $\partial_{\infty}M$, thanks to the Riemann Uniformization Theorem $f$ is interpreted as a locally injective holomorphic $(\Gamma', \Gamma_X)$-equivariant map from $\mathbb{H}^2$ to $\mathbb{CP}^1$, with $\Gamma'< \mathbb{P}SL(2, \mathbb{R})$ such that $\partial_{\infty}M$ is biholomorphic to $\mathbb{H}^2/\Gamma'$, and $\Gamma_X<\Gamma$ the stabilizer of the connected component of $\widetilde{X}\subseteq \Omega(\Gamma)$ image of $f$. We denote by $f_M$ the union of these maps, and we will simply call it the \textit{developing map of $\partial_{\infty}M$}. The \textit{Schwarzian derivative} of $f$ is the holomorphic quadratic differential \[\mathcal S(f)=\left(\left(\dfrac{f^{''}}{f^{'}}\right)^{'}-\dfrac{1}{2}\left(\dfrac{f^{''}}{f^{'}}\right)^{2}\right)dz^2~.\] An important property of the Schwarzian derivative is that $\mathcal S (f)$ is zero if and only if $f$ is the restriction of a M\"obius transformation. 
	We also need to point out that the space of \textit{holomorphic quadratic differentials} $Q(X)$ on a Riemann surfaces $X\in \mathcal T(\partial\bar{M})$ is identified with the cotangent bundle of the Teichm\"uller space $T^{*}_{X}\mathcal T(\partial\bar{M})$ (see \cite[Sections $6.6$ and $6.7$]{Hubbard2016}, or \cite[Chapter $6$]{QTT}). Moreover, the natural pairing $\langle\cdot, \cdot\rangle$ between $Q(X)$ and the space of \textit{harmonic Beltrami differentials} $B(X)$ on $X$, identifies the latter with the tangent bundle $T_{X}\mathcal T(\partial\bar{M})$.
	
	\bthm\label{dVRdiff}(\cite{TZUnif}, \cite[Corollary $8.6$]{S2008}, \cite[Corollary $3.11$]{compare}) \label{dVR}
	Let $M$ be a convex co-compact hyperbolic $3$-manifold, and let $\mathcal S(f_M)$ be the holomorphic quadratic differential on $\partial_{\infty}M$ given by the Schwarzian derivative of the developing map $f_M$ of $\partial_{\infty}M$. Then, the differential of the renormalized volume at $[\partial_{\infty}M]\in\mathcal{T}(\partial\bar M)$ is
	\[ dV_R(\cdot) =\re\langle \mathcal S(f_M), \cdot \rangle~.\] 
	\ethm
	
	As a consequence of this theorem, when $M$ has incompressible boundary, the Kraus-Nehari bound \cite{Ne1949} for the supremum norm of the Schwarzian derivative of a univalent map implies that the \textit{Teichm\"uller norm} of (minus) the \textit{Weil-Petersson gradient} of the renormalized volume is uniformly bounded. Since $\mathcal{T}(\partial\bar M)$ equipped with the Teichm\"uller norm is complete, the gradient flow exists at any time. In the (relatively) acylindrical case, this has been widely studied in \cite{BBB2018},  \cite{bridgeman-brock-bromberg:gradient} and \cite{bridgeman-bromberg-pallete:convergence}, where it is proven, in particular, that the flowlines always converge to the point whose associated convex cocompact structure on $M$ has totally geodesic convex core boundary. On the other hand, the norm of the Schwarzian derivative of non-univalent map $f$ with image a hyperbolic domain in $\mathbb{CP}^1$ is bounded from below by a diverging function of the hyperbolic length of the shortest non-trivial simple closed curve in its image \cite{KM}. When $f$ is the developing map associated to a convex co-compact manifold $M$, this corresponds to the situation where the boundary at infinity $\partial_{\infty}M$ contains a compressible simple closed curve of small hyperbolic length $\ell$. 
	Inspired by the asymptotic behaviour of the renormalized volume under the pinching of a compressible curve in the boundary (see \cite[Theorem $A.15$]{averages}, or also \cite[Theorem $1.4$]{CG}), we define the adapted renormalized volume by subtracting to the standard one the divergent terms as follows.

	\begin{customdefi}{\ref{adaptedVR}}
	
	Given $M$ a convex co-compact hyperbolic 3-manifold,
	we define the \textit{adapted renormalized volume} as the function \[\widetilde{V_R}\colon \mathcal{T}(\partial\bar{M})\rightarrow \mathbb{R}\] such that \[\widetilde{V_R}(X)=V_R(X)+ \max_{m\in\Gamma^{comp}(\partial\bar M)} L(X, m)\] where  $\Gamma^{comp}(\partial\bar M)$ denotes the set of compressible multicurves on $\partial\bar M$ and 
	\[ L(X, m)=\pi^3\sum_{\substack{\gamma\in m }} \f{1}{\ell_{\gamma}(X)} \]
	with $\ell_{\gamma}(\cdot)$ the hyperbolic length function of $\gamma$.

	\end{customdefi}
	
	We say that a \textit{multicurve} is \textit{compressible} if each of its components is compressible.
	We start by observing that the adapted renormalized volume is well defined, continuous, and Lipschitz with respect to the Weil-Petersson metric (see Lemma \ref{maximumexistence}, Corollary \ref{contLip}, Remark \ref{welldef}, and Corollary \ref{WPdistance}). Moreover, outside the subset of Teichn\"uller space where the maximum of $L(X, \cdot)$ as in the definition is realized by multiple multicurves, it is also smooth.
	
	Exploiting \cite[Theorem $1.4$]{CG}, we are able to bound the adapted renormalized volume from below.
	
	\begin{customthm}{\ref{AVRbounded}}
		For every convex co-compact hyperbolic $3$-manifold $M$, the adapted renormalized volume $\widetilde{V_R}(\cdot)$ is bounded from below by a constant depending just on the topology of the boundary $\partial\bar M$.
	\end{customthm}
	
	%Recall that the differential of the renormalized volume at a point $X\in\mathcal{T}(\partial\bar M)$ coincides with the real part of the Schwarzian derivative of $f_M$, with $f_M$ the developing map of the natural projective structure on the boundary at infinity $\partial_{\infty}M$, with $M=M(X)$ the convex co-compact manifold associated to $X$ (see Section \ref{cb}, Section \ref{SD} , and Theorem \ref{UnifCCM}).
	%As already discussed in Chapter \ref{3}, the infinity norm of the Schwarzian derivative of a non-univalent map, with image a hyperbolic domain in $\mathbb{CP}^1$, diverges whenever the hyperbolic length of a non-contractible simple closed curve in its image goes to zero, \cite{KM}. When $f$ is the developing map associated to a convex co-compact manifold $M$, this corresponds to the situation where the hyperbolic length $\ell$ of a compressible simple closed curve in $\partial_{\infty}M$ is going to zero. 
	From the explicit behaviour of the Schwarzian derivative on long compressible tubes furnished by \cite[Theorem $1.1$]{CG}, one can easily see that both the $L^1$ norm and the $L^{\infty}$ norm of the differential of the renormalized volume diverge as $\ell\rightarrow 0$, with $\ell$ the hyperbolic length of a compressible simple closed curve in $\partial_{\infty}M$. This divergence does not occur for the differential of the adapted renormalized volume.
	
	\bthm(Theorems \ref{boundedGradient} and \ref{boundedGradientinfty})\label{L1}
		Let $M$ be a convex co-compact hyperbolic manifold.
		At the points where it exists, the differential of the adapted renormalized volume is bounded in $L^1$ and $L^{\infty}$ norms by a constant that depends only on the topology of the boundary $\partial\bar M$.
	\ethm
	
	The definition of the $L^p$ norms of holomorphic quadratic differential is furnished in Section \ref{norms}.
	
	\begin{customcor}{\ref{boundedGradientWP}}
		Let $M$ be a convex co-compact hyperbolic manifold.
		At the points where it exists, the Weil-Petersson gradient of the adapted renormalized volume has Weil-Petersson norm bounded by a constant that depends only on the topology of the boundary $\partial\bar M$.
	\end{customcor}
	
	As another corollary of Theorem \ref{L1}, we get an analogous bound to the one given in \cite[Therem $1.2$]{compare} for acylindrical manifolds, for the difference of the renormalized volumes of two convex co-compact structures in terms of the Weil-Petersson distance between their boundaries.
	
		\begin{customcor}{\ref{WPdistance}}
	Let M be a convex co-compact hyperbolic manifold. For any couple of points $X_0, X_1\in \mathcal{T}(\partial\bar M)$, the difference of the corresponding adapted renormalized volumes is bounded by the Weil-Petersson distance between $X_0$ and $X_1$ as \[\abs{\widetilde{V_R}(X_1)-\widetilde{V_R}(X_0)}\leq \sqrt{2\pi\abs{\chi(\partial\bar M)}}C(\partial\bar M)d_{WP}(X_0, X_1)~,\]
	where $C(\partial\bar M)$ is a constant depending only on the topology of $\partial\bar M$.
	\end{customcor}
	
	 In Theorem \ref{Vcdistance}, we also prove that the difference between the adapted renormalized volume and the convex core volume function is bounded by a constant, again depending only on the topology of the boundary. As a corollary of this theorem, we obtain a similar upper-bound to the one in \cite[Theorem $C$]{bridgeman-brock-bromberg:gradient} for the convex core volume $V_C(M(\cdot))$ of acylindrical manifolds, in the case of a handlebody.
	 
	 	\begin{customcor}{\ref{Vcestimate}}
	 Let $M$ be a convex co-compact genus $g$ handlebody, and let $X\in\mathcal{T}(\partial\bar M)$. Let $\Gamma_0^{comp}(\partial \bar M)$ denote the set of compressible multicurves on $\partial \bar M$ that decompose $M$ into a union of punctured tori. Then, for any $\ell>0$ small enough, the following inequality holds \[V_C(M(X))\leq K(g)d_{WP}(X, H_{\ell})+B(g, \ell)~,\] where $K$ and $B$ are constants depending, respectively, linearly on $g$ and on $(g, \ell, g\ell)$, and $H_{\ell}$ is the following subset of the Teichm\"uller space \[ H_{\ell}=\{X\in \mathcal{T}(\partial\bar M)\ |\  \exists m\in\Gamma_0^{comp}(\partial\bar M) \text{ s.t. } \ell_{\gamma}(X)\leq \ell\ \forall \gamma\in m\}~.\] 
	 \end{customcor}
	
	The renormalized volume can actually be defined and studied in the more general setting of \textit{geometrically finite} hyperbolic $3$-manifolds, as done in \cite{GGRVRfingeom}, \cite{VP2016}, and \cite{VP2017}. These works establish several continuity results for the renormalized volume, among which is the continuity under pinching an incompressible curve to obtain a geometrically finite limit.
	%that the renormalized volume is continuous under a sequence of convex co-compact acylindrical manifolds limiting to a geometrically finite one \cite[Theorem 1]{VP2016}, or, more generally, under (\textit{admissibly}) pinching an incompressible curve to obtain a rank-$1$ cusp (\cite[Theorem 4]{GGRVRfingeom}).
	
	The pinching of a compressible simple closed curve in the boundary limits, instead, in the pointed Gromov-Hausdorff topology, to a convex co-compact manifold $\textit{marked}$ at one or two points corresponding to the curve that has been pinched, as described in \cite[Appendix $A.10$]{averages}. The marked points are two when the disk $D$ compressing the pinched curve \textit{separates} $M$, that is $M\setminus D$ is disconnected. A point in the boundary of the Weil-Petersson completion of a closed surface $S$ is the data of a multicurve $m$ and a complete hyperbolic metric on $S\setminus m$, see Section \ref{boundarywpclosure}. Given $M$ a convex co-compact manifold, we will call \textit{compressible} a stratum in $\overline{\mathcal{T}(\partial\bar M)}^{\scriptscriptstyle{WP}}$ corresponding to a compressible multicurve.
	
	%In the Teichm\"uller space, a sequence of pinching a (multi)curve converges to a point in the strata of the same (multi)curve in the boundary of the Weil-Petersson completion (see the last subsection in Section \ref{norms}).
	
	A natural question to ask is whether the adapted renormalized volume, being bounded, converges, under the pinching of a compressible curve, to the sum of the renormalized volumes of the convex co-compact manifolds arising as Gromov-Hausdorff limits. The answer is negative. The key reason behind this lies in the fact that the pointed geometric limit of the convex cores does not coincide with the convex core of the limit manifold(s), but it is instead the convex hull of their marked point(s) at infinity union the convex core \cite[Lemma $A.8$]{averages}. In fact, the Gromov-Hausdorff limits retains memory of the (multi)curve that has been pinched, encoded through the marked point. Consequently, the limit of the renormalized volume depends on the marked points, in addition to the limit manifold.
	
	In this work then, after analyzing the behavior of the Epstein surface in a neighborhood of hyperbolic cusp singularities, we define the \textit{renormalized volume of pointed convex co-compact manifolds} (see Section \ref{SecRVM} and in particular Definition \ref{cuspedVR}).  	
	We then prove that, under pinching a compressible (multi)curve, the adapted renormalized volume converges to the sum of the adapted renormalized volumes of the pointed limit convex co-compact manifolds. In this sense, the adapted renormalized volume extends continuously to the strata of the boundary of the Teichm\"uller space corresponding to compressible multicurves.
	
	\begin{customthm}{\ref{AVRContinuity}}
		Let $M_t=(M, g_t)$ be a path of convex co-compact hyperbolic $3$-manifolds obtained by pinching a compressible multicurve $m$ in the conformal boundary at infinity. Let $D(m)$ be a union of disks compressing $m$, and let $(M_i, g_i, P_i)$, for $i=1,\dots,k$, be the pointed convex co-compact limits of $(M_t, y_{i}(t))$ in the Gromov-Hausdorff topology, with $ y_{i}(t)$ in the thick part of the $i$-th connected component of $C(M_t)\setminus D(m)$. Then \[\lim_{t\rightarrow\infty} \widetilde{V_R}(M_t)=\sum_{i=1}^k \widetilde{V_R}(M_i, g_i, P_i)~.\]
	\end{customthm}
	
	%The admissibility assumption on the multicurve ensures the Gromov-Hausdorff limits to be convex co-compact (see Theorem \ref{GromovHausd}, Remark \ref{Remadmissible} and Definition \ref{admissible}).
	
	\subsection{Outline of the paper} 
	In Section \ref{background}, we provide the main objects and preliminary knowledge needed. In Section \ref{RVM}, we define the renormalized volume associated to a convex co-compact manifold with marked points at infinity (Definition \ref{cuspedVR} and Remark \ref{P}). The main idea is to associate a renormalized volume to the unique hyperbolic representative conformal to the boundary at infinity with a cusp singularity at each marked point. To this aim, we first study the divergence of the $W$-volume of a truncated hyperbolic cusp (Proposition \ref{VRcusp1}) by explicitly constructing the associated Epstein surface (in Section \ref{SecI0}). 
	
	In Section \ref{RVA}, we start by defining the adapted renormalized volume and prove that this is uniformly bounded from below (Theorem \ref{AVRbounded}), that its differential has bounded infinity norm (Theorem \ref{boundedGradientinfty}), deducing the Weil-Petersson Lipschitz continuity in Corollary \ref{WPdistance}, and that its gradient has bounded Weil-Petersson norm (Corollary \ref{boundedGradientWP}). We also show the bounds for the distance from the convex core volume function (Theorem \ref{Vcdistance}), and the estimate in terms of Weil-Petersson distance for the convex core volume function of an handlebody (Corollary \ref{Vcestimate}). Moreover, we prove that the variation of the adapted renormalized volume is arbitrarily small on Bers regions of compressible pants decompositions (Theorem \ref{BersRegion}). We then proceed by defining the adapted version of the renormalized volume of convex co-compact manifolds with marked points at infinity (Definition \ref{AVRpointed}), and finally proving the continuity of $\widetilde{V_R}$ under pinching a compressible multicurve (Theorem \ref{conv}). This is done by exploiting Section \ref{VRtube}, in which we study the divergence of the $W$-volume of long hyperbolic tubes (Theorem \ref{Wtubegeneral}). 
	
	\subsection*{Acknowledgements} The author would like to thank her advisor Jean-Marc Schlenker for proposing this project and for all the essential discussions. She also thanks Francesco Bonsante for his insightful comments following a presentation on an earlier version of this work, and Kenneth Bromberg for suggesting the problem leading to Corollary $4.17$. Finally, she thanks David Fisac for useful conversations and for creating Figure $4$.
	The author was supported by FNR AFR grant 15777.

	\section{Preliminaries}\label{background}

	\subsection{Epstein surfaces and W-volumes}\label{SecEpst}
	%Epstein surfaces were already introduced in Section \ref{SecEpst}, however, this chapter requires a more in-depth exploration of the topic. 
	%	Epstein surfaces are the main ingredient of this paper. These had been introduced for the first time in \cite{epsteinSurf} (the original pre-print is from the year $1984$, but we reference a LaTex version transcribed by Bridgeman), and they produce a surface in $\mathbb{H}^3$ starting from a projective domain. When the domain under consideration is a discontinuity domain $\Omega(\Gamma)$, the associated Epstein surface is $\Gamma$-invariant, and it descends to a surface in the quotient $3$-manifold $\mathbb{H}^3/\Gamma$.
		Epstein surfaces are essential for defining the $W$-volume of a convex co-compact manifold $M$ with respect to a Riemannian metric $g$ that is conformal to $\partial_{\infty}M$. 
	These were introduced for the first time in \cite{epsteinSurf} (the original pre-print dates to $1984$, but we reference a LaTex transcription by Bridgeman of $2024$). A comprehensive and elegant treatment of Epstein surfaces, the Osgood-Stowe differential and $W$-volumes, along with their interrelations, is provided in \cite{BBOS}. 
	
	\subsubsection{Epstein surfaces: construction and properties}\label{SecEpst}
	The construction of Epstein produces a surface in $\mathbb{H}^{n+1}$ as an envelope of horospheres, starting from any domain $\Omega$ in $\mathbb{S}^n=\partial_{\infty}\mathbb{H}^{n+1}$ equipped with a conformal metric $g=e^{2\phi}\abs{dz}^2$: \[\text{Eps}_{(\Omega, \phi)}\colon \Omega\rightarrow\mathbb{H}^{n+1}~,\]
	identifying the notations $(\Omega, g)$ and $(\Omega, \phi)$.
	Here, we are interested in the case $n=2$, and mainly when the domain in question is a $\Gamma$-invariant subset of the discontinuity domain $\Omega(\Gamma)$ of some convex co-compact manifold $M=\mathbb{H}^3/\Gamma$. In this setting, the resulting Epstein surface is $\Gamma$-invariant, so that it descends to a surface in the quotient $3$-manifold $\mathbb{H}^3/\Gamma$.
	
		Briefly, the construction is as follows. Let us first fix some notations. Let $z=x+iy$ be the complex coordinate of $\mathbb{C}$, then we use the following: 
	\begin{align*}
		&dz^2=dz\otimes dz= dx^2-dy^2+i(dx\otimes dy + dy\otimes dx)\\
		&d\bar{z}^2=d\bar{z}\otimes d\bar{z}= dx^2-dy^2-i(dx\otimes dy + dy\otimes dx)\\
		&|dz|^2= S(dz\otimes d\bar z)=dx^2+ dy^2~.
	\end{align*}
	Given $\Omega\subseteq \mathbb{C}$ an open domain, a \textit{conformal metric} on $\Omega$ can be expressed as \[g=e^{2\phi}|dz|^2~,\] with \[\phi\colon \Omega\rightarrow \mathbb{R}\] a $C^k$ function, which, if the metric is required to be Riemannian, it must be smooth.
%	In the next section, we will sometimes consider conformal factors $\phi$ with \textit{cusp singularities} at a finite set of points $P\subseteq \Omega$. In this setting, and in light of Section \ref{ConfMarked}, it is possible to define cusp singularities for Riemann surfaces as follows.
%	\bdefi
%	Let $g_P=e^{2\phi}\abs{dz}^2$ be a conformal metric on a domain $\Omega\subseteq \mathbb{C}$, where $\phi$ is smooth (or just $C^k$) out of a finite set of points $P\subseteq\Omega$. Then $g_P$ has a \textit{cusp singularity} at $p\in P$ if, for a local coordinate such that $z(p)=0$, the function $\phi+\log{(\abs{z}\log{\abs{z}})}$ extends continuously to $p$.
%	\edefi
	Let then $g=e^{2\phi}\abs{dz}^2$ be a conformal metric on a domain $\Omega\subseteq\mathbb{CP}^1=\partial_{\infty}\mathbb{H}^3$. For each $z\in \Omega$, consider the \textit{horosphere} $H(z, e^{-\phi(z)})$ in $\mathbb{H}^3$ pointed at $z$ and of Euclidean radius $e^{-\phi(z)}$. The \textit{Epstein surface associated to $(\Omega, \phi)$} is then obtained by taking the boundary of the convex envelope of the union of this family of horospheres. In particular
	\[D_z\text{Eps}_{(\Omega, \phi)}\subset T_{\text{Eps}_{(\Omega, \phi)}(z)}H(z, e^{-\phi(z)})~,\] where $D_z$ denotes the image of the tangent map at $z$.
	 The horosphere $H(z, e^{-\phi(z)})$ coincides with the set of points in $\mathbb{H}^3$ whose \textit{visual metric} at $z$ coincides with $e^{2\phi}\abs{dz}^2$. The \textit{visual metric from $p$} is obtained by pulling back the spherical metric on $\partial_{\infty}\mathbb{H}^3=\mathbb{C}\cup\{\infty\}$ via any isometry of $\mathbb{H}^3$ sending $p$ to $(0,0,1)$.
	 The Epstein map is then natural in the sense that for any M\"obius transformation $f\in\mathbb{P}SL(2,\mathbb{C})=\text{Isom}^+(\mathbb{H}^3)$ \[\text{Eps}_{(\Omega, \phi)}=f\circ\text{Eps}_{(f^*\Omega, \phi\circ f+\log(\abs{df}))}~.\]
	The metric $g$ does not need to be smooth. If $g$ is $C^k$, then the Epstein map is $C^{k-1}$. On the other hand, even when g is smooth, the resulting surface may fail to be immersed. For example, the image of $\text{Eps}_{(\mathbb{CP}^1, \phi)}$, with $\phi$ such that $e^{2\phi}\abs{dz}^2$ is the spherical metric on $\mathbb{C}\cup \{\infty\}$, is the point $p=(0,0,1)$.
	 To avoid the singular cases, the strategy
	is to consider the Epstein surfaces associated to rescalings of the conformal metric on $\Omega$. More precisely, it is possible to define the map to the unit tangent bundle \[\widetilde{\text{Eps}}_{(\Omega, \phi)}\colon \Omega\rightarrow T^1\mathbb{H}^3~\] associating to $z\in\Omega$ the unit normal vector at $\text{Eps}_{(\Omega, \phi)}(z)$ whose geodesic ray ends at $z$, so that, if $\pi\colon T^1\mathbb{H}^3\rightarrow\mathbb{H}^3 $ is the bundle projection, then \[\text{Eps}_{(\Omega, \phi)}=\pi\circ\widetilde{\text{Eps}}_{(\Omega, \phi)}~.\]
	When the conformal metric at infinity is smooth, the map $\widetilde{\text{Eps}}_{(\Omega, \phi)}$ is a smooth immersion. 
	
	\bprop(\cite{epsteinSurf}, \cite[Section $5.1$, Theorem $5.8$]{BBOS}, \cite[Lemma $5.6$]{S2008}, or \cite[Theorem $7.4$]{BBOS})\label{Epstflow}
	Let $F_t\colon T^1\mathbb{H}^3\rightarrow T^1\mathbb{H}^3 $ be the geodesic flow, then \[F_t\circ\widetilde{ \text{Eps}}_{(\Omega, \phi)}= \widetilde{\text{Eps}}_{(\Omega, \phi+t)}~.\] If both the principal curvatures of the image of $\text{Eps}_{(\Omega, \phi)}$ are different from $-1$, then it is an immersion. Moreover, in this case, for $t$ big enough it is a locally convex immersion.
	\eprop
	
	Note that the Epstein surface image of $\text{Eps}_{(\Omega, \phi+t)}$ is a connected component of the set of $t$-equidistant points from the image of $\text{Eps}_{(\Omega, \phi)}$.
	
	\bthm(\cite{S2008}, \cite[Lemma $8.2$]{BBOS})\label{EpstEmb}
	If the conformal metric on $\Omega$ is smooth, and the associated Epstein surface is a locally convex immersion, then it is a convex embedding.
	\ethm

	\subsubsection{Epstein surfaces coordiantes}\label{EpstSurf}
%	Recall that to any smooth conformal metric $g=e^{2\phi}\abs{dz}^2$ on a domain $\Omega\subseteq\mathbb{CP}^1=\partial_{\infty}\mathbb{H}^3$ it is possible to associate a convex embedded surface $\Sigma(g)$ in $\mathbb{H}^3$ as the image of the Epstein map $\text{Eps}_{\Omega, \phi}\colon\Omega\rightarrow\mathbb{H}^3$ (see Section \ref{SecEpst}). Sometimes, to ensure the regularity of $\Sigma(g)$, it is necessary to rescale $g$ by a large enough constant (see Proposition \ref{Epstflow} and Theorem \ref{EpstEmb}). The Epstein surface $\Sigma(g)$ is given by the boundary of the convex envelope of the union of horospheres centered at points $z\in\Omega$ and of euclidean radii equal to $e^{-\phi(z)}$.
	% We summarize here the main points of the Epstein construction, and we report the main facts that will be needed, but we refer to \cite{epsteinSurf}, \cite[Section $3$]{DumasEpst}, \cite{BBOS}, and \cite{BrockPalleteWbound} for a complete treatment.
	%	Briefly, to any $z\in D$ it is associated a horosphere in $\mathbb{H}^3$ pointed at $z$ and of Euclidean radius $e^{-\phi(z)}$. The \textit{Epstein surface} associated to the Liouville field $\phi$ is then the boundary of the convex envelop of the union of all the horospheres. Let us denote it by $\Sigma(g)$. The surface $\Sigma(g)$ it is then parameterized by $D$, in the upper-half space model $\mathbb{C}\times\mathbb{R}^{+}$ with coordinates $(z,t)$, as the following: 
	Given a conformal metric $g=e^{2\phi}\abs{dz}^2$ on a domain $\Omega\subseteq \mathbb{CP}^1$, thanks to the aforementioned construction, it is possible to find the explicit coordinates of the Epstein surface $\Sigma(g)$ image of $\text{Eps}_{(\Omega, \phi)}$ in terms of $\phi$ and its first derivatives:
	\begin{equation}\label{Sigma}\Sigma(g)(z)= \left(z+\f{2e^{-2\phi}\nabla\phi}{1+e^{-2\phi}\abs{\nabla\phi}^2},\  \f{2e^{-\phi}}{1+e^{-2\phi}\abs{\nabla\phi}^2}\right)~,\end{equation}
	where $\nabla(\cdot)$ is the gradient with respect to the Euclidean metric on $\mathbb{C}$. In this notation, for any $r>0$, the Epstein surface $\Sigma(e^{2r}g)$ is $r$-equidistant from the one produced by the original metric $g$. 
	%Moreover, for $r$ big enough, the Epstein surface is convex immersed.

	\subsubsection{Fundamental forms}
%	As already seen in Section \ref{I*}, there is a beautiful correspondence between the fundamental forms on the Epstein surface induced by $\mathbb{H}^3$ and those of the domain at infinity $\Omega\subseteq \mathbb{CP}^1$ (see Definition \ref{DefI*}). We use here the same notations as in Section \ref{I*}. The key properties that will be used in this chapter are summarized below.
	%Denoting by $I^*$ the conformal metric at infinity, and by $I$, $B$, $\II=IB$ the first, second and third fundamental forms and the shape operator induced on the Epstein surface associated to $I^*$ by $\mathbb{H}^3$, these satisfy: 
	%\[I^*=\f12 (Id+B)^*I\]
	%Moreover, the \textit{second fundamental form and shape operator at infinity} $\II^*$ and $B^*$ can be defined as \[B^*=(Id+B)^{-1}(Id-B)~\]
	%and
	%\[\II^*=I^*B^*\]
	%so that the couple $(I, B)$ satisfies the Gauss-Codazzi equations in $\mathbb{H}^3$ if and only if $(I^*, B^*)$ satisfies the \textit{projective Gauss-Codazzi equations}.
	
	%he second fundamental form at infinity can be described in terms of the Liouville field in the following way. 
	
	Let us fix $\Omega\subseteq\mathbb{CP}^1$ a complex projective domain, and let us denote by $\Sigma(g)$ the Epstein surface associated to $(\Omega, \phi)$, with $g=e^{2\phi}\abs{dz}^2$. At any immersed point, the \textit{first fundamental form} $I$ of $\Sigma(g)$ is the pull-back via $\text{Eps}_{(\Omega, \phi)}$ to $\Omega$ of the Riemannian metric induced on $\Sigma(g)$ by $\mathbb{H}^3$. The \textit{shape operator} $B\colon T\Sigma(g)\rightarrow T\Sigma(g)$ is the bundle morphism defined as \[B(X)=\nabla_{X}(N)~,\] where $N$ is the outer unit normal to $\Sigma(g)$, and $\nabla$ is the Levi-Civita connection of $\mathbb{H}^3$. The \textit{second fundamental form} of $\Sigma(g)$ is \[\II(\cdot, \cdot)=I(B(\cdot), \cdot)~,\] and the \textit{third fundamental form} is \[\III(\cdot, \cdot)=I(B(\cdot), B(\cdot))~.\] We use the same notation for their pull-backs via $\text{Eps}_{(\Omega, \phi)}$ to $\Omega$. Similarly, we denote by $I_t$, $B_t$, $\II_t$ and $\III_t$, the pull-back via $\text{Eps}_{(\Omega, \phi+t)}$ on $\Omega$ of the fundamental forms and the shape operator of $\Sigma(e^{2t}g)$. 
	
	\bdefi\label{DefI*}(\cite{S2008}, \cite[Section $5$]{BBOS})
	In the notations above, the \textit{fundamental forms at infinity} of $(\Omega, \phi)$ are defined as \begin{align*}\hat I&= I+2\II+\III= I((Id+B)(\cdot),  (Id+B)(\cdot))\\
		\hat{\II} &= I-\III=I((Id+B)(\cdot),  (Id-B)(\cdot))\\
		\hat{\III}&= I-2\II+\III= I((Id-B)(\cdot),  (Id-B)(\cdot))\end{align*} and the \textit{shape operator at infinity} as \[\hat B= (Id+B)^{-1}(Id-B)~.\]
	\edefi
	 There is a beautiful correspondence between the fundamental forms on the Epstein surface induced by $\mathbb{H}^3$ and those of the domain at infinity.
	\bprop(\cite{epsteinSurf}, \cite[Lemma $5.7$, Theorem $5.8$]{S2008})
	In the notations above, let $I_t$ denote the pull-back on $\Omega$ of the induced metric on $\Sigma(e^{2t}g)$ via $\text{Eps}_{(\Omega, \phi+t)}$. Then \[I_t=\f14\left(e^{2t}\hat I+2\hat{\II}+e^{-2t}\hat\III\right)~,\] in particular \[\lim_{t\rightarrow \infty}4e^{-2t}I_t=\hat I~.\] Moreover \[\hat I = g~.\]
	\eprop
	
	Furthermore, the induced fundamental forms on a Epstein surface can be expressed in terms of the ones at infinity as in Definition \ref{DefI*}, as done in \cite[Lemma $5.6$]{S2008}.
	
	%\blem[\cite{S2008}, ]\label{fromII*}
	%The fundamental forms $I, \II, \III$ of an embedded Epstein surface are obtained from $\hat{I}, \hat{II}, \hat{\III}$ as follows 
	%\begin{align*}I&= \f14\left(\hat I+2\hat{\II}+\hat{\III}\right)= \f14\hat I((Id+\hat B)(\cdot),  (Id+\hat B)(\cdot))\\
	%	\II &= \f14\left(\hat I-\hat{\III}\right)=\f14\hat I((Id+\hat B)(\cdot),  (Id-\hat B)(\cdot))\\
	%	\III&= \f14\left(\hat I-2\hat{\II}+\hat{\III}\right)=\f14 \hat I((Id-\hat B)(\cdot),  (Id-\hat B)(\cdot))\end{align*} and the \textit{shape operator} as \[B= (Id+\hat B)^{-1}(Id-\hat B)~.\]
	%\elem
	
	\brem
	Note that any result on Epstein surfaces has an equivariant version. When $\Omega$ is the discontinuity domain of $M=\mathbb{H}^3/\Gamma$, then $\Sigma(g)$ is $\Gamma$-invariant and descends to a surface in $M$. We will keep denoting by $I$, $B$, $\II$ and $\III$, respectively, the Riemannian metric, shape operator, second and third fundamental forms induced on $\Sigma(g)/\Gamma$ by $M$.
	\erem

	Additional key properties that will be used later are summarized below.
	\blem\label{secondfun}(\cite{TZUnif, S2008})
	Let $\hat I$ be a Riemannian metric on $\Omega\subseteq\mathbb{CP}^1$, and let $\phi\colon \Omega\rightarrow \mathbb{R}$ be a smooth function such that $\hat I=e^{2\phi}\abs{dz}^2$. Then, the second fundamental form $\hat{\II}$ can be written in terms of $\phi$ as \[\hat{\II}=2 q dz^2 + 2\bar q d\bar{z}^2+4\phi_{z\bar{z}}|dz|^2\] with 
	\[q= \phi_{zz}- (\phi_z)^2~.\]
	Moreover, the curvature $K$ of $\hat I$ satisfies \[K=-4\phi_{z\bar z}e^{-2\phi}~.\]
	\elem
	
	\brem\label{secondfun2}
	In the coordinates $(x,y)$ such that $z=x+iy$, the equation in Lemma \ref{secondfun} becomes \[\hat{\II}= (4\phi_{z\bar{z}}+4\re(q))dx^2+(4\phi_{z\bar{z}}-4\re(q))dy^2-4\im(q)dxdy\] 
	where $dxdy$ denotes twice the symmetric tensor product between $dx$ and $dy$.
	\erem
	
	\blem\label{meancurv}(\cite{S2008})
	The mean curvature $H= \tr(B)/2$ of an Epstein surface, with $B$ the shape operator, can be expressed in terms of the data at infinity as \[H=\f{1-\det(\hat B)}{1+\tr(\hat{B})+\det(\hat{B})}~.\]
	\elem

	\subsubsection{Epstein surfaces of domains with boundary} \label{BoundaryEps}
	In \cite{BrockPalleteWbound}, Brock and Pallete extended the definition of Epstein surfaces associated to a $k$-dimensional open domain in the boundary at infinity $\partial_{\infty}\mathbb{H}^n$ of $\mathbb{H}^n$ (as introduced in Section \ref{SecEpst}) to domains with boundary. Here, we are interested in the dimensions $k=n=2$. We report the main facts we will need, restricting to the case of round boundary components, which is also the same setting where it is possible do define an associated $W$-volume, again introduced in \cite{BrockPalleteWbound}, and which we discuss in the next subsection. 
	
	We now briefly outline the construction in \cite{BrockPalleteWbound} (see, in particular, Definition $2.1$) and fix some notations. Let $\Omega\subseteq \mathbb{CP}^1$ be a domain with boundary $\partial \Omega$ given by a union of $k$ round circles, and let $e^{2\phi}|dz|^2$ be a $C^{2,\alpha}$ conformal metric on $\Omega$. It is then possible to associate to the interior of ${\Omega}$ its Epstein surface $\Sigma(\Omega, \phi )$ in $\mathbb{H}^3$, as in Section \ref{SecEpst}, and whose coordinates are furnished by Equation (\ref{Sigma}). Moreover, it can be defined a surface $\Sigma(\partial \Omega, \phi )$ parameterized by the normal bundle of $\partial \Omega$ and such that at each point $p\in\partial \Omega$ is assigned the unique horocycle in the horosphere centered at $p$ and Euclidean radius $e^{-\phi(p)}$, that is also tangent to the Epstein surface $\Sigma(\Omega, \phi)$.
	
	We can then construct a compact region $N(\Omega, \phi )$ in $\mathbb{H}^3$ as follows: consider the union $H$ of hyperplanes with boundary at infinity in $\partial \Omega$, and connect then $\Sigma(\Omega, \phi )$ to $H$ via $\Sigma(\partial \Omega, \phi )$, which intersects $\Sigma(\Omega, \phi )$ tangentially and meets $H$ orthogonally. This, by gluing, forms a sphere with piece-wise smooth boundary, which can be filled to obtain the compact region $N(\Omega, \phi)$. Each of the $k$ connected components of the region connecting the Epstein surface $\Sigma(\Omega, \phi)$ to the hyperplanes in $H$ is called a \textit{caterpillar region}, and we will denote their union by $C=\cup_{i=1}^{k}C_i$. Let us remark that $C$ can be parameterized by $\partial \Omega \times I$, for $I$ some compact interval (see \cite[Section $2.2$]{BrockPalleteWbound}): \begin{equation}\label{Pcaterpillar}C(s,v)=\end{equation} \[\left(\gamma(s)+\gamma'(s)\f{2\phi'(s)}{e^{2\phi(s)}+\abs{\phi'(s)}^2+v^2}+i\gamma'(s)\f{2v}{e^{2\phi(s)}+\abs{\phi'(s)}^2+v^2}, \f{2e^{\phi(s)}}{e^{2\phi(s)}+\abs{\phi'(s)}^2+v^2} \right)~,\]
	for $\gamma(s)$ a unit velocity parameterization of $\partial\Omega$ with respect to $\abs{dz}^2$ . 
	We will refer to the union of the Epstein surface of the interior and the caterpillar region as the \textit{Epstein surface of the domain with boundary}.
	
	Whenever we have a metric $g=e^{2\phi}|dz|^2$, we will often denote the pair $(\Omega, \phi)$ as $(\Omega, g)$.
	
	\subsubsection{W-volume and Renormalized volume}\label{WvolVR}
	We can now finally define $W$-volumes and the renormalized volume of convex co-compact hyperbolic $3$-manifolds. We start by define the $W$-volume of a compact region.
	\bdefi\label{DefiWvol}
	Let $M$ be convex co-compact and $C\subseteq M$ be a compact subset with (piece-wise) smooth boundary. The {\em $W$-Volume of $C$} is defined as 
	\[W(C)=\text{Vol}(C)-\dfrac{1}{2}\int_{\partial C} H dA_{\partial C}~,\]
	where $\text{Vol}(C)$ is the hyperbolic volume of $C$ with respect to the metric of $M$, $H=\text{Tr}_I(\II)/2$ denotes the \textit{mean curvature} of $\partial C$, and $dA_{\partial C}$ is the area form of the induced metric on the boundary $\partial C$.
	\edefi
	Note that we can also take $M=\mathbb{H}^3$ in the definition.
When the boundary is piece-wise smooth, the integral of the mean curvature term decomposes into the sum of the integrals over the smooth faces plus the contributions from the bending lines $b_j$ of their intersections. In the case the exterior angles $\theta_j$ at each $b_j$ are constant, these terms are given by \begin{equation}\label{b_j}-\f14\sum_{j=1}^k \theta_j\ell(b_j)~,\end{equation}
with $\ell(b_j)$ the induced length of $b_j$. A special case is when $C$ is the \textit{convex core of $M$} \[C(M)=\text{Hull}(\Lambda(\Gamma))/\Gamma~,\] where $\text{Hull}(\cdot)$ denotes the convex envelope inside $\mathbb{H}^3$, and $\Gamma<\mathbb{P}SL(2, \mathbb{C})$ the subgroup s.t. $M=\mathbb{H}^3/\Gamma$. The integral of the mean curvature on its boundary becomes one-half the \textit{length of the bending lamination} $L(\beta)/2$ (see \cite{thurston-notes,epstein-marden}), so that
\begin{equation}\label{WCC}
	W(C(M))=\text{Vol}(C(M))-\f{L(\beta)}{4}~.
\end{equation}

Given now a smooth metric $g\in[\partial_{\infty}M]$, up to rescaling $g$ by a constant $e^{2r}$, the Epstein construction furnishes an embedded convex surface in $M$ homeomorphic to $\partial\bar M$ (see Section \ref{SecEpst} and Theorem \ref{EpstEmb}), which therefore bound a \textit{strongly geodesically convex} subset $C_r(g)$ of $M$, that is, a convex subset homotopically equivalent to $M$. By \cite[Lemma $4.2$]{S2008}, for every $t>0$ \[W(C_r(g))=W(C_{r+t}(g))+\pi t\chi(\partial\bar M)~.\] Thank to this, it is possible to define the $W$-volume of any smooth metric $g$ in $[\partial_{\infty}M]$ as follows. 

\bdefi
Let $M$ be a convex co-compact manifold, and let $g\in[\partial_{\infty}M]$ be a smooth conformal metric. Let also $r>0$ be large enough such that the Epstein surface of $e^{2r}g$ bounds a strongly geodesically convex subset $C_r$ of $M$. Then, the \textit{W-volume of $M$ with respect to $g$} is defined as \[W(M,g)=W(C_r)+\pi r \chi(\partial \bar M)\]
\edefi

We are now ready to define the renormalized volume of $M$.

\bdefi\label{DefiVR}
For any convex co-compact hyperbolic $3$-manifold $M$, the {\em renormalized volume} is defined as 
\[V_R(M)=W(M,h)~,\] 
with $h\in [\partial_{\infty}M]$ the hyperbolic representative in the conformal boundary at infinity.
\edefi
Thanks to the parameterization of the deformation space of convex co-compact structures $CC(M)$ provided by Theorem \ref{UnifCCM}, the renormalized volume can be interpreted as a function on the Teichm\"uller space of $\partial\bar M$	\[V_R\colon\mathcal{T}(\partial \bar{M})\longrightarrow \mathbb{R}~.\]

	\subsubsection{W-volume for domains with round boundary}
%	During this chapter, we will often talk about $W$-volumes of compact regions, so we recall the definition here.
	
%	\bdefi\label{Wvol}
%	Let $M$ be a convex co-compact manifold, and let $N\subseteq M$ be a compact subset of $M$ with smooth (or piece-wise smooth) boundary. The \textit{W-volume of $N$}  is defined as \[W(N)=\text{Vol}(N)-\f12 \int_{\partial N} H da_{\partial N}~,\]
%	where $\text{Vol}(N)$ is the hyperbolic volume of $N$ induced by $M$, $H$ is the mean curvature on the boundary $\partial N$ of $N$, and $da_{\partial N}$ is the induced area form on $\partial N$.
%	\edefi

	%	\bdefi\label{W}
	%	Let $g$ be a metric conformal to the boundary at infinity $\partial_{\infty}M$ of a convex co-compact hyperbolic $3$-manifold. Let then $\widetilde{\Sigma}(g)$ the $\gamma$-invariant Epstein surface in $\mathbb{H}^3$ associated to $(\Omega(\Gamma), g)$, and assume that this is convex immersed. Let then also $\Sigma(g)$ be its quotient inside $M$, and let $N(g)\subseteq M$ be the compact detected by $\Sigma(g)$. We define the \textit{$W$-volume of $g$} as \[W(M,g)=W(N(g))~.\] More, we define the \textit{renormalized volume of M} as the $W$-volume of the unique hyperbolic representative $h$ conformal to $\partial_{\infty}M$ \[V_R(M)=W(M,h).\]
	%	\edefi

	\bdefi\label{Wboundary}
	In the notations above, the \textit{$W$-volume associated to a domain with $k$ round boundary components} and equipped with a conformal metric $(\Omega, g)$, with $g=e^{2\phi}|dz|^2$, is defined as \[W(\Omega, g)= W(N(\Omega, \phi))=Vol(N(\Omega, \phi ))-\f12\int_{ S} Hda_{S}~,\]
	where $S=\partial N(\Omega, \phi )$, $H$ is the mean curvature of $S$, and $da_S$ denotes the area form induced on the boundary of $N(\Omega, \phi )$ (as in Section \ref{BoundaryEps}) by $\mathbb{H}^3$.
	\edefi
	
	Since $S$ is piece-wise smooth with constant exterior angles along its co-dimension one faces, and since $H=0$ on a hyperplane, the integral of the mean curvature in the definition splits into the sum of the integral over the Epstein surface $\Sigma(\Omega, \phi)$ of the interior of $\Omega$, the ones over the caterpillar regions $C_j\subseteq \Sigma(\partial \Omega, \phi)$, and the ones over the faces $b_j$ as in Equation \eqref{b_j}. We refer to the last two types of contributions as those arising from the $j$-th boundary component of the domain $\Omega$.
	
	\brem
	When $\Omega=\Omega(\Gamma)$ is a discontinuity domain of a convex co-compact manifold $M=\mathbb{H}^3/\Gamma$, the $W$-volume of the quotient in Definition \ref{Wboundary}  coincides with the one in Definition \ref{DefiWvol}, as the boundary term on the Caterpillar regions $C_j\subseteq \Sigma(\Omega, \phi)$ cancel out when paired by the action of $\Gamma$. 
	\erem
	
	It will be essential for our purposes to understand how the $W$-volume of a domain with round boundary components changes under a conformal change of the metric. To this aim, we state the following version of the \textit{Polyakov formula} (see \cite{AP}, \cite{BrockPalleteWbound}).
	
	\bthm[\cite{BrockPalleteWbound}, Theorem $2.1$]\label{Poly}
	Let $g=e^{2\phi}\abs{dz}^2$ be a smooth conformal metric on a domain $\Omega\subseteq \mathbb{C}$ with round boundary components, and let $u\colon \Omega\rightarrow \mathbb{R}$ be another smooth function, then \[W(\Omega, e^{2u}g)-W(\Omega, g)=-\f14\int_{\Omega}\left(\abs{\nabla_g u}^2+K(g)u\right)da(g)-\f12\int_{\partial \Omega}k(g) u ds(g)~,\]
	where $\nabla_g(\cdot)$ is the gradient with respect to $g$, $K(\cdot)$ and $k(\cdot)$ are, respectively, the scalar curvature on $\Omega$ and the induced geodesic curvature of the boundary $\partial \Omega$, and $da(\cdot)$ and $ds(\cdot)$ denote the area and length forms. 
	\ethm 
	
	Requiring $\phi$ and $u$ to be of class $C^{2,\alpha}$ is actually sufficient for the equality in the theorem to hold.
	
	\brem\label{lastterm}
	In the original definition of $W$-volume for domain with round boundary components (see Definition $2.2$ in \cite{BrockPalleteWbound}), the following additional integral on the caterpillar region $C$ appears: \[-\f32 \int_{\partial \Omega\times [0,1]}(1+H)da_C~,\]
	whose variation is \[-\f34\int_{\partial \Omega}\partial_n u ds(g)~,\]
	where $\partial_{n}(\cdot)$ denotes the derivative with respect to the outer normal to the boundary $\partial \Omega$.
	In this formulation, the relative Polyakov formula stated as Theorem $2.1$ in \cite{BrockPalleteWbound} coincides with the variation of the \textit{determinant of the Laplacian} for domain with boundary (see \cite{AP}). The three terms in Theorem \ref{Poly}, instead, arise from an application of the Schl\"afli formula for the variation of the volume of compact regions with piece-wise smooth boundary (see \cite{SchlafliPiesewise}, \cite{RivinJM}, and Section $4$ in \cite{KSJM} in which the Schl\"afli formula is expressed in terms of the fundamental forms at infinity).  
	Here, we make use of the definition of $W$-volume for domains with round boundary components to cut a domain equipped with a complete hyperbolic metric. What is nice about this construction, is that the caterpillar regions ensure the intersection between the Epstein surface of the domain and the hyperplanes to always have angle $\pi/2$, so that we do not have to account for its variation.
	\erem
	
	In the next sections, it will be crucial to have some form of additivity property for the $W$-volume. 
	
	\blem\label{add}
	Let $(\Omega, g)$ be a domain in $\mathbb{C}$ with round boundary components equipped with a conformal metric $g$. Suppose $\Omega=\Omega_1\cup \Omega_2$ to be a decomposition of $\Omega$ such that $\Omega_1\cap\Omega_2$ consists of a finite union of round circles $\gamma_i$. Then \[W(\Omega, g)=W(\Omega_1, g)+W(\Omega_2, g)~.\]
	
	\elem
	
	\bpf 
	We follow the notations introduced in this section.
	The key observation is that the caterpillar regions of $\partial\Omega_1$ and $\partial \Omega_2$ on the common components $\gamma_i$ coincide, with reversed outer normals, as also the hyperplanes with boundary in $\cup\gamma_i$. Then, the terms of the caterpillar regions associated to $\gamma_i$ simplify as they are opposite. Moreover, the sum of the hyperbolic volumes of $N(\Omega_1, g)$ and $N(\Omega_2, g)$ is equal to the hyperbolic volume of $N(\Omega, g)$, and the same holds for the integrals of the mean curvature on the Epstein surfaces $\Sigma(\Omega_i, g)$ of the interiors (the exterior angles are considered modulo $2\pi$ or $\pi$, depending on whether they lie in the interior or on the boundary of $N(\Omega, g)$).
	\epf
	
	\subsection{Norms on quadratic differential and the Weil-Petersson pairing}\label{norms}
	Let $X$ be a Riemann surface structure on a closed surface $S$ of genus $g\geq 2$, and let $z=x+iy$ and $\bar{z}=x-iy$ denote its local holomorphic and anti-holomorphic coordinates.  A \textit{quadratic differential} $\phi$ is a symmetric $(0,2)$-tensor that can be locally written as \[\phi(z)=q(z)dz\otimes dz~.\] When $q(z)$ is a holomorphic function, $\phi$ is called a \textit{holomorphic quadratic differential}. The space of holomorphic quadratic differential on $X$ is denoted by $Q(X)$. It is customary to use the notations \[dz\otimes dz=  dz^2~,\] and \[\abs{dz^2}=-\f{1}{2i}dz\wedge d\bar{z}= dx\wedge dy~.\] 
	 A \textit{Beltrami differential} on $X$ is a $(1,-1)$-tensor $\mu$ which can be expressed in local coordinate as 
	\[\mu(z)= \eta(z)\partial z\otimes d\bar{z}~,\]
	where $\eta$ is a measurable complex-valued function, and such that \[\norm{\mu}_{\infty}=\text{ess}\sup_{z\in X} \abs{\eta(z)}< \infty~,\]
	where the essential supremum means that the set of points in which the norm of $\eta$ is unbounded has measure zero.
	
	\brem\label{touniversal}
	Given a Riemann surface $X=\mathbb{H}^2/\Gamma$, the space $Q(X)$ can be identified with the space of $\Gamma$-invariant holomorphic quadratic differentials on the universal cover $\mathbb{H}^2$. The same holds for Beltrami differentials.
	\erem
	
	For any $\mu=\eta \partial_z\otimes d\bar{z}$ and $\phi=qdz^2$, there is a natural pairing defined as \begin{equation}\label{wppairing}\langle \phi, \mu\rangle = \int_{X} (q\eta) dx\wedge dy~,\end{equation}
	with $z=x+iy$ the local complex coordinate of $X$. This is called the \textit{Weil-Petersson pairing}, and
	establishes a duality between the two spaces.
	The pairing can be restricted to be between the space of \textit{holomorphic} quadratic differentials $Q(X)$, and the one of \textit{harmonic} Beltrami differentials $B(X)$, i.e. the quotient of the space of Beltrami differentials by the subspace of elements whose pairing with every holomorphic quadratic differential is null. 
	
	 Let $\rho_X=\rho\abs{dz}^2$ be the hyperbolic metric conformal to $X$, with $z$ a local conformal coordinate of $X$. Then, given $\phi\in Q(X)$, the quotient \[\norm{\phi(z)}=\abs{\phi(z)}/\rho_X(z)~,\] is a function which defines $L^p$ norms by integration on $X$ with respect to $\rho_X$.
	 In particular, the squared $L^2$ norm on holomorphic quadratic differentials at $X$ is \[\norm{\phi}_2^2=\int_X \f{\abs{\phi}^2}{\rho_X}~,\] and it is induced by the inner product  \[\re\langle \phi, \psi\rangle =\re\int_X \f{\phi\bar\psi}{\rho_X}~,\]
	 which furnishes a Riemannian metric on $T^*\mathcal{T}(S)$, and then, by duality, also on $T\mathcal{T}(S)$ (see \cite[Sections $6.6$ and $6.7$]{Hubbard2016}, or \cite[Chapter $6$]{QTT}) as \[\re\langle\mu, \nu\rangle= \int_X \mu\bar\nu\rho_X~.\] This is called the \textit{Weil-Petersson metric}.

	\subsection{The boundary of the Weil-Petersson completion}\label{boundarywpclosure}

	With respect to the Weil-Petersson metric, the Teichm\"uller space is not complete: there exist paths \textit{pinching} a simple closed curve $\gamma$ - that is, paths along which the hyperbolic length of $\gamma$ goes to zero - that have finite Weil-Petersson length and exit every compact subset (see, for example, \cite[Theorem $3$]{Wol06:Convexity}). This leads to define the \textit{augmented Teichm\"uller space}, \cite{BersWPNode}, \cite{AugmentedTeich}, \cite[Introduction]{WolpertGeometryWP}. Briefly, as a set, the augmented Teichm\"uller space of $S$ is the union \[\mathcal{T}(S)\bigcup_{m\in\Gamma(S)} \mathcal{S}(m)\] with \[\mathcal{S}(m)= \{X\ |\ X\text{ is a marked complete finite area hyperbolic structure on $S\setminus m$}\}~,\]
	where $\Gamma(S)$ denotes the set of isotopy classes of multicurves on $S$ (or, equivalently, the set of \textit{$0$-skeletons of simplices in the curve complex} of $S$). The sets $\mathcal{S}(m)$ are called \textit{strata}, and a point in $\mathcal{S}(m)$ can be interpreted as a \textit{nodal surface}, or, equivalently, as a hyperbolic surface where all the curve of $m$ have been \textit{pinched}: \[\mathcal{S}(m)=\{X\ |\ \ell_{\gamma}(X)=0\ \text{\ iff }\ \gamma\in m\}~.\]
	% The topology on the augmented Teichm\"uller space is furnished by extending the Fenchel-Nielsen coordinates in the following way. Given a multicurve $m$, complete it to a pants decomposition $P$. Then, the set $\mathcal{T}(S)\cup \mathcal{S}(m)$ is described by allowing the length parameters of the simple closed curves in $m\subseteq P$ to take values in $\mathbb{R}_{\geq 0}$. The union $\mathcal{T}(S)\cup \mathcal{S}(m)$ is equipped with the coarsest topology which makes the Fenchel-Nielsen coordinates map, forgetting about the twisting parameters of the curves in $m$, continuous \[\mathcal{T}(S)\cup \mathcal{S}(m)\longrightarrow \left(\mathbb{R}_{>0}\times \mathbb{R}\right)^{3g-3-|m|}\times \mathbb{R}_{\geq 0}^{|m|}~.\] It is possible to show that this topology does not depend on the choice of the pants decomposition $P$. Repeating the procedure for each $m\in\mathcal{C}(S)$ defines the topology on the whole space. 
	
	In \cite{BersWPNode} and \cite{Masur}, it is shown that the Weil-Petersson completion of the Teichm\"uller space coincides with the augmented Teichm\"uller space. Thus, the boundary of the Weil-Petersson completion of the Teichm\"uller space $\partial\overline{\mathcal{T}(S)}^{\scriptscriptstyle{WP}}$ is stratified by the space of multicurves $\Gamma(S)$: a point in the boundary is the data of a multicurve $m$ on $S$ and a complete finite area hyperbolic metric on $S\setminus m$. A \textit{stratum} in $\partial\overline{\mathcal{T}(S)}^{\scriptscriptstyle{WP}}$ is then the product of lower dimensional Teichm\"uller spaces of the connected components of $S\setminus m$. Moreover, the induced metric on $\partial\overline{\mathcal{T}(S)}^{\scriptscriptstyle{WP}}$ coincides with the (product of) Weil-Petersson metrics on the strata.
	
	\bdefi\label{pinchingacurve}
	We say that a path of Riemann surfaces $X_t$ in the Teichm\"uller space of $S$ is obtained by \textit{pinching} a multicurve $m\subseteq S$, if the hyperbolic lengths of all the components of $m$ tend to zero, and $X_t$ converges to a point in the Weil-Petersson completion of $\mathcal{T}(S)$.
	\edefi

	\subsection{Conformally equivalent surfaces with marked points}\label{ConfMarked} 
	For any pair $(X, D)$ with $S$ a closed Riemann surface and $D\subseteq S$ a finite set of points, by Uniformization Theorem, there exists a unique hyperbolic metric which is conformal to $X\setminus D$ and has cusps at each point in $P$, \cite[Chapter $10$, parabolic case]{ConfInv}, \cite{Heins}. Moreover, the behaviour of this conformal metric in a neighborhood of $p\in D$ can be described explicitly, and we are particularly interested in the local expression of the metric tensor. To this end, we report the version from \cite{MR3929543}, whose proof techniques are well suited to our context.
	
	\bthm (\cite[Theorem $1.2$]{MR3929543})\label{cuspshape}
	Let $\mathbb{D}^*$ be the punctured disk in $\mathbb{C}$ of radius one centered at the origin, and let $d\rho^2$ be a hyperbolic conformal metric on $\mathbb{D}^*$ with a cusp singularity at $0$. Then, there exists a complex coordinate $w$ on $\mathbb{D}(\varepsilon)=\{z\in\mathbb{C}\ |\ \abs{z}<\varepsilon\}$ for some $\varepsilon>0$ and with $w(0)=0$, such that \[d\rho^2=\f{1}{|w|^2\log^2(|w|)}dw^2~,\]
	Moreover, the coordinate $w$ is unique up to rotation.
	\ethm

	Note that Theorem \ref{cuspshape} does not assume the hyperbolic metric $d\rho^2$ on $\mathbb{D}^*$ to be complete, in which case it would not be necessary to restrict to a smaller disk or to change coordinate (see beginning of Section \ref{SecI0}).
	We remark that, up to composing with a M\"obius transformation, it is always possible to assume the cusp singularity to be at the origin. On the other hand, the change of coordinate from the standard one $z$ to $w$ is conformal but not necessarily projective, that is, it may not be realized by a M\"obius transformation. Since we will work with Epstein surfaces, which depend crucially on the complex projective structure (see Section \ref{SecEpst}), it is important to understand how far the aforementioned change of coordinate is from being projective. Examining the proof of Theorem $1.2$ in \cite{MR3929543}, we find that the answer is `not much', provided we restrict to a sufficiently small neighborhood of the cusp singularity.
	There is an equivalent analytic definition for a Riemann surface to have a cusp: the Schwarzian of the uniformization map is of the form  \[\f{1}{2z^2}+\f{d}{z}+\psi(z)~,\] where $d$ is a number and $\psi(z)$ a holomorphic function,  both depending on the choice of coordinate (see \cite[Lemma $2.1$]{CuspAnalyticDef}). Thanks to this, as done in \cite{MR3929543}, by integrating the Schwarzian, up to M\"obius transformation, one obtains that the developing map has the form \[g(z)=\log(z)+\psi_1(z)~,\] with $\psi_1(z)$ a holomorphic function such that $\psi_1(0)=0$. Then, in the coordinate $w$ such that $\log(w)=\log(z)+\psi_1(z)$, one gets \[g(w)=\log(w)~.\]
	Note that we can multiply $g$ by $-i$, which, being a M\"obius transformation, does not change the Schwarzian derivative, to get the developing map \[g_1(w)=-i\log(w)~,\]
	through which the pull back of the hyperbolic metric of $\mathbb{H}^2$ gives the one in Theorem \ref{cuspshape}. The only non projective change of coordinate needed is then \begin{equation}\label{changeofvar}w=ze^{\psi_1(z)}~,\end{equation}
	which, since $\psi_1(z)$ is holomorphic and such that $\psi_1(0)=0$, behaves at the first order at $0$ again like $z$. 
	
	\subsection{Thin Tubes}\label{margulistube}
	An important role in this work is played by \textit{long tubes}, i.e., conformal tubes with \textit{large modulus}. We now introduce a particularly important class of long tubes: those furnished by the Collar Lemma, \cite[Theorem $4.1.1$]{Bu1992}. A crucial property of the tubes in this class is that they are pairwise disjoint when their \textit{core} curves are sufficiently \textit{short}, \cite[Theorem $4.1.6$]{Bu1992}. This leads to the \textit{thin-thick decomposition} of hyperbolic surfaces, which in fact is a more general result holding for hyperbolic manifolds (see \cite[Chapter $4$]{Mar16}). We first fix the notation for the two dimensional \textit{Margulis constant} \begin{equation}\label{e0}\varepsilon_0=2\argsinh(1)~,\end{equation} which will appear frequently in the next chapters.  
	
	\bdefi \label{margulitubedefi}
	A \emph{thin tube} in a hyperbolic surface X, is the set of points $\mathcal A(\ell)$ around a simple closed geodesic $\gamma$ of length $\ell\leq \epsilon_0$ that are at a distance at most \[L\eqdef\argsinh\left(\f1{\sinh\left(\f\ell 2\right)}\right)~.\]
	\edefi
	
	We remark that there is no need to assume the length $\ell$ to be shorter than $\varepsilon_0$ in order for a collar of width $L$ as above to exist. However, the condition is both necessary and sufficient to ensure them to be disjoint, \cite[Theorem $4.1.6$]{Bu1992}.
	
	 The hyperbolic metric on the thin tube $\mathcal A(\ell)$ of core length $\ell$ can be written as 
	\[d\rho^2+\left(\f{\ell}{2\pi}\right)^2\cosh^2(\rho)d\theta^2~,\] 
	in the coordinates $\theta\in[0,2\pi]$ and $\rho\in[-L,L]$, see \cite[Collar Lemma: Theorem $4.1.1$]{Bu1992}.
	
	\subsection{Geometric convergence}\label{Geomconv}
	There are several equivalent definitions of geometric convergence for \textit{Kleinian groups}, i.e., discrete subgroups of $\mathbb{P}SL(2, \mathbb{C})$, and for their corresponding hyperbolic quotient $3$-manifolds (see for example \cite[Section $2.2$]{Mc1996}, \cite[Chapter $4$]{Mar16}, \cite[Chapter $7$]{MT1998}). Here, we present the one we will use.
	
	\bdefi
	A sequence of pointed hyperbolic $3$-manifold $(M_n, x_n)$ \textit{geometrically converges} to a pointed hyperbolic $3$-manifold $(M, x)$ if and only if, for every compact subset $K\subseteq M$ containing $x$, and for all sufficiently large $n$, there exists a smooth embedding $f_n\colon K\rightarrow M_n$ such that $f_n(x)=x_n$ and $f_n$ converges in the $C^{\infty}$ topology to an isometric embedding. The induced topology is also called the \textit{pointed Gromov-Hausdorff topology}.
	\edefi

	\section{Renormalized volume for pointed manifolds}\label{RVM}
	In the next section, we will define the adapted renormalized volume. In particular, we will show that this admits limit when approaching points on the boundary of the Weil-Petersson completion of the Teichm\"uller space (see Section \ref{boundarywpclosure}) in the strata corresponding to compressible multicurves. Additionally, we aim to provide a geometric interpretation of the limit quantity in such a way that the adapted renormalized volume extends continuously to such boundary points. The Gromov-Hausdorff limit (see Section \ref{Geomconv}) of a sequence of convex co-compact hyperbolic $3$-manifolds obtained by pinching a compressible multicurve is a finite union of new \textit{pointed convex co-compact manifolds}, that is, a disjoint union of pairs $(M_i, P_i)$, where $M_i$ is convex co-compact and $P_i\subseteq \partial_{\infty}M_i$ is a finite set of points (see \cite[Appendix $A.10$]{averages}). The goal of this section is then to define the renormalized volume for these objects. The key idea in constructing this extension is to consider the Epstein surface associated to the unique hyperbolic metric conformal to $\partial_{\infty}M_i$ with a cusp singularity at each $p_i\in P_i$. We begin by analyzing the case where $P_i$ has cardinality one, and then observe that the definition naturally extends to the general case.
	
	\subsection{Renormalized volume behaviour at cusp singularities}\label{cuspsing}
	In what follows, we analyze the divergence of the $W$-volume associated to a domain \[(\Omega\setminus\{p\}, h_p)\subseteq \mathbb{C}\setminus \{p\}\subseteq\mathbb{CP}^1\] equipped with a conformal hyperbolic metric with a cusp at $p$ (as in Section \ref{ConfMarked}), in a neighborhood of the singularity. The divergence is caused by the fact that at the cusp singularity the Epstein surface associated to $(\Omega\setminus \{p\}, h_p)$ touches the boundary at infinity $\partial_{\infty}\mathbb{H}^3$, definitively exiting any $r$-neighborhood of the geodesic of $\mathbb{H}^3$ through $p$ and $\infty$ (in the upper half-space model). We start by studying a toy model case: the \textit{hyperbolic cusp on a round punctured disk}. In this case, the metric $h_p$ is explicit, and, up to composing by a M\"obius transformation, we can assume $p=0\in\mathbb{C}$ and $\Omega=\mathbb{D}^*$, where $\mathbb{D}^*=\mathbb{D}\setminus \{0\}$ with $\mathbb{D}$ the unit disk centered at the origin. This will allow us to explicitly compute the $W$-volume of any compact sub-domain of $\mathbb{D}^*$ with two concentric round boundary components (as in Definition \ref{Wboundary}), and see how this diverges under shrinking one boundary component to $0$. Then, by applying the Polyakov formula of Theorem \ref{Poly} together with the argument of Section \ref{ConfMarked}, we show that the result obtained still holds for a general cusp singularity. This furnishes the right renormalization process needed in order to define the $W$-volume of a pointed convex co-compact hyperbolic $3$-manifold, which will be done in the next section.

	\subsubsection{Epstein surface of the hyperbolic cusp on a round punctured disk}\label{SecI0}
	
	We fix the notation $\mathbb{D}\subseteq\mathbb{C}\subseteq\mathbb{CP}^1$ for the subset of complex numbers $z\in\mathbb{C}$ of norm $\sqrt{z\bar{z}}$ less than one, and $|dz|^2=dx^2+dy^2$, with $z=x+iy$, for the Euclidean flat metric on $\mathbb{C}$. We will use the polar coordinates $(\rho, \theta)\in \mathbb{R}_{\geq 0}\times [0,2\pi]$ such that $z=\rho e^{i\theta}$. 
	Let us consider the metric tensor \begin{equation}\label{I}
		\hat{I_0}(z)=\f{1}{\rho^2\log^2{(\rho)}}|dz|^2
	\end{equation} on the unit disk $\mathbb{D}^{*}=\mathbb{D}\setminus 
	\{0\}\subseteq\mathbb{C}$ punctured at the origin, with $\rho$ the norm of $z=\rho e^{i\theta}$. This is a hyperbolic cusp: it can be obtained by pushing the hyperbolic metric on the cusp $\mathbb{H}^2/\langle \psi\rangle$, with $\psi$ the parabolic isometry $\psi(z)=z+2\pi$, via the diffeomorphism $f\colon \mathbb{H}^2/\langle \psi\rangle\rightarrow \mathbb{D}^{*}$ given by $f(z)=e^{iz}$.
	
	The conformal factor $\phi_0$ (also called the \textit{Liouville field}) associated to $\hat{I_0}$, such that $\hat I_0=e^{2\phi_0}\abs{dz}^2$, is:
	\begin{align}
		\phi_0(z)&=-\log(\abs{\rho\log(\rho)}) \label{firho} \\
		&=\log(2)-\f12\log(z\bar{z})-\log(\abs{\log(z\bar{z})}) \label{fiz}
	\end{align}
	for any $z\in\mathbb{D}^{*}$.
	
	We consider the upper half-space model $\mathbb{H}^3\cong\mathbb{C}\times \mathbb{R}^{+}$ for the $3$-hyperbolic space, whose boundary at infinity is conformally identified with $\mathbb{C}\cup \{\infty\}\cong \mathbb{CP}^1$. A point $p\in \mathbb{H}^3$ has then coordinates $p=(z,t)$, with $z=\rho e^{i\theta}\in \mathbb{C}$ its projection to the boundary at infinity, and $t\in \mathbb{R}^+$ its height.
	
	\begin{lem}\label{Epscoord}
		Let $\Sigma(\hat{I_0})\subseteq\mathbb{H}^3$ be the Epstein surface associated to $(\mathbb{D}^*, \hat{I_0})$. The point $s=s(z)=(r_0(z), \theta_0(z), t_0(z))$ in $\Sigma(\hat{I_0})$ whose projection through the geodesic ray based at $s$ and perpendicular to $\Sigma(\hat{I_0})$ is $z=\rho e^{i\theta}\in \mathbb{D}^{*}$ has the following coordinates:  
		
		\begin{align*}
			&r_0(z)=\left|\f{\rho(\log^2(\rho)-2)}{A(\rho)}\right|~,\\
			&\theta_0(z)=\theta+\pi\mathbbm{1}_{\rho\leq e^{-\sqrt{2}}}~, \\
			&t_0(z)=\f{-2\rho\log(\rho)}{A(\rho)}~,
		\end{align*}
		where \[A(\rho)=\log^2(\rho)+2\log(\rho)+2~,\]
		and $\mathbbm{1}_{\cdot}$ is the indicator function.
	\end{lem}
	
	\begin{figure}
		\begin{center}
			\includegraphics[width=.8\linewidth]{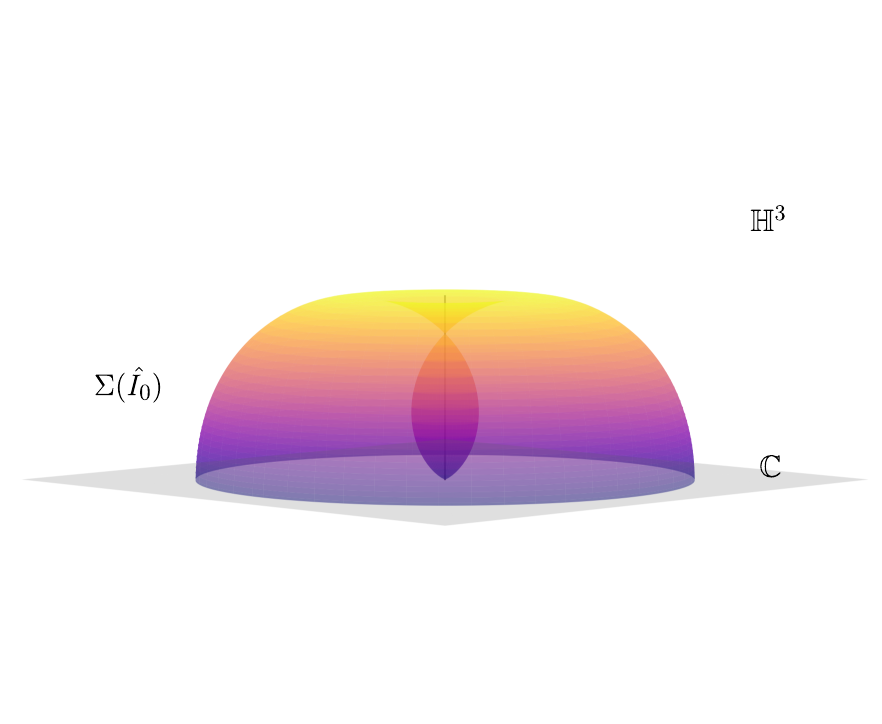}
			\caption{Epstein surface of the infinite cusp on the punctured disk with round boundary.}\label{Epsttot}
		\end{center}
	\end{figure}
	
	\begin{proof}
		First, note that since the metric at infinity $\hat{I_0}$ is invariant under rotations at the origin, the surface $\Sigma(\hat{I_0})$ is invariant under rotations along the vertical geodesic $\{0\}\times \mathbb{R}^{+}$ in $\mathbb{H}^3$. Therefore $r_0(z)=r_0(\rho)$ and $t_0(z)=t_0(\rho)$, for any $z=\rho e^{i\theta}\in \mathbb{D}^*$.
		Thanks to Equation \eqref{Sigma}, we can explicitly compute these two coordinates. We use the notation \[\phi_{0, \bar{z}}=\f{\partial \phi_0}{\partial\bar{z}}~,\] where we recall that in polar coordiantes \[\f{\partial}{\partial\bar{z}}=\f{e^{i\theta}}{2} \f{\partial}{\partial \rho}+ \f{ie^{i\theta}}{2\rho}\f{\partial}{\partial\theta}~,\] and, since $\phi_0$ depends just on $\rho$ \[\phi_{0, \bar{z}}=\f{e^{i\theta}}{2}\f{\partial \phi_0}{\partial \rho}=-\f{e^{i\theta}}2\f{\log(\rho)+1}{\rho\log(\rho)}=-\f{e^{i\theta}}2(\log(\rho)+1)e^{\phi_0}~.\]
		Therefore, from Equation (\ref{Sigma}), and using the identity $\nabla\phi_0=2\phi_{0,\bar{z}}$
		\begin{align*}
			r_0(z)&=\left|z+\f{4\phi_{0, {\bar{z}}}e^{-2\phi_0}}{1+|2\phi_{0,{\bar{z}}}|^2e^{-2\phi_0}}\right|=\left|\rho e^{i\theta}- \f{2e^{i\theta}(\log(\rho)+1)e^{-\phi_0}}{1+(\log(\rho)+1)^2}\right|\\
			&=|e^{i\theta}|\left|\rho - \f{2(\log(\rho)+1)\rho\log(\rho)}{1+(\log(\rho)+1)^2} \right|=\left|\f{\rho(2-\log^2(\rho))}{A(\rho)}\right|~,
		\end{align*}
		and
		\begin{align*}
			t_0(z)=\f{2e^{-\phi_0}}{1+|2\phi_{0, \bar{z}}|^2e^{-2\phi_0}}=\f{-2\rho\log(\rho)}{1+(\log(\rho)+1)^2}=\f{-2\rho\log(\rho)}{A(\rho)}~.
		\end{align*}
		Finally, we determine the angular coordinate, again via (\ref{Sigma}), simply by observing that \[\f{\rho(2-\log^2(\rho))}{A(\rho)}>0\ \text{\ iff\ }\ \rho\geq e^{-\sqrt{2}}~.\]

	\end{proof}
	
	%\blem\label{tangent}
	%The tangent $\text{Tg}(\Sigma(\hat{I_0}))$ to $\Sigma(\hat{I_0})$ at the point $(r_0(z), \theta_0(z), t_0(z))$, with $z=\rho e^{i\theta}$, satisfies \[\text{Tg}(\Sigma(\hat{I_0}))(\rho) = \f{2|1+\log(\rho)|}{A(\rho)}~,\]
	%where $A(\rho)=\log^2(\rho)+2\log(\rho)+2$.
	%\elem

	\brem
	Note that in a neighborhood of $0$ in $\mathbb{C}$, i.e. for $\rho\sim 0$, the coordinates of the Epstein surface $\Sigma(\hat{I_0})$ behave like \[(r_0(z), t_0(z)) \sim \left(\rho, \f{-2\rho}{\log(\rho)}\right)~.\]
	In particular, the surface $\Sigma(\hat{I_0})$ exits any $r$-neighborhood of the vertical geodesic of $\mathbb{H}^3$ based at $0\in\mathbb{C}$, with tangent going to zero as $-2/\log(\rho)$. Then, it is easy to see that the volume of the region in $\mathbb{H}^3$ between a hyperplane centered at the origin with small radius and $\Sigma(\hat{I_0})$ is infinite. Moreover, since the angular coordinate satisfies \[\theta_0(z)=\theta+\pi\mathbbm{1}_{\rho\leq e^{-\sqrt{2}}}~,\] for radii $\rho\leq e^{-\sqrt{2}}$, the outer normal is pointing inside the region bounded by the surface, see Figure \ref{Epsttot}.
	\erem

	\bdefi\label{Epstcut}
	Let $\Sigma(\hat{I_0})$ be the Epstein surface associated to $(\mathbb{D}^*, \hat{I_0})$ and $0<\rho_1<\rho_2<e^{-\sqrt{2}}$ be two radii. We denote by $\Sigma_{\rho_1}^{\rho_2}(\hat{I_0})$ the subset of $\Sigma(\hat{I_0})$ obtained by restricting its parameterization to the interior of the annulus \[\mathbb{D}_{\rho_1}^{\rho_2}=\{z\in\mathbb{C}\ |\ \rho_1\leq |z|\leq \rho_2\}~.\]
	As a subset of $\mathbb{H}^3$, the surface $\Sigma_{\rho_1}^{\rho_2}(\hat{I_0})$ is obtained by cutting $\Sigma(\hat{I_0})$ with the two hyperplanes $H_i$ in $\mathbb{H}^3$ with boundary at infinity coinciding with the two circles centered at $0$ and radii \[ h_i=\sqrt{r_0(\rho_i)^2+t_0(\rho_i)^2}\] respectively for $i=1,2$, i.e. \begin{equation}\label{sigma12}\Sigma_{\rho_1}^{\rho_2}(\hat{I_0})= \Sigma(\hat{I_0}) \cap \{(z,t)\in \mathbb{H}^3\  |\  h_1^2 \leq |z|^2+t^2 \leq h_2^2\}~.\end{equation}
	We will also denote by $N_{h_1}^{h_2}(\hat{I_0})$ the compact subset of $\mathbb{H}^3$ delimited by $\Sigma\rrho(\hat{I_0})$ and the hyperplanes $H_{i}$ of radii $h_i$, for $i=1,2$.
	\edefi

	\begin{figure}[h]
		\begin{center}
			\includegraphics[width=.8\linewidth]{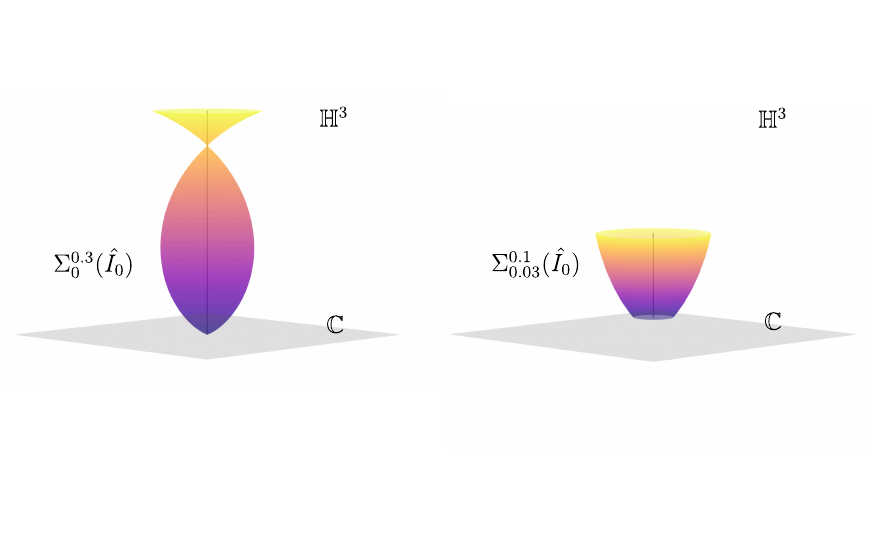}
			\caption{The Epstein surfaces associated to the interior of $\mathbb{D}_{\rho_1}^{\rho_2}$ equipped with the hyperbolic metric $ \hat{I_0}$, for two different examples of radii $\rho_i$. Note that $0.3 > e^{-\sqrt{2}} > 0.1$.}\label{cut12}
		\end{center}
	\end{figure}
	
	In what follows, we first compute the integral of the mean curvature $H_0$ over the surface $\Sigma_{\rho_1}^{\rho_2}(\hat{I_0})$, and then the volume of the compact region $N_{h_1}^{h_2}(\hat{I_0})$. This will be useful in the next section to study the behaviour of the $W$-volumes associated to the sequence of annuli with round boundary components $(\mathbb{D}_{\rho_1}^{\rho_2}, \hat{I_0})$ as $\rho_1\rightarrow 0$, which approximate a truncation of the infinite cusp $(\mathbb{D}^*, \hat{I_0})$.
	
	First, we recall how to express the mean curvature $H_0=\tr(I^{-1}_0\II_0)$ in terms of the first and second fundamental forms at infinity $\hat{I_0}$ and $\hat{\II_0}$, with $I_0$ and $\II_0$, respectively, the push-forward via the normal projection of the first and second fundamental forms on $\Sigma(\hat{I_0})$ induced by $\mathbb{H}^3$, (see Section \ref{DefI*} and Lemma \ref{meancurv}) : \[H_0=\f{1-\det(\hat{B_0})}{1+\tr(\hat{B_0})+\det(\hat{B_0})}~,\]
	where $\hat{B_0}=(\hat I)^{-1}_0\hat{\II_0}$.
	
	\blem\label{Meancurv0}
	In the notations above, the mean curvature $H_0$ of the surface $\Sigma_{\rho_1}^{\rho_2}(\hat{I_0})$ satisfies \[\f12 \int_{\Sigma_{\rho_1}^{\rho_2}(\hat{I_0})} H_0 da_{I_0} = \f{\pi}{12}\left[\log^3(\rho)\right]_{\rho_1}^{\rho_2}:= \f{\pi}{12}\log^3\left(\f{\rho_2}{\rho_1}\right)~,\]
	where $da_{I_0}$ denotes the area form associated to the first fundamental from $I_0$ induced by $\mathbb{H}^3$ on $\Sigma_{\rho_1}^{\rho_2}(\hat{I_0})$.
	\elem
	
	\begin{proof}
		We start by expressing the area form associated to $I_0$ in terms of that of $\hat{I_0}$. The pushed-forward through the normal projection of the first fundamental form on the Epstein surface is related to the data at infinity (see Section \ref{DefI*}) by the identity \cite[Lemma $5.6$]{S2008} \[I_0=\f14 \hat{I_0}((\id+\hat{B_0})(\cdot), (\id+\hat{B_0})(\cdot ))~,\] so that \[\sqrt{\det(I_0)}=\f 14 \sqrt{\det(\hat{I_0})}\det(\id+\hat{B_0})~,\]
		and therefore \[da_{I_0}=\f14 \det(\id+\hat{B_0})da_{\hat{I_0}}~.\]
		Since also \[\det(\id+\hat{B_0})=1+\tr(\hat{B_0})+\det(\hat{B_0})~,\] then, from the previous two identities and Lemma \ref{meancurv} \begin{equation}\label{HB}\f12 \int_{\Sigma_{\rho_1}^{\rho_2}(\hat{I_0})} H_0 da_{I_0}=\f18 \int_{\mathbb{D}_{\rho_1}^{\rho_2}} (1-\det(\hat{B_0})) da_{\hat{I_0}}~.\end{equation}
		It then remains to compute the determinant of the endomorphism $\hat{B_0}$. To this end, we write the second fundamental form at infinity in terms of the derivatives of the Liouville field $\phi_0$ (see Lemma \ref{secondfun}): \[\hat{\II_0}= 2q_0 dz^2 + 2\bar{q_0} d\bar{z}^2+4\phi_{0,{z\bar{z}}}|dz|^2\] with 
		\[q_0= \phi_{0,{zz}}-(\phi_{0,z})^2~.\]
		Since now (\ref{fiz}) \[\phi_0(z)=\log(2)-\f12\log(z\bar{z})-\log(\log(z\bar{z}))\]
		we can compute 
		\begin{align*}
			&\phi_{0,z} (z) = -\f{1}{2z} \left(1+\f{1}{\log(|z|)}\right)~, \\
			&\phi_{0,z\bar{z}} (z) = \f{1}{4|z|^2\log^2(|z|)}~,\\
			&\phi_{0,zz} (z) = \f{1}{2z^2} \left(1+ \f{1}{\log(|z|)}+\f{1}{2\log^2(|z|)}\right)~,
		\end{align*}
		then also \[q_0 (z) = \f{1}{4z^2}~.\]
		We observe that \[\hat{I_0}=\f{1}{|z|^2\log^2(|z|)}|dz|^2=4\phi_{0,z\bar{z}}|dz|^2~,\]
		moreover, by Remark \ref{secondfun2}
		\[\hat{\II_0}= (4\phi_{0,z\bar{z}}+4\re(q_0))dx^2+(4\phi_{0,z\bar{z}}-4\re(q_0))dy^2-4\im(q_0)dxdy\] 
		therefore 
		\[\hat{B_0}=(\hat{I_0})^{-1}I\!\hat{I_0}=\left(1+\f{\re(q_0)}{\phi_{0,z\bar{z}}}\right)dx^2+\left(1-\f{\re(q_0)}{\phi_{0,z\bar{z}}}\right)dy^2-\f{\im(q_0)}{\phi_{0,z\bar{z}}}dxdy\]
		and so 
		\[\det(\hat{B_0})=1-\f{1}{(\phi_{0,z\bar{z}})^2}\left(\re(q_0)^2+\im(q_0)^2\right)=1-\f{|q_0|^2}{(\phi_{0,z\bar{z}})^2}=1-\log^4(|z|)~.\]
		Finally, thanks to Equation (\ref{HB})
		\begin{align*}
			\f12 \int_{\Sigma_{\rho_1}^{\rho_2}(\hat{I_0})} H_0 da_{I_0}&=
			\f18 \int_{\mathbb{D}_{\rho_1}^{\rho_2}} (1-\det(\hat{B_0})) da_{\hat{I_0}}\\
			&=\f{\pi}{4}\int_{\rho_1}^{\rho_2} \log^4(\rho)\f{1}{\rho^2\log^2(\rho)} \rho d\rho \\
			&= \f{\pi}{4}\int_{\rho_1}^{\rho_2} \f{\log^2(\rho)}{\rho} d\rho \\
			&= \f{\pi}{4} \left[\f{\log^3(\rho)}{3}\right]_{\rho_1}^{\rho_2} = \f{\pi}{12}\log^3\left(\f{\rho_2}{\rho_1}\right)~.
		\end{align*}
	\end{proof}
	
	In the next lemma we compute the hyperbolic volume of the compact convex domain $N_{h_1}^{h_2}(\hat{I_0})$ defined in \ref{Epstcut}, where we remember that $h_i=\sqrt{r_0(\rho_i)^2+t_0(\rho_i)^2}$ with $0<\rho_1<\rho_2\leq e^{-\sqrt{2}}$. In particular, we focus on  how this diverges for $\rho_1\rightarrow0 $. The proof just consist in an elementary explicit computation and can be skipped.
	
	\blem\label{hypvol0}\label{Volume}
	The hyperbolic volume of the compact domain $N_{h_1}^{h_2}(\hat{I_0})$ is given by \[\text{Vol}(N_{h_1}^{h_2}(\hat{I_0}))=\f{\pi}{12}\log^3\left(\f{\rho_2}{\rho_1}\right)-\f{\pi}{2}\log\left(\f{\rho_2}{\rho_1}\right)+c(\rho_2)-c(\rho_1)~.\]
	with $c\colon (0,1)\rightarrow \mathbb{R}$ a smooth function, bounded on any sub-interval $(0,a)$ with $a<1$, and such that at $0$ \[ c(\rho) = O\left(\f{1}{\log(\rho)}\right)~,\] 
	where $O(x)$ stands for a real function such that the limit $\lim_{x\rightarrow 0} O(x)/x$ exists and it is finite.
	\elem
	
	\begin{proof}
		First, recall that we use the notation $(r_0(\rho), \theta+\pi, t_0(\rho))$ for the point in $\Sigma(\hat{I_0})$ associated through the normal projection to the point at infinity $z=\rho e^{i\theta}\in\mathbb{D}^{*}\subseteq \partial\bar{\mathbb{H}^3}$ (see Section \ref{EpstSurf}).
		We compute the volume of $N_{h_1}^{h_2}(\hat{I_0})$ as the sum of the signed volumes of the following three regions:
		\begin{align*}
			&O(\rho_1, \rho_2)=\{(z(\rho), t_0(\rho))\in \mathbb{H}^3\ |\ \rho_1\leq \rho \leq \rho_2,\ 0\leq |z(\rho)|\leq r_0(\rho) \}\\
			&S(\rho_i)=\{(z, t)\in \mathbb{H}^3\ |\ 0\leq |z|+t^2\leq h_i^2, \ t_0(\rho_i)\leq t \leq h_i\}
		\end{align*}
		for $i=1,2$, where \[ h_i= \sqrt{t_0(\rho_i)^2+r_0(\rho_i)^2}\] is the euclidean ray of the hyperplane whose intersection with $\Sigma(\hat{I_0})$ is the circle $\{(z, t_0(\rho_i))\in\mathbb{H}^3\ |\ |z|=r_0(\rho_i)\}$. Then, the hyperbolic volumes satisfy
		\[\text{Vol}(N_{h_1}^{h_2}(\hat{I_0}))=\text{Vol}(O(\rho_1, \rho_2))-\text{Vol}(S(\rho_1))+\text{Vol}(S(\rho_2)) ~.\]
		We now explicitly calculate the three terms. 
		The volume form of $\mathbb{H}^3$ in the cylindrical coordinates $(\rho, \theta, t)\in \mathbb{C}\times \mathbb{R}^+$ is \[d\text{Vol}_{\mathbb{H}^3}=\f{1}{t^3}\rho d\rho \wedge d\theta \wedge dt~, \] then, thanks to Lemma \ref{Epscoord}:\begin{align*}
			\text{Vol}(O(\rho_1, \rho_2))&=2\pi \int_{t_0(\rho_1)}^{t_0(\rho_2)} \int_{0}^{r_0(\rho)} \f{1}{t^3} \rho d\rho dt\\
			&=2\pi \int_{\rho_1}^{\rho_2} \f{r_0(\rho)^2}{2t_0(\rho)^3}t_0^{'}(\rho)d\rho \\
			&= \pi \int_{\rho_1}^{\rho_2} \f{A(\rho)^3}{8\rho^3\log^3(\rho)}\f{\rho^2(2-\log^2(\rho))^2}{A(\rho)^2}\f{2(\log(\rho)+1)(\log^2(\rho)+2)}{A(\rho)^2}d\rho\\
			&= \f{\pi}{4} \int_{\rho_1}^{\rho_2} \f{(2-\log^2(\rho))^2(\log(\rho)+1)(\log^2(\rho)+2)}{\rho\log^3(\rho)(\log^2(\rho)+2\log(\rho)+2)}d\rho\\
			&= \f{\pi}{4}\left[\f{\log^3(\rho)}{3}-\f{\log^2(\rho)}{2}-2\log(\rho)-4\log\left(\f{\log(\rho)}{\log^2(\rho)+2\log(\rho)+2}\right)-\f{2}{\log^2(\rho)}\right]_{\rho_1}^{\rho_2} \\
			&= \f{\pi}{4}\left[\f{\log^3(\rho)}{3}-\f{\log^2(\rho)}{2}-2\log(\rho)+4\log(\log(\rho))\right]_{\rho_1}^{\rho_2} \\  
			&+\f{\pi}{4}\left[4\log\left(1+\f{2}{\log(\rho)}+\f{2}{\log^2(\rho)}\right)-\f{2}{\log^2(\rho)}\right]_{\rho_1}^{\rho_2} ~;
		\end{align*}
		
		\begin{align*}
			\text{Vol}(S(\rho_i))&= 2\pi\int_{t_0(\rho_i)}^{ h_i} \int_0^{\sqrt{ h_i^2-t^2}} \f{1}{t^3}\rho d\rho dt \\
			&= \pi\int_{t_0(\rho_i)}^{ h_i} \f{ h_i^2-t^2}{t^3} dt \\
			&= \pi \left[-\f{ h_i^2}{2t^2}-\log(t)\right]_{t_0(\rho_i)}^{ h_i}\\
			&= -\f{\pi}{2}+\f{\pi}{2}\left(\f{ h_i}{t_0(\rho_i)}\right)^2 -\pi \log\left(\f{ h_i}{t_0(\rho_i)}\right) \\
			&= -\f{\pi}{2}+\f{\pi}{2}\left(1+\left(\f{r_0(\rho_i)}{t_0(\rho_i)}\right) ^2\right) - \f{\pi}{2} \log\left(1+\left(\f{r_0(\rho_i)}{t_0(\rho_i)}\right) ^2\right) \\
			&= \f{\pi(\log^2(\rho_i)-2)^2}{8\log^2(\rho_i)}-\f{\pi}{2}\log\left(\f{\log^4(\rho_i)+4}{4\log^2(\rho_i)}\right) \\
			&= \f{\pi}{8}\log^2(\rho_i)-\f{\pi}{2}+\f{\pi}{2\log^2(\rho_i)}-\f{\pi}{2}\log\left(\log^2(\rho_i)\left(1+\f{4}{\log^4(\rho_i)}\right)\right)+\f{\pi}{2}\log(4) \\
			&= \f{\pi}{8}\log^2(\rho_i)-\pi\log(\log(\rho_i))+\f{\pi}{2}\log(4)-\f{\pi}{2}+\f{\pi}{2\log^2(\rho_i)}-\f{\pi}{2}\log\left(1+\f{4}{\log^4(\rho_i)}\right)~;
		\end{align*}
		we then note that \[\text{Vol}(S_2)-\text{Vol}(S_1)=\left[\f{\pi}{8}\log^2(\rho)-\pi\log(\log(\rho))+\f{\pi}{2\log^2(\rho)}-\f{\pi}{2}\log\left(1+\f{4}{\log^4(\rho)}\right)\right]_{\rho_1}^{\rho_2}\]    
		and therefore \[\text{Vol}(N_{h_1}^{h_2}(\hat{I_0}))=\pi\left[\f{\log^3(\rho)}{12}-\f{\log(\rho)}{2}+c(\rho)\right]_{\rho_1}^{\rho_2}\]
		with \begin{equation}\label{c}
			c(\rho)=\pi\log\left(1+\f{2}{\log(\rho)}+\f{2}{\log^2(\rho)}\right)-\f{\pi}{2}\log\left(1+\f{4}{\log^4(\rho)}\right)~,
		\end{equation}
		which concludes the proof.
	\end{proof}

	\subsubsection{W-volume of a truncated cusp on a round annulus}
	%From Definition \ref{Wvol} and Lemmas \ref{Meancurv0}, \ref{hypvol0} and \ref{boundarycurvature}, we get the $W$-volume for $N_{h_1}^{h_2}(\hat{I_0})$.
	
	%\bprop\label{Wvol}
	%The following equality holds\marginnote{boundary term missing. Also, probably not needed} \[W(N_{h_1}^{h_2}(\hat{I_0}))=-\f{\pi}{2}\log\left(\f{\rho_2}{\rho_1}\right)+c(\rho_2)-c(\rho_1)+\alpha(\rho_1)\ell(\rho_1)-\alpha(\rho_2)\ell(\rho_2)~,\]
	%where $c$ is the function as in Lemma \ref{hypvol0} (defined explicitly by (\ref{c})) and is such that $c(x)$ is $o(x)$, and $\alpha$ and $\ell$ the ones at Lemma \ref{boundarycurvature}.
	%In particular \[\lim_{\rho_1\to 0}W(N_{h_1}^{h_2}(\hat{I_0}))= \f{\pi}{2}\log(\rho_1)+ C(\rho_2) \]
	%with \[C(\rho_2)=\pi/2\log(\rho_2)+c(\rho_2)-\alpha(\rho_2)\ell(\rho_2)\]
	%a constant that can be uniformly bounded in $\rho_2 < (1-\delta)$, for any fixed $\delta >0$.
	%\eprop
	
	We now consider the truncated cusp of the previous section as a domain with two round boundary components equipped with a hyperbolic metric, and compute its $W$-volume as in Definition \ref{Wboundary}.
	Let us then again consider the annulus \begin{equation}\label{D12}\mathbb{D}_{\rho_1}^{\rho_2}= \{z\in \mathbb{C}\ |\ \rho_1\leq |z|\leq \rho_2\}\end{equation} equipped with the restriction of the conformal metric $\hat{I_0}$, and also the associated sphere $S\mrho$ for domain with boundary defined as in Section \ref{BoundaryEps}: \[S\mrho=\partial N\left(\mathbb{D}_{\rho_1}^{\rho_2}, \phi_0\right)~,\]
	where $\phi_0$ is the Liouville field of $\hat{I_0}$, defined by Equation \eqref{firho}.
	Note that $S\mrho$ is obtained by connecting $\Sigma_{\rho_1}^{\rho_2}(\hat{I_0})$ (as in Definition \ref{Epstcut}) through the caterpillar region to the two hyperplanes whose boundary at infinity coincides with a component of $\partial\mathbb{D}_{\rho_1}^{\rho_2}$, and then capping it with these.
	
	We use the symbol $\sim$ to denote two equivalent functions at $0$, i.e. \[f(x)\sim g(x) \quad \Longleftrightarrow \quad \lim_{x\rightarrow 0}\f{f(x)}{g(x)}=1~.\]

	In the following statement, even though it diverges, we separate and highlight the mean curvature terms arising from the boundary of the domain, as this cancels out anytime we sum $W$-volumes of a decomposition of a domain along round circles (see Lemma \ref{add}), as we will do later.
	\begin{prop}\label{CuspConv} The $W$-volume associated to $(\mathbb{D}\rrho, \hat{I_0})$ has the following behaviour as $\rho_1\sim 0$: 
	\[W(\mathbb{D}\rrho, \hat{I_0}) = \f{\pi}{2}\log(\rho_1)-b(\rho_1)+C'(\rho_2)\]
	where $b(\cdot)$ represents the mean curvature boundary term, it is a continuous function, and at $\rho_1\sim 0$ satisfies \[b(\rho_1)=-\f{\pi^2}{8}\log(\rho_1)-\f{\pi}{2}\log(\rho_1)+O(1/\log \rho_1)~,\]
	and $C'(\rho_2)$ is a constant depending just on $\rho_2$. In particular, as $\rho_2\sim 0$ \[C'(\rho_2)= -\f\pi2\log(\rho_2)+b(\rho_2) ~.\]
	
	\end{prop}
	
	\begin{proof}
	In what follows, to streamline the notation, we will sometimes omit the dependence on $\hat{I_0}$, writing, for example, simply $S_{\rho_1}^{\rho_2}$. Also recall that we denote by $\Sigma(\cdot)$ the Epstein surface parameterized by the interior of the domain, while by $S(\cdot)$ the sphere obtained as the union of $\Sigma(\cdot)$, the caterpillar regions, and the hyperplanes whose boundaries match those of the domain. The $W$-volume associated to a domain with round boundary components is defined as (see Section \ref{BoundaryEps}, Definition \ref{Wboundary}): 
	\[ W(\mathbb{D}\rrho, \hat{I_0})= Vol(N\left(\mathbb{D}_{\rho_1}^{\rho_2}, \phi_0\right))-\f{1}{2}\int_{S\rrho} H da_{S\rrho}\]
	where $N(\mathbb{D}_{\rho_1}^{\rho_2}, \phi_0)$ is the ball filling $S_{\rho_1}^{\rho_2}$, so that $ S_{\rho_1}^{\rho_2}=\partial N(\mathbb{D}_{\rho_1}^{\rho_2}, \phi_0)$.
	First, we note that \[Vol(N(\mathbb{D}_{\rho_1}^{\rho_2}, \phi_0))= Vol(N_{\rho_1}^{\rho_2})+Vol(P_2)-Vol(P_1)\]
	where $N_{\rho_1}^{\rho_2}$ is the compact delimited by the Epstein surface $\Sigma(\hat{I_0})$ and the hyperplanes $H_i$ of euclidean radii $\rho_i$, and $P_i$, for $i=1,2$, is the region of $\mathbb{H}^3$ detected by $C_i$, the hyperplane $H_i$, and $\Sigma(\hat{I_0})$, such that $N(\mathbb{D}_{\rho_1}^{\rho_2}, \phi_0)=(N_{\rho_1}^{\rho_2}\cup P_2)\setminus P_1$.
	By Lemma \ref{Volume}:  
	\[\text{Vol}(N_{\rho_1}^{\rho_2})=\f{\pi}{12}\log^3\left(\f{x_2}{x_1}\right)-\f{\pi}{2}\log\left(\f{x_2}{x_1}\right)+c(x_2)-c(x_1)~.\]
	with $x_i=x_i(\rho_i)$ such that $\rho_i^2=h^2(x_i)=r_0^2(x_i)+t_0^2(x_i)$, and $c(x)$ a function such that at $0$ is $O(1/\log x)$. Since also, by Lemma \ref{Epscoord} \[r_0^2(x)+t_0^2(x)=\f{x^2(4+\log^4x)}{(\log^2x+2\log x+2)^2}~,\] then \[\rho_1=x_1(1+O(1/\log x_1))\] and therefore \begin{equation}\label{V}\text{Vol}(N_{\rho_1}^{\rho_2})=\f{\pi}{12}\log^3\left(\f{x_2}{\rho_1}\right)-\f{\pi}{2}\log\left(\f{x_2}{\rho_1}\right)+c(x_2)-O(1/\log \rho_1)~.\end{equation} 
	It is also easy to show (for example, looking at the area of the smallest rectangle containing a vertical section of $P_1$, seen as a revolution solid) that \[\lim_{\rho_1\rightarrow 0}Vol(P_1)=0~,\]
	while $Vol(P_2)$ is a finite number depending just on $\rho_2$, which again limits to zero at $\rho_2\sim 0$.
	We also notice that the integral of the mean curvature term splits in the following sum: \[\f{1}{2}\int_{S\rrho}Hda_{S\rrho}=\f{1}{2}\int_{\Sigma\rrho}Hda_{I_0}+\f12\int_{C_1}Hda_{C_1}-\f12\int_{C_2}Hda_{C_2}+\f\pi8\ell(\partial_1 C_1)-\f\pi8\ell(\partial_1 C_2)~,\]
	where $\Sigma_{\rho_1}^{\rho_2}$ denotes the Epstein surface $\Sigma(int(\mathbb{D}_{\rho_1}^{\rho_2}), \hat{I_0})$, and $C_i$, for $i=1,2$, is the caterpillar region, which meets the hyperplane $H_i$ with an exterior angle of $\pi/2$, and whose intersection with $H_i$ has hyperbolic length $\ell(\partial_1C_i)$.
	By Lemma \ref{Meancurv0}, arguing as for the volume, we already know that  \begin{equation}\label{H}\f{1}{2}\int_{\Sigma\rrho}Hda_{I_0}=\f{\pi}{12}\log^3\left(\f{x_2}{\rho_1}\right)+O(1/\log\rho_1)~.\end{equation}
	%	We are now going to prove that the sum of the terms in the mean curvature on $C_i$ and their area, compensate the additional volumes $Vol(P_i)$. Then the thesis will follow from Proposition \ref{Wvol}.
	What remains is to handle the mean curvature terms coming from the boundary $\partial\mathbb{D}_{\rho_1}^{\rho_2}$. We start by studying the integral of the mean curvature on the caterpillar regions.
	Thanks to the explicit parameterization of the Epstein surface associated to $\partial \mathbb{D}\rrho$ furnished by Equation \eqref{Pcaterpillar}, the caterpillar region $C_i$ is parameterized by the subset of the normal bundle of $\partial \mathbb{D}\rrho$ obtained by imposing $s\in [0, 2\pi\rho_i]$ and $v\in \left[ \f{\log(\rho_i)+1}{\rho_i\log(\rho_i)}, \f{1}{\rho_i}\right]$. The extremes for the parameter $v$ are computed by explicitly solving for the values where the caterpillar region intersects the hyperplane $H_i$ and the Epstein surface $\Sigma(\hat{I_0})$, respectively.
	% The coordinates of the caterpillar region are (see \cite[Section $2.2$]{BrockPalleteWbound}, or Section \ref{BoundaryEps}):
	%\begin{align*}&r_b(z)=-\rho+\f{2v}{e^{2\phi_0}+v^2}~,\\
	%	&t_b(z)=\f{2e^{\phi_0}}{e^{2\phi_0}+v^2}~, 
	%\end{align*} with \[e^{2\phi_0(z)}=\f{1}{\rho^2\log^2(\rho)}~.\]
	Recalling that $\hat{I_0}=e^{2\phi_0}|dz|^2$ and that $k_i=1/\rho_i$ is the euclidean curvature of the boundary circles of $\partial\mathbb{D}_{\rho_1}^{\rho_2}$, we have \cite[Section $2.2$]{BrockPalleteWbound}: \[Hda_{C_i}=\f{1}{2}e^{-2\phi}(2k_iv-v^2)dvds~.\]  
	%	\[-\f12\int_{C_i}Hda_{C_i}-\f32\int_{C_i}(1+H)da_{C_i}=-\f12\int_{C_i}Hda_{C_i}-2\pi\rho_i\f34\int_{1/\rho_i}^{\f{\log(\rho_i)+1}{\rho_i\log(\rho_i)}}dv\] whose behaviour at $\rho_1\sim 0$ is 
	Therefore,
	\begin{align*}
		-\f12\int_{C_i}Hda_{C_i}=& \int_{0}^{2\pi\rho_i}\int_{1/\rho_i}^{\f{log(\rho_i)+1}{\rho_i log(\rho_i)}}\f{1}{4}\rho_i^2\log^2(\rho_i)\left(\f{2v}{\rho_i}-v^2\right)dvds \\
		=& 	\f{\pi}{2}\rho_i^3\log^2(\rho_i)\left[\f{v^2}{\rho_i}-\f{v^3}{3}\right]_{1/\rho_i}^{\f{log(\rho_i)+1}{\rho_i log(\rho_i)}} \\
		=& \f{\pi}{2}\log(\rho_i)-\f{\pi}{6\log(\rho_i)}~.
	\end{align*} 
	We remark that, in the $W$-volume, the term relative to the second boundary component has inverted sign, as the caterpillar region has the opposite induced orientation.
	It only remains to take care of the mean curvature on the co-dimension one faces of $S_{\rho_1}^{\rho_2}$: \[\pm\f14\alpha_i\ell(\partial_1 C_i)~,\]
	where $\alpha_i=\pi/2$ is the exterior angle between the lower boundary $\partial_1C_i$ of the caterpillar region and the hyperplane $H_i$, and $\ell(\partial_1C_i)$ denotes its induced hyperbolic length. Again thanks to Equation \eqref{Pcaterpillar} (see also \cite[Section $2.2$]{BrockPalleteWbound}), and noting that $\phi_0$ is constant on each boundary component of $\partial\mathbb{D}_{\rho_1}^{\rho_2}$, it is possible to find the explicit radial and height coordinates in $\mathbb{H}^3$ of the Epstein surface associated to the boundary:
	\begin{align*}&r_{b}(s,v)=-\rho_i+\f{2v}{e^{2\phi_0}+v^2}~,\\
		&t_{b}(s,v)=\f{2e^{\phi_0}}{e^{2\phi_0}+v^2}~, 
	\end{align*} where $e^{2\phi_0(s)}=1/(\rho_i^2\log^2(\rho_i))$.
	It is now easy to calculate $\ell(\partial_1C_i)$:  \[\ell(\partial_1C_i)=\f{2\pi r_{b}(v_i)}{t_{b}(v_i)}=-\pi\log(\rho_i)+O(1/\log(\rho_i))~,\]
	where the last inequality holds since $v_i=1/\rho_i$.
	The statement now follows by defining $b(\rho_i)$ as \[b(\rho_i)=-\f{\pi^2}{8}\log(\rho_i)-\f{\pi}{2}\log(\rho_i)+O(1/\log(\rho_i))~,\] and summing up with (\ref{V}) and (\ref{H}). 
	\end{proof}
	
	\begin{figure}[h]
	\begin{center}
		\includegraphics[width=.8\linewidth]{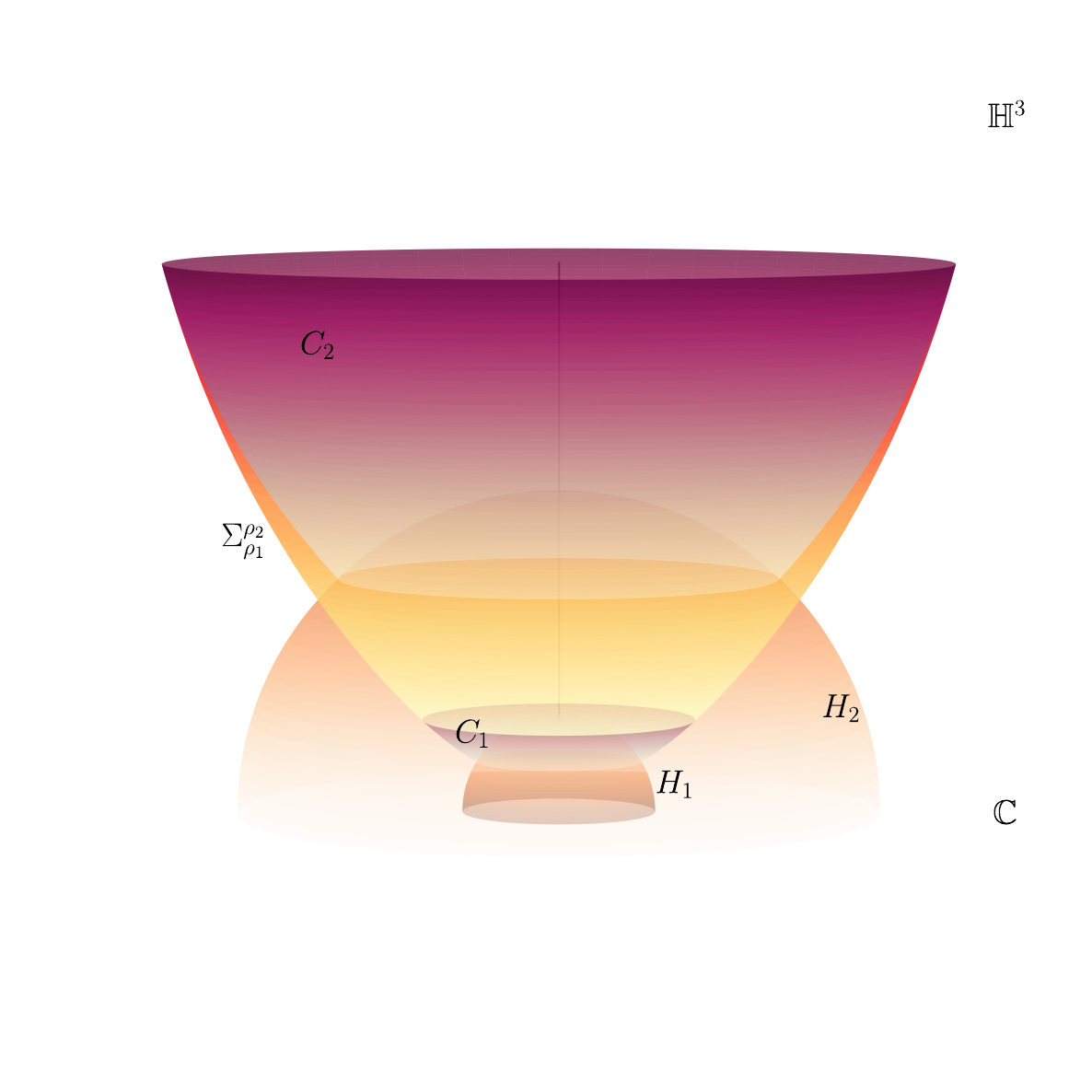}
		\caption{The Epstein surface associated to the domain with boundary $(\mathbb{D}_{\rho_1}^{\rho_2}, \hat{I_0})$: in yellow the Epstein surface $\Sigma_{\rho_1}^{\rho_2}(\hat{I_0})$ of the interior of the annulus, in purple the two caterpillar regions $C_i$, and in light brown the two hyperplanes $H_i$ of euclidean radii $\rho_i$.}\label{Epsteincut_boundaries}
	\end{center}
	\end{figure}
	
	\brem\label{1+H}
	The area form of the caterpillar region $C_1$, and its product with the mean curvature of $C_1$ diverge with the same rate, as (see \cite{BrockPalleteWbound}) \[Hda_{C_i}=\f{1}{2}e^{-2\phi}(-2k_iv+v^2)dv\wedge ds\] and \[da_{C_i}=\left(\f{1}{2}+\f{1}{2}e^{-2\phi}(2k_iv-v^2)\right)dv\wedge ds~.\]
	Thanks to this, the additional term appearing in Definition $2.2$ in \cite{BrockPalleteWbound} of $W$-volume coming from the innermost boundary of $(\mathbb{D}_{\rho_1}^{\rho_2}, \hat{I_0})$ vanishes in the limit $\rho_1\rightarrow 0$ as: \[\f32\int_{C_1}(1+H)da_{C_1}=2\pi\rho_1\f34\int_{1/\rho_1}^{\f{\log(\rho_1)+1}{\log(\rho_1)}}dv=\f{3\pi}{2\rho_1\log(\rho_1)}\rightarrow 0~.\]
	\erem
	
	Let us fix some notation. The hyperbolic cusp $[0, 2\pi]\times \mathbb{R}^{+}\subseteq \mathbb{H}^2$ is foliated by the curves given by the intersection between the horocycles at height $t\in\mathbb{R}^{+}$ centered at infinity and the cusp itself. This isometrically corresponds to foliating $(\mathbb{D}^{*}, \hat{I_0})$ by the circumferences centered at the origin of $\mathbb{C}$ and Euclidean radii $\rho=e^{-t}$. Equivalently, we are foliating a cusp in simple closed curves of hyperbolic length $\varepsilon=2\pi/\abs{\log(\rho)}=2\pi/t$.
	Then, the domain $\mathbb{D}(\varepsilon)\subseteq\mathbb{D}^{*}$ defined as \begin{equation}\label{DE} \mathbb{D}(\varepsilon)=\{\rho e^{i\theta}\in\mathbb{C} | 0<\rho\leq e^{-2\pi/\varepsilon} \}\end{equation} equipped with the restriction of the metric tensor $\hat{I_0}$ is a \textit{truncated cusp} with boundary of length $\epsilon$. We also define the function \begin{equation}\label{e}\rho(\varepsilon)=e^{-2\pi/\varepsilon}~,\end{equation}
	and recall the notation 
	\[\mathbb{D}_{\rho_1}^{\rho_2}= \{z\in \mathbb{C}\ |\ \rho_1\leq |z|\leq \rho_2\}~.\]
	
	\bdefi\label{VRcuspDef}
	The \textit{renormalized volume of a truncated cusp on a punctured disk with round boundary} of length $\bar\varepsilon$ is defined as \[V_R(\mathbb{D}(\bar{\varepsilon}))= \lim_{\rho\rightarrow 0}\left( W(\mathbb{D}_{\rho}^{\bar\rho}, \hat{I_0})-\f{\pi}{2}\log(\rho)+b(\rho)\right)~,\]
	with $\bar\rho=\rho(\bar\varepsilon)=e^{-2\pi/\bar\varepsilon}$, and $b(\rho)= -\pi^2\log(\rho)/8-\pi\log(\rho)/2+O(1/\log \rho)$ the boundary term as in Proposition \ref{CuspConv}.
	\edefi
	
	\bprop\label{RVcusp}
	For any $\bar{\varepsilon}>0$, the renormalized volume of the truncated cusp with round boundary $V_R(\mathbb{D}(\bar{\varepsilon}))$ exists and it is finite.
	\eprop
	
	\begin{proof}	
	This is a direct corollary of Proposition \ref{CuspConv}.
	\end{proof}
	
	\brem\label{RemVRcusp}
	Using the identity of Equation \eqref{e}, the renormalized volume of Definition \ref{VRcuspDef} can be expressed as a limit in the hyperbolic length $\varepsilon$ of the horocycles as \[V_R(\mathbb{D}(\bar{\varepsilon}))= \lim_{\epsilon\rightarrow 0}\left( W\left(\mathbb{D}_{\rho(\varepsilon)}^{\rho(\bar{\varepsilon} )},\hat{I_0}\right)+\f{\pi^2}{\varepsilon}+b(\rho(\varepsilon))\right)~,\]
	where \[b(\rho(\varepsilon))=\f{\pi^3}{4\varepsilon}+\f{\pi^2}{\varepsilon}+O(\varepsilon)~.\]
	\erem
	
	\subsubsection{General case}
	In the previous section, we analyzed the divergence of the $W$-volume of a cusp on a punctured disk with round boundary, assuming the conformal metric to be equal to the restriction of $\hat{I_0}$. We turn here to the general case of a cusp singularity at infinity, proving that, if we restrict to a small enough neighborhood, the same behaviour holds. Up to M\"obius transformation, we can assume the cusp singularity being at $p=0$.
	
	\bprop\label{VRcusp1}
	Let $\Omega\subseteq\mathbb{C} $ be a hyperbolic complex projective domain with $0\in \Omega$, and let $h_0$ be the unique conformal hyperbolic metric with a cusp singularity at $0$. Then, for $\bar{\rho}> 0$ small enough, the limit \[\lim_{\rho\rightarrow 0}\left(W(\mathbb{D}_{\rho}^{\bar{\rho}}, h_0)-\f{\pi}{2}\log(\rho)+b(\rho)\right)~,\]
	with $b(\rho)$ the boundary term as in Proposition \ref{CuspConv}, 
	exists and it is finite.
	\eprop
	
	\bpf
	As seen in Section \ref{ConfMarked}, in a small enough neighborhood $U\subseteq \Omega$ of $0$, there exists a (non projective) change of coordinate $w=ze^{\psi(z)}$, with $\psi$ a holomorphic function such that $\psi(0)=0$, in which the metric $h_0$ coincides with $\hat{I_0}$ (see Theorem \ref{cuspshape}, the discussion below, and Equation \eqref{changeofvar}). Then, in the projective coordinate $z$ of $\Omega$, in $U$, we can write $h_0$ as \begin{align*}
	h_0(z)&=\f{|e^{\psi(z)}+ze^{\psi(z)}\psi'(z)|^2}{|z|^2|e^{2\psi(z)}|\log^2(|ze^{\psi(z)}|)}|dz|^2\\
	&= \f{|1+z\psi'(z)|^2}{(1+Re(\psi(z))/(\log(|z|))^2}\hat{I_0}(z)~.
	\end{align*}
	Let us define $\nu\colon: \Omega\rightarrow \mathbb{R}$ as \begin{equation}\label{NU}
	\nu(z)=\log\left( \f{|1+z\psi'(z)|}{1+Re(\psi(z))/(\log(|z|)}\right)~, 
	\end{equation}
	so that \[h_0=e^{2\nu}\hat{I_0}~.\]
	By applying the Polyakov formula to estimate the difference of the $W$-volumes $W(\mathbb{D}_{\rho}^{\bar\rho}, e^{2\nu}\hat{I_0})$ and $W(\mathbb{D}_{\rho}^{\bar\rho}, \hat{I_0})$  (see Theorem \ref{Poly}), since $\nu(z)=O(\abs{z})$ and $\nu'(z)=\phi'(0)+O(\abs{z})$, the statement follows from Proposition \ref{RVcusp}.
	\epf

	\subsection{Renormalized volume of pointed hyperbolic 3-manifolds}\label{SecRVM}
	Building on the previous section, we are going to define the renormalized volume of the triple $(M, g, P)$, where $(M,g)=\mathbb{H}^3/\Gamma$ is a convex co-compact hyperbolic $3$-manifold, and $P\subseteq\Omega(\Gamma)/\Gamma=\partial_{\infty}M$ a finite set of marked points at infinity. We call the triple $(M, g, P)$ a \textit{pointed convex co-compact manifold}. We start with the case of a single point $P=\{p\}$, and then notice that the definition can be easily extended. Let $h$ be the hyperbolic metric in the conformal class of the boundary at infinity $[\partial_{\infty}M]$, and let $h_{p}$ be the unique conformal hyperbolic metric on $\partial\bar M\setminus \{p\}$ with a cusp singularity at $p$ (see Section \ref{ConfMarked}). With some abuse of notation, we denote again by $p$, $h$ and $h_{p}$ their $\Gamma$-invariant lifts respectively to $\Omega({\Gamma})$ and $\Omega({\Gamma})\setminus\Gamma\cdot p$. Up to Mobius transformation, we can assume that $p=0\in\mathbb{C}\subseteq\mathbb{CP}^{1}=\partial_{\infty}\mathbb{H}^3$, so that the restriction of $h_0$ to a small enough neighborhood of $0$ coincides with $\hat{I_0}$ (as in (\ref{I})), up to lower order terms described by Equation (\ref{NU}): $h_0=e^{2\nu}\hat{I_0}~$.
	Given a sufficiently small $\varepsilon>0$, we can cut $M$ with the hyperplane $H(\varepsilon)$ centered at $0$ and such that the intersection between the boundaries at infinity of each of its translates and the domain of discontinuity $\Omega(\Gamma)$ have length $\varepsilon$ with respect to $h_0$. We sometimes denote by $D(\varepsilon)$, with an abuse of notation, both the ball with boundary $\partial H(\varepsilon)$ containing $0$ in $\Omega(\Gamma)$ and in $\partial_{\infty}M$. This disk coincides with $\mathbb{D}(\varepsilon')=\{\rho e^{i\theta}\in\mathbb{C} | 0<\rho\leq e^{-2\pi/\varepsilon}\}$ where \[\varepsilon=\int_{\partial\mathbb{D}(\varepsilon')}e^{\nu}ds(\hat{I_0})~,\]
	with $\nu$ defined by Equation \eqref{NU}. We fix the following notations:  \[\Omega(\Gamma, \varepsilon):=\Omega(\Gamma) \setminus (\Gamma\cdot  D(\varepsilon))~,\] and analogously \[\partial_{\infty}M(\varepsilon):=\Omega(\Gamma, \varepsilon)/\Gamma~.\] We then consider the $\Gamma$-invariant Epstein surface $\widetilde{\Sigma}(h_0)$ in $\mathbb{H}^3$ associated to $(\Omega({\Gamma}), h_0)$, and the $\Gamma$-invariant Epstein surface $\widetilde{\Sigma}(h_0, \varepsilon)$ of $(\Omega({\Gamma}, \varepsilon), h_0)$ obtained by restricting the parameterization of $\widetilde{\Sigma}(h_0)$ to the open domain $\Omega({\Gamma}, \varepsilon)$, and connecting it to $\Gamma\cdot H(\varepsilon)$ through the caterpillar regions (as in \ref{BoundaryEps}). We also denote by $\Sigma(h_0)$ and $\Sigma(h_0, \varepsilon)$ the corresponding quotient Epstein surfaces in $M$. Finally, we define $W(\partial_{\infty}M(\varepsilon), h_0)$ to be the $W$-volume of the compact region in $M$ bounded by $\Sigma(h_0, \varepsilon)$ and $H(\varepsilon)$. The same construction could be performed in order to get the $W$-volume $W(\partial_{\infty}M(\varepsilon), g)$ for another conformal metric $g$, possibly with a puncture at the marked point $p$, where this time $\partial_{\infty}M(\varepsilon)$ is obtained by removing a round disk centered at $p$ of boundary length $\varepsilon$ with respect to $g$.

	\bdefi\label{cuspedVR}
	Let $(M, g, p)$ be a convex co-compact hyperbolic $3$-manifold pointed at $p\in \partial_{\infty}M$. Let also $\Gamma< PSL(2,\mathbb{C})$ be a representation of $\pi_1(M)$ such that $M$ is isometric to $\mathbb{H}^3/\Gamma$ and $p=0\in \Omega(\Gamma)$, and let $h_0$ be the unique conformal hyperbolic metric on $\partial_{\infty}M$ with a cusp singularity at $0$. In the notations above, we define the \textit{renormalized volume of $(M,g,p)$} as \[V_R(M,g, p)=\lim_{\varepsilon\rightarrow 0}\left( W(\partial_{\infty}M(\varepsilon), h_0)+\f{\pi^2}{\varepsilon} +b(\varepsilon)\right)~,\]
	with $b(\varepsilon)=\pi^3/(4\varepsilon)+\pi^2/\varepsilon+O(\varepsilon)$.
	\edefi
	
	\brem\label{pointdependence}
	Even if it is always possible to assume the marking point to be $0$, given two different marking points $p_1$ and $p_2$ in $\Omega(\Gamma)$, except in very special cases, there is no M\"obius transformation that restricts to a homeomorphism of $\Omega(\Gamma)$ and sends $p_1$ to $p_2$. Even if two pointed convex co-compact manifold $(M,g,p_1)$ and $(M,g,p_2)$ are isometric, their renormalized volume strictly depends on the choice of the marking point.
	\erem
	
	\bthm\label{Finiteness}
	The renormalized volume of a convex co-compact manifold $(M,g)$ pointed at $p\in\partial_{\infty}M$ exists and it is finite.
	\ethm
	
	\begin{proof}
	Let us fix a small enough $\bar{\varepsilon}>0$, so that, thanks to the additive property of the $W$-volume of Lemma \ref{add}, for any $0<\varepsilon<\bar\varepsilon$, we can split $W(\partial_{\infty}M(\varepsilon), h_0)$ as \[W(\partial_{\infty}M(\varepsilon), h_0)=W\left(D(\varepsilon, \bar{\varepsilon}), h_0\right)+W(\partial_{\infty}M(\bar{\varepsilon}), h_0)~,\]
	with $D(\varepsilon, \bar\varepsilon)=D(\bar\varepsilon)\setminus\text{int}(D(\varepsilon))$, and where, we recall, the boundary of $D(x)$ has length $x$ with respect to $h_0$. Then, using the notations introduced in (\ref{D12}) and (\ref{DE}), there exist two radii $0<\rho<\bar\rho$ depending respectively on $\varepsilon$ and $\bar\varepsilon$, such that \[D(\varepsilon, \bar\varepsilon)=\mathbb{D}_{\rho}^{\bar\rho}~,\]
	and \[D(\varepsilon)=\mathbb{D}(-2\pi/\log(\rho))~.\]
	Moreover, since $h_0=e^{2\nu}\hat{I_0}$, with $\abs{\nu(\rho e^{i\theta})}=O(\abs{\rho})$ at $\rho\sim 0$, then, for small radii \[\varepsilon=\int_{\partial\mathbb{D}(-2\pi/\log(\rho))}e^{\nu}ds(\hat{I_0}) \sim -\f{2\pi}{\log(\rho)}~.\] Therefore, since also the boundary term $b(\cdot)$ defined in Proposition \ref{CuspConv} is continuous \[\lim_{\varepsilon\rightarrow 0}\left( W(D(\varepsilon, \bar\varepsilon), h_0)+\f{\pi^2}{\varepsilon}+b(\varepsilon)\right)=\lim_{\rho\rightarrow 0}\left(W(\mathbb{D}_{\rho}^{\bar\rho}, h_0)-\f{\pi}{2}\log(\rho)+b(\rho)\right)~,\]
	and the right hand side exists and is finite by Proposition \ref{VRcusp1}.
	The renormalized volume of $(M,g,p)$ can then be expressed as
	\[V_R(M,g, p)=W(\partial_{\infty}M(\bar{\varepsilon}), h_0)+\lim_{\varepsilon\rightarrow 0}\left(W\left(D(\varepsilon, \bar{\varepsilon}), h_0\right)+\f{\pi^2}{\varepsilon}+b(\varepsilon)\right)
	%=\widetilde{V_R}(\mathbb{D}(\varepsilon))+ W(\partial_{\infty}M(\bar{\varepsilon}), \hat{I_0})
	~.\]
	It then remains to show that the first term in the equality is finite. First, note that $\bar{\varepsilon}$ can also be chosen such that $h_0$ and $h$ stay at uniformly bounded small distance on $\partial_{\infty}M(\bar{\varepsilon})$, where $h$ is the unique hyperbolic representative in $[\partial_{\infty}M]$ (see Section \ref{ConfMarked}). Therefore, the difference between $W(\partial_{\infty}M(\bar{\varepsilon}), h_0)$ and $W(\partial_{\infty}M(\bar{\varepsilon}), h)$ is bounded. Moreover, the last one is less than the renormalized volume of $M$, plus the boundary term $b(\bar\varepsilon)$ coming from $\partial D(\bar\varepsilon)$, which, for a fixed $\bar{\varepsilon}$, is also finite.
	\end{proof}
	
	\brem
	We point out that, in the classical definition of the $W$-volume, it is often necessary to rescale the metric at infinity by a constant factor to ensure the associated Epstein surface is embedded (see Proposition \ref{Epstflow}, Theorem \ref{EpstEmb}, and Section \ref{WvolVR}), and the well-definedness is guaranteed by the identity $W(M,g)=W(M, e^{2r}g)+\pi r\chi(\partial_{\infty}M)$. This continues to work in the setting of $W$-volumes for domain with round boundary components $\Omega$, since, thanks to the Polyakov formula (Theorem \ref{Poly}), for any $r>0$ \[W(\Omega, e^{2r}g)-W(\Omega, g)=-\f14\int_{\Omega} K(g)rda_{\Omega}-\f12\int_{\partial\Omega}k(g)rds_{\partial\Omega}~,\]
	and therefore, by the Gauss-Bonnet Theorem and the fact that the scalar curvature $K(g)$ is twice the Gaussian curvature \[W(\partial_{\infty}M(e^{r}\varepsilon), e^{2r}g)-W(\partial_{\infty}M(\varepsilon),g)=-\pi r\chi(\partial_{\infty}M(\varepsilon))~.\]
Note that, anyways, the singular point on the Epstein surface associated with a metric having a hyperbolic cusp cannot be removed by rescaling (see Figure \ref{Epsttot}).
	\erem
	
	\brem\label{P}
	We can extend Definition \ref{cuspedVR} to a finite set of marked points $P\subseteq \partial_{\infty}M$ of cardinality $|P|$ as \[V_R(M, g, P)=\lim_{\varepsilon\rightarrow 0}\left( W(\partial_{\infty}M(\varepsilon), h_P) +|P|\left(\f{\pi^2}{\varepsilon}+b(\varepsilon)\right)\right)~,\] where $h_P$ is the unique complete hyperbolic metric conformal to $\partial_{\infty}M$ with cusps at each point in $P$, and \[\partial_{\infty}M(\varepsilon)=\Omega(\Gamma, \varepsilon)/
	\Gamma\quad \text{with}\quad\Omega(\Gamma, \varepsilon)=\Omega(\Gamma)-\bigcup_{i=1}^{|P|}\Gamma\cdot D_i(\varepsilon)~,\]
	where $D_i(\varepsilon)$ denotes the euclidean ball of boundary coinciding with the one of the hyperplane $H(\varepsilon)$ used to cut the $i$-th puncture. The proof of Theorem \ref{Finiteness} for the finitedness of $V_R(M, g, P)$ applies straightforwardly. 
	\erem
	
	\bcor
	The renormalized volume of a convex co-compact manifold $(M,g)$ pointed at a finite set of points $P\subseteq \partial_{\infty}M$ exists and it is finite.
	\ecor

	\section{Adapted renormalized volume}\label{RVA}
	In this section, we finally define the adapted renormalized volume for hyperbolic manifold with compressible boundary. In Theorem \ref{AVRbounded}, we show that, unlike the classical version, it is bounded from below. Moreover, we prove that its differential has uniformly bounded norm (Theorems \ref{boundedGradient} and \ref{boundedGradientWP}). We then proceed to prove that the adapted renormalized volume extends by continuity to points of the boundary of the Teichm\"uller space representing the pinching of a compressible multicurve $-$ whose renormalized volume was defined in the previous section (see Definition \ref{cuspedVR}), and for which we define the adapted version later in this section.
	
	We recall that the deformation space of the convex co-compact structures on a tame manifold $\bar M$ is identified with the quotient of the Teichm\"uller space of the boundary $\partial\bar M $ by the subgroup generated by all Dehn twists along compressible simple closed curves (see Theorem \ref{UnifCCM}). Accordingly, we will use the notation $V_R(X)$ for $V_R(M(X))$, where $M(X)$ is the convex co-compact hyperbolic $3$-manifold associated to the Riemann surface $X$ via Theorem \ref{UnifCCM}.

	\bdefi\label{adaptedVR}
	Given $M$ a convex co-compact hyperbolic 3-manifold,
	we define the \textit{adapted renormalized volume} as the function \[\widetilde{V_R}\colon \mathcal{T}(\partial\bar{M})\rightarrow \mathbb{R}\] such that \[\widetilde{V_R}(X)=V_R(X)+ \max_{m\in\Gamma^{comp}(\partial\bar M)} L(X, m)\] with 
	\[ L(X, m)=\pi^3\sum_{\substack{\gamma\in m }} \f{1}{\ell_{\gamma}(X)} \]
	where $\Gamma^{comp}(\partial\bar M)$ denotes the set of compressible multicurves on $\partial\bar M$, and $\ell_{\gamma}(\cdot)$ the length function of $\gamma$.
	\edefi
	
	\brem\label{maximality}
	In the definition above, we can also assume that the multicurves $m$ running under the maximum are \textit{maximal}, that is, any other compressible simple closed curve in $\partial \bar M$ intersects $m\in\Gamma^{comp}(\partial\bar M)$ non trivially, so that $m$ cannot be extended.
	\erem
	
	\blem\label{maximumexistence}
	For every $X\in\mathcal{T}(\partial\bar M)$, the supremum of $L(X, m)$ over all the compressible multicurves $m\in\Gamma^{comp}(\partial\bar M)$ is attained. Therefore, the adapted renormalized volume function is well defined.
	\elem
	
	\bpf
	By contradiction, let us assume $\{m_i\}_i\subseteq \Gamma^{comp}(\partial\bar M)$ to be an infinite sequence of pairwise distinct multicurves such that \[L(X, m_i)\nearrow \sup_{m\in\Gamma^{comp}(\partial\bar M)}L(X,m)~.\] Let then $\gamma_0$ be a simple closed compressible curve of $\partial \bar M$, and let $\ell_0$ denote its hyperbolic length in $X$. Let $g$ be the sum of the geni of the connected components of $\partial\bar M$. Every multicurve on $\partial\bar M$ has at most of $3g-3$ components (see \cite[Theorem $4.1.1$]{Bu1992}). Choose now an $L>0$ such that $(3g-3)/L < 1/\ell_0$. Since in every hyperbolic surface, for any fixed $L$, there are finitely many simple closed curves of length less than $L$, there exists a large enough $j\in\mathbb{N}$ such that for all $i\geq j$ every curve in $m_i$ has length greater than $L$. Then \[L(X, m_i)\leq \f{\pi^3(3g-3)}{L}< \f{\pi^3}{\ell_0}=L(X, \gamma_0)~,\]
	contradicting the assumption that $\{m_i\}_i$ is a sequence realizing the supremum. 
	\epf
	
	\blem\label{finitemaximum}
	For every compact subset $K\subseteq\mathcal{T}(\partial\bar M)$, the maximum of Definition \ref{adaptedVR} can be taken over a finite subset $\sigma(K)\subseteq \Gamma^{comp}(\partial \bar M)$. In particular, for every $X\in \mathcal{T}(\partial\bar M)$ there are finitely many multicurves realizing the maximum.
	\elem
	
	\bpf
	Fix an $m\in\Gamma^{comp}(\partial \bar M)$, and consider the function \[L(\cdot, m)\colon K\rightarrow \mathbb{R}_{\geq 0}~,\]
	as in Definition \ref{adaptedVR}, which is continuous and it then attains maximum and minimum on $K$. Therefore, there exists an $L_0>0$ such that, for every $X\in K$ \[L_0\geq L(X,m)>0~.\] Choose an $L>0$ such that $(3g-3)/L \leq 1/L_0$, and for any subset $U\subseteq \mathcal{T}(\partial\bar M)$ define the subset of simple closed curves \[\sigma_0(U)=\{\gamma\in \pi_1(\partial\bar M) \text{ compressible } |\ \exists X\in U\ \colon \ell_{\gamma}(X)\leq L\}~.\] We claim that $\sigma_0(K)$ is finite. Let us consider the Thurston (asymmetric) distance $d_{Th}$, which is (forward) compatible with the topology of $\mathcal{T}(\partial \bar M)$, and fix a small $\delta>0$. We can then cover $K$ with a finite number of balls $B_{\delta}(X_i)$ of radius $\delta$, for $i=1,\dots, k$. Since \[d_{Th}(X,Y)=\sup_{\gamma}\log\left(\f{\ell_{\gamma}(X)}{\ell_{\gamma}(Y)}\right)~,\] for all $Y\in B_{\delta}(X)=\{Y\in\mathcal{T}(\partial\bar M)\ |\ d_{Th}(X,Y)\leq\delta\}$ and for all $\gamma\in\pi_1(\partial \bar M)$ there exists $\delta_1>0$ such that  \[\abs{\ell_{\gamma}(Y)-\ell_{\gamma}(X)}\leq \delta_1~.\] Therefore, as for every $X_i$ there are a finite number of simple closed curves of bounded length,
	the set $\sigma_0(B_{\delta}(X_i))$ is finite. Since the union of the balls  $B_{\delta}(X_i)$ covers $K$, then also $\sigma_0(K)$ is finite. We conclude by taking $\sigma(K)$ as the multicurves formed by subsets of curves in $\sigma_0(K)$.
			\epf
	
	\bcor\label{contLip}
	The adapted renormalized volume function is continuous.
	\ecor
	
	\bpf
	The adapted renormalized volume is locally defined as the sum of the renormalized volume, which is continuous, and, by Lemma \ref{finitemaximum} and local compactness of $\mathcal{T}(\partial\bar M)$, the maximum of finitely many continuous functions, which is again continuous. 
	\epf
	
	\brem\label{failsmooth}
	The adapted renormalized volume is smooth outside the subset of Teichm\"uller space where the multicurve realizing the maximum in Definition \ref{adaptedVR} is not unique.
	\erem
	
	\brem\label{welldef}
	The simple length spectrum is invariant under the action of the mapping class group, and Dehn twists along compressible simple closed curves preserve the property of a curve being compressible or not. Therefore, the adapted renormalized volume of a convex co-compact manifold is well defined, through Theorem \ref{UnifCCM}, as \[\widetilde{V_R}(M(X))=\widetilde{V_R}(X)~.\]
	\erem
	
	%\brem\label{discontinuity}
	%Outside of the union of codimension-one manifolds in $\mathcal{T}(\partial\bar M)$, each corresponding to the set of surfaces with a compressible simple closed geodesic of length equal to $\varepsilon_0$, the adapted renormalized volume is real analytic.
	%\erem
	
	\bdefi
	Let $M$ be a convex co-compact manifold. For any $\delta>0$, we call \textit{$\delta$-compressible thick part} of the Teichm\"uller space, and we denote it by $\mathcal{T}_{\delta}^{comp}(\partial\bar M)\subseteq \mathcal{T}(\partial \bar M)$, the subset of marked Riemann surfaces having length of the shortest (with respect to the hyperbolic metric) compressible simple closed curve bounded below by $\delta$. We call \textit{$\delta$-compressible thin part} the complement of $\mathcal{T}_{\delta}^{comp}(\partial\bar M)$ in $\mathcal{T}(\partial \bar M)$.
	\edefi
	
	Recall, that any simple closed curve on a hyperbolic surface of length less than the Margulis constant $\varepsilon_0=2\argsinh(1)$ admits a long collar, outside of which the surface has injectivity radius greater than $\varepsilon_0$ \cite[Theorem $4.1.6$]{Bu1992}.
	
	\blem\label{var1}
	There exists a small enough $\varepsilon_1>0$ depending only on the genus of $\partial \bar M$ such that $\varepsilon_1< \varepsilon_0$, and, for any $X\in\mathcal{T}(\partial\bar M)$, if $m\in\Gamma^{comp}(\partial \bar M)$ realizes the maximum of $L(X,\cdot)$, then $m$ contains all the simple closed compressible curves of length less than $\varepsilon_1$ in $X$. 
	\elem
	
	\bpf
	By contradiction, suppose that for all $\varepsilon>0$ there exist a Riemann surface $X\in\mathcal{T}(\partial\bar M)$, a multicurve $m\in\Gamma^{comp}(\partial \bar M)$ realizing the maximum of $L(X, \cdot)$, and a simple closed compressible curve $\gamma$ of length less than $\varepsilon$ in $X$ and such that $\gamma\notin m$. By Remark \ref{maximality}, we can assume $m$ to be maximal. Then, there exists $\alpha\in m$ such that $i(\gamma, \alpha)\neq 0$, for $i(\cdot, \cdot)$ the geometric intersection. Let us consider now the multicurve $m_1$ obtained from $m$ by adding $\gamma$ and removing the set $m_0\subseteq m$ of all the curves of $m$ intersecting $\gamma$. In this way \[L(X,m)=L(X, m_1)-\f{\pi^3}{\ell_{\gamma}(X)}+L(X, m_0)~.\] Since $\ell_{\gamma}(X)\leq \varepsilon_1$, and every $\alpha\in m_0$ intersects $\gamma$ non trivially, by the Collar Lemma \[L(X,m)-L(X, m_1)\leq -\f{\pi^3}{\varepsilon}+\f{\pi^3(3g-3)}{\argsinh(1/\sinh(\varepsilon/2))}\] that, for $\varepsilon$ small enough, is strictly negative , contradicting the fact that $m$ realizes the maximum. Let us choose an $\varepsilon_1=\varepsilon_1(g)>0$ such that the member of the inequality above is negative. Since $\varepsilon_0=2\argsinh(1)$, it is straightforward to verify that $\varepsilon_1<\varepsilon_0$.
	\epf

	\bthm\label{AVRbounded}
	For every convex co-compact hyperbolic $3$-manifold $M$, the adapted renormalized volume $\widetilde{V_R}(\cdot)$ is bounded from below by a constant depending just on the topology of the boundary $\partial\bar M$.
	\ethm
	
	\bpf
	First, note that for all $X\in \mathcal{T}(\partial\bar M)$ \[\widetilde{V_R}(X)\geq V_R(X)~.\]
	Fixed a small $\delta>0$, we begin by considering the $\delta$-compressible thick part of the Teichm\"uller space $\mathcal{T}_{\delta}^{comp}(\partial\bar M)$.
	% For all $X\in \mathcal{T}_{\delta}^{comp}(\partial\bar M)$, and for any multicurve $m\in\Gamma^{comp}(\partial \bar M)$ \[0<L(X, m)\leq \f{\pi^3(3g-3)}{\delta}~.\] In particular, $\widetilde{V_R}(X)\geq V_R(X)$.
	 Merging Theorem $2.16$ in \cite{BBB2018}, which furnishes an upper bound on the length of the bending lamination, and Theorem $A.6$ in \cite{averages}, which states that the $W$-volume of the convex core (see Section \ref{WvolVR}) and the renormalized volume have distance bounded by a constant depending just on the genus of the boundary, for any $X\in\mathcal{T}_{\delta}^{comp}(\partial\bar M)$ \[\widetilde{V_R}(X)\geq V_R(X)\geq V_C(X)-\f32\pi|\chi(\partial\bar M)|\coth^2(\delta/4)-C~,\]
	where $V_C(\cdot)$ is the function associating to the convex co-compact manifold $M(\cdot)$ the volume of its convex core, and $C$ is a constant depending just on the Euler characteristic of the boundary $\chi(\partial\bar M)$. Let us now fix $\delta=\varepsilon_1$ as in Lemma \ref{var1}, and let $X$ be a point in the complement of $\mathcal{T}_{\varepsilon_1}^{comp}(\partial\bar M)$. Let then $m_1=m_1(X)$ be the set of compressible curves in $X$ of length less than $\varepsilon_1$. By Lemma \ref{var1}, every compressible multicurve $m$ realizing the maximum of $L(X, \cdot)$ contains $m_1$. Recall that $m$ has cardinality at most $3g-3$, with $g$ the sum of the geni of the components of the surface. We order the $k$ components $\gamma_i$ of $m_1$, and we consider the path in $\mathcal{T}(\partial\bar M)$ given by the concatenation of $k$ (negative) grafting paths, each of whom on $\gamma_j$ and terminating in a surface $X_j$ in which the geodesic representatives of the first $j$ curves in $m$ have length $\varepsilon_0$, and $X_0=X$. Then, the ending point $X_k$ belongs to $\mathcal{T}_{\varepsilon_1}^{comp}(\partial\bar M)$. Note that for any $X,Y\in\mathcal{T}(\partial\bar M)$ and any couple of multicurves $m,p\in\Gamma^{comp}(\partial\bar M)$ realizing the maximum of $L(X, \cdot)$ and $L(Y, \cdot)$ respectively \[\abs{L(X, m\setminus m_1(X))-L(Y, p\setminus m_1(Y))}< c_{g}(\varepsilon_1)\eqdef\f{(3g-3)\pi^3}{\varepsilon_1}~.\] Thanks to Theorem $1.4$ in \cite{CG}
	\[V_R(X_{j-1})=V_R(X_j)-\f{\pi^3}{\ell_{j}(X_{j-1})}+\f{\pi^3}{\ell_{j}(X_j)}-\f\pi4\left(\ell_{j}(X_j)-\ell_{j}(X_{j-1})\right)+O\left(e^{-\pi s_j/(2\ell_j(X_j))}s_j^3\right)~,\]
		where $\ell_j(X_i)$ denotes the length of the curve $\gamma_j$ with respect to the hyperbolic metric in $X_i$, and $s_j$ is the parameter of the $j$-th grafting path.
	 Therefore, we can estimate the variation of the adapted renormalized volume under the just mentioned grafting path as: \begin{align*}
	\abs{\widetilde{V_R}(X_k)-\widetilde{V_R}(X)}&\leq\sum_{j=1}^{k-1} \abs{\widetilde{V_R}(X_j)-\widetilde{V_R}(X_{j-1})}=\sum_{j=1}^{k-1} \abs{V_R(X_j)-V_R(X_{j-1})-\f{\pi^3}{\ell_{j}(X_{j-1})}}+c_{g}(\varepsilon_1) \\
	&\leq \sum_{j=1}^{k-1}\abs{\f{\pi^3}{\ell_{j}(X_j)}-\f\pi4\left(\ell_{j}(X_j)-\ell_{j}(X_{j-1})\right)+O\left(e^{-\pi s_j/(2\ell_j(X_j))}s_j^3\right)}+c_{g}(\varepsilon_1) \\
	&=\sum_{j=1}^{k-1}\abs{\f{\pi^3}{\varepsilon_0}-\f\pi4\left(\varepsilon_0-\ell_{j}(X_{j-1})\right)+O\left(e^{-\pi s_j/(2\varepsilon_1)}s_j^3\right)}+c_{g}(\varepsilon_1)~.
	\end{align*}
 It is now easy to see that the last expression is uniformly bounded, in particular it is always smaller than \[(3g-3)\left(\f{2\pi^3}{\varepsilon_1}+\f{\pi\varepsilon_1}4+A\left(\f{6\varepsilon_1}{e\pi}\right)^3\right)~,\]
	for some universal constant $A$ furnished by \cite[Remark $5.23$]{CG} bounding the terms in the $O(\cdot)$. Since for $X_k$ holds the first estimate, the statement follows.
	\epf
	
	In the next theorem, we prove that the $L^1$-norm of the differential of the adapted renormalized volume is bounded. We recall that the $L^1$-norm of a covector of the Teichm\"uller space at $X=\mathbb{H}^2/\Gamma$, i.e. of a $\Gamma$-invariant holomorphic quadratic differential $q(z)dz^2$ on $\mathbb{H}^2$, is defined as \begin{equation}\label{inftynorm}\norm{q(z)dz^2}_{1}=\int_X \abs{q}~.\end{equation}
	%The dual of this norm is the \textit{Teichm\"uller norm}  on the space of harmonic Beltrami differential (see Section \ref{norms}).
	
	\bthm\label{boundedGradient} Let $M$ be a convex co-compact hyperbolic manifold.
	At the points where it exists, the differential of the adapted renormalized volume is bounded in the $L^1$ norm by a constant that depends only on the topology of the boundary $\partial\bar M$.
	\ethm 
	
	\bpf
	%	We recall that, at any $X\in\mathcal{T}(\partial\bar M)$, the Weil-Petersson isomorphism between $\mathcal{T}_X^*(\partial\bar M)$ and $\mathcal{T}_X(\partial\bar M)$ sends $\phi$ to $-\bar{\phi}/\rho_X$, with $\rho_X$ the hyperbolic metric on $X$. Then, the statement of the theorem is equivalent to asserts that the differential of the adapted renormalized volume has uniformly bounded $L^1$-norm as a holomoprhic quadratic differential.
	
	Similarly to the proof of Theorem \ref{AVRbounded}, first, we consider any point in the $\varepsilon_1$-compressible thick part $X\in \mathcal{T}^{comp}_{\varepsilon_1}(\partial\bar M)$, where $\varepsilon_1$ is as in Lemma \ref{var1}. We also assume that there exists a unique multicurve $\bar m\in\Gamma^{comp}(\partial\bar M)$ realizing the maximum of $L(X, \cdot)$. The adapted renormalized volume is smooth at $X$ and there exists a neighborhood $U$ of $X$ in $\mathcal{T}(\partial\bar M)$ such that for all $Y\in U$ \[L(Y, \bar m)=\max_{m\in\Gamma^{comp}(\partial\bar M)} L(Y, m)~.\] Therefore, if $\bar m=\{\gamma_j\}_{1}^k$, with $k\leq (3g-3)$, denoting by $\ell_{j}(\cdot)$ the hyperbolic length function of $\gamma_j$, in the neighborhood $U$ of $X$
	\[\widetilde{V_R}(\cdot)=V_R(\cdot)+\sum_{j=1}^{k}\f{\pi^3}{\ell_{j}(\cdot)}~,\] so that \begin{equation}\label{dv1}(d\widetilde{V_R})_X=\re(\mathcal S(f_X))-\sum_{j=1}^{k}\f{\pi^3}{\ell^2_{j}(X)}(d\ell_j)_X~,\end{equation} 
	where $f_X$ is the developing map of the boundary at infinity of $M(X)$. By Lemma $5.1$ in \cite{KM} (see also Corollary $2.12$ in \cite{BBB2018}): \[\norm{\re(\mathcal S(f_X))}_1\leq 2\pi\abs{\chi(X)}\norm{\re(\mathcal S(f_X))}_{\infty}\leq 3\pi\abs{\chi(X)}\coth^2(\varepsilon_1/4)~.\] Moreover, since $X\in \mathcal{T}^{comp}_{\varepsilon_1}(\partial\bar M)$, for every $j=1,\cdots, k$, the lengths $\ell_j(X)$ are bounded below by $\varepsilon_1$. In the second part of the proof, we will show the following claim \begin{equation}\label{toshow}
		\norm{\f{\pi^3}{\ell^2_{j}(X)}(d\ell_j)_X}_1 \leq \f{2\pi^3}{\ell_j(X)}\leq \f{2\pi^3}{\varepsilon_1}~,
	\end{equation}
	from which, together with the bound above, it follows the uniform bound in the $\varepsilon_1$-compressible thick part of Teichm\"uller space.

	Let us now study the differential of the adapted renormalized volume at a point $X$ in the complement of $\mathcal{T}^{comp}_{\varepsilon_1}(\partial\bar M)$, and of the subset where $\widetilde{V_R}$ is not smooth (see Remark \ref{failsmooth}). As before, there exists a non-empty compressible multicurve $\{\gamma_j\}_{1}^{k}$ with $k\leq 3g-g$ components, containing the multicurve $\{\gamma_j\}_{1}^{k_1}$ formed by all the $k_1\leq k$ compressible curves of length strictly less than $\varepsilon_1$ in $X$ (by Lemma \ref{var1}), such that \begin{equation}\label{dv1}(d\widetilde{V_R})_X=\re(\mathcal S(f_X))-\sum_{j=1}^{k}\f{\pi^3}{\ell^2_{j}(X)}(d\ell_j)_X~.\end{equation}
	Let us fix some notation. We denote by $X_{thick}\subseteq X$ the thick part of $X$, that is, the maximal subsurface in $X$ that has injectivity radius bigger than $\varepsilon_1/2$ (since $\varepsilon_1<\varepsilon_0$, this splitting is always possible, see for example Chapter $4$ in \cite{Mar16}). The subsurface $X_{thick}$ contains the complement of the union of the thin tubes around any simple closed curve of lenght $< \varepsilon_1$ (see Definition \ref{margulitubedefi}, and \cite[Theorem $4.1.6$]{Bu1992}). We now define $X_{\varepsilon_1}$ as the subsurface of $X$ given by the union of the thin tubes $\mathcal{A}_j$ around $\gamma_j$, for  $1\leq j \leq k_1\leq k$. Let also $D$ be a fundamental domain for $\Gamma$ such that $X=\mathbb{H}^2/\Gamma$. For any $j=1,\cdots, k$, we denote by $\varphi_{j}$ the hyperbolic isometry corresponding to $\gamma_j$, and by $X_{j}=\mathbb{H}^2/\langle \varphi_{j}\rangle$ the respective quotient annulus, which is a covering of $X$. Anytime $j$ is fixed and we focus on a tubular neighborhood $\mathcal{A}_j$ of $\gamma_j$, up to conjugation, we can choose the lift $\widetilde{\gamma_j}$ of $\gamma_j$ given by the geodesic of $\mathbb{H}^2$ thorough $\{0, \infty\}$. Moreover, we can identify $X_j$ with the annulus \[\{z\in\mathbb{H}^2\ |\ 1\leq \abs{z}\leq e^{\ell_{j}(X)}\}~,\] and assume $D$ to contain $ \widetilde{\mathcal{A}_j}$, with $\widetilde{\mathcal{A}_j}$ the lift of $\mathcal{A}_j$ through $\pi_j\colon X_j\rightarrow X$ intersecting $\widetilde{\gamma_j}$.
	
	The $L^1$-norm of $(d\widetilde{V_R})_X$ can be split into the sum of the $L^1$-norms of its restriction to $X_{\varepsilon_1}$ and its complement $X_{\varepsilon_1}^c$: \begin{equation}\label{01}\norm{(d\widetilde{V_R})_X}_1=\int_X \abs{(d\widetilde{V_R})_X}= \int_{X_{\varepsilon_1}}\abs{(d\widetilde{V_R})_X}+\int_{X_{\varepsilon_1}^c}\abs{(d\widetilde{V_R})_X}~. \end{equation}
	
	We start by bounding the $L^1$-norm in $X_{\varepsilon_1}^c$, which is contained in the union of the thick part $X_{thick}$ and the incompressible thin tubes. The norm of the Schwarzian in this region is bounded as before by $3\pi\abs{\chi(X)}\coth^2(\varepsilon_1/4)$. We need to bound the $L^1$-norm of $(d\ell_j)_X/\ell^2_j(X)$.
	To this end, let us recall the Gardiner's formula for the differential of the length function of $\gamma_j$ (see \cite{Gard}):
	\begin{equation}\label{dl} d\ell_j(\mu)= \f{2}{\pi}\re\left\langle \f{dz^2}{z^2}, \widetilde{\mu }\right\rangle_{X_{j}}\end{equation}
	where $\mu$ is a harmonic Beltrami differential in $T_{X}\mathcal{T}(\partial\bar{M})$, and $\widetilde{\mu}$ denotes its lift to the annulus $X_{j}$. Note that, in the setting fixed above, the differential $2dz^2/(\pi z^2)$ is $\phi_j$-invariant.
	It might be tempting to assume $(d\ell_j)_X$ to be equal to $2dz^2/(\pi z^2)$ by the non-degeneracy of the Weil-Petersson pairing, but, since the pairing in (\ref{dl}) is on the whole annulus $X_j$, this is not true. In fact, the differential $dz^2/z^2$ is not invariant under $\Gamma$, as instead (the lift of) $(d\ell_j)_X$ is. Gardiner's formula can though be rephrased as a pairing on the fundamental domain $D$ by taking the holomorphic quadratic differential giving by the $\Theta$-series of $2dz^2/(\pi z^2)$ with respect to the projection $\pi_j\colon X_j\rightarrow X$ (see \cite[Chapter $5$]{Hubbard2016}, \cite[Theorem $3$]{Gard}), so that: \begin{equation}\label{thetaseries}(d\ell_j)_X=(\pi_j)_*\left(\f{2dz^2}{\pi z^2}\right)=\f{2}{\pi}\sum_{\beta\in\Gamma/\langle\phi_j\rangle} \beta^*\left(\f{dz^2}{z^2}\right)~.\end{equation}
	It is now easy to see that the $L^1$-norm of $(d\ell_j)_X$ on $D$, and therefore on $X$, is bounded by the $L^1$-norm of $2dz^2/(\pi z^2)$ on $X_j$: 
	\begin{align*}\norm{(\pi_j)_*\left(\f{dz^2}{z^2}\right)}_1&=\int_D \abs{\sum_{\beta\in\Gamma/\langle\phi_j\rangle} \beta^*\left(\f{dz^2}{z^2}\right)}\leq \sum_{\beta\in\Gamma/\langle\phi_j\rangle}\int_D \abs{\beta^*\left(\f{dz^2}{z^2}\right)}=\\
	&=\sum_{\beta\in\Gamma/\langle\phi_j\rangle}\int_{\beta(D)} \abs{\f{dz^2}{z^2}}=\int_{X_j} \abs{\f{dz^2}{z^2}}=\norm{\left(\f{dz^2}{z^2}\right)}_1~. \end{align*}
	When $j>k_1$, we can just explicitly compute $\norm{2dz^2/(\pi z^2)}_1$ on $X_j$, finding $2\ell_j(X)$, and obtaining the claimed Equation (\refeq{toshow}).
	Otherwise, it is enough to bound the norm of $dz^2/(\ell^2_j(X)z^2)$ on the complement of $\widetilde{\mathcal{A}_j}$ in $X_j$, as this contains the lift of $X^c_{\varepsilon_1}$ via $\pi_j$. More precisely, this lift is contained in the subset (see \cite[Lemma $5.2.7$]{Mar16}) \begin{equation}\label{sintheta}D\setminus \widetilde{\mathcal{A}_j}=\{z=\rho e^{i\theta}\in\mathbb{H}^2\ |\ 1\leq \rho\leq e^{\ell(\gamma_j)},\ 0<\sin(\theta_j)<1/\cosh(L_j)\}~,\end{equation} where \[L_j=\argsinh\left(1/\sinh(\ell_j(X)/2)\right)~,\]
	is the width of the thin tube $\mathcal{A}_j$.
	%Then, the infinity norm of the restricion of $dz^2/z^2$ to $X_j\setminus \widetilde{\mathcal{A}_j}$, with $z=x+iy=\rho e^{i\theta}$, is bounded as 
	%\[\sup_{z\in X_j\setminus \widetilde{\mathcal{A}_j}} \abs{\f{y^2}{z^2}}=\sup_{z\in X_j\setminus \widetilde{\mathcal{A}_j}} \abs{\f{y^2}{z^2}}\abs{\f{\rho^2\sin^2(\theta)}{\rho^2}}\leq \left(\left(\f{1}{\sinh(\ell/2)}\right)^2+1\right)^{-1}=\left(\f{\sinh(\ell/2)}{\cosh(\ell/2)}\right)^2~,\] so that
	%\[\f{1}{\ell_X^2(\gamma_j)}\sup_{z\in X_j\setminus \widetilde{\mathcal{A}_j}}\norm{\f{dz^2}{z^2}}\leq \f12~,\]
	Then, omitting the dependence by $X$ of the lengths \begin{equation}\label{L1thick}\int_{X_{\varepsilon_1}^c}\f{\abs{(d\ell_j)_X}}{\ell_j^2}\leq \f{1}{\ell_j^2}\int_{D\setminus\widetilde{\mathcal{A}_j} }\abs{(d\ell_j)_X} \leq\f{2}{\pi\ell_j^2}\int_{X_j\setminus\widetilde{\mathcal{A}_j}}\abs{\f{dz^2}{z^2}}=\f{4\theta_j}{\pi\ell_j}\leq \f{2}{\pi}\end{equation} where the last inequality follows from
	\[\theta=\arcsin(1/\cosh(L_j))\leq \ell_j/2 ~.\] Therefore, from Equation (\ref{dv1})
	\begin{equation}\label{02}\int_{X_{\varepsilon_1}^c}\abs{(d\widetilde{V_R})_X}\leq 3\pi\abs{\chi(X)}\coth^2(\varepsilon_1/4)+ 2\pi^2k_1+(k-k_1)\f{2\pi^3}{\varepsilon_1}~,\end{equation}
	where we recall that $k_1\leq k\leq (3g-3)$.
	
	We now proceed to estimate the $L^1$-norm of $(d\widetilde{V_R})_X$ on the union $X_{\varepsilon_1}$ of the compressible thin tubes. Let us fix a $j$ and use the same notations as above. Thanks to Theorem $1.1$ in \cite{CG}, again denoting $\ell_{\gamma_j}(X)$ simply by $\ell_j$, for any $j\leq k_1$, the Schwarzian derivative on $\widetilde{\mathcal{A}_j}$ satisfies \[\mathcal S(f)=\f{1}{2z^2}\left(1+\f{4\pi^2}{\ell_j^2}\right)dz^2+q_{\ell_j}(z)dz^2~,\] with \[\abs{q_{\ell_j}}=O\left(\f{e^{-\pi^2/(2\ell_j)}}{\ell_j^2}\right)\leq A~,\]
	where $A$ is a universal constant furnished by \cite[Remark $5.23$]{CG}.
	Therefore, by splitting $(d\ell_j)_X$ as 
	\[(d\ell_j)_X=\f{2}{\pi}\sum_{\substack{\beta\in\Gamma/\langle\phi_j\rangle \\ \beta\neq [\phi_j]}} \beta^*\left(\f{dz^2}{z^2}\right)+\f{2}{\pi}\f{dz^2}{z^2}~,\]
	on $\widetilde{\mathcal{A}_j}$ Equation \eqref{dv1} can be written as  \begin{equation}\label{dVrthin}(d\widetilde{V_R})_X=\f{1}{2}\f{dz^2}{z^2}+q_{\ell_j}(z)dz^2-\f{2\pi^2}{\ell_j^2}\sum_{\substack{\beta\in\Gamma/\langle\phi_j\rangle \\ \beta\neq [\phi_j]}} \beta^*\left(\f{dz^2}{z^2}\right)-\sum_{\substack{i=1 \\ i\neq j}}^{k}\f{\pi^3}{\ell_i^2}(d\ell_i)_X~.\end{equation}
	It is easy to see that the $L^1$-norm of the first two terms is bounded on $\widetilde{\mathcal{A}_j}$ by a universal constant $B$. Moreover, the $L^1$-norm on $\widetilde{\mathcal{A}_j}$ of the first sum, that we denote now by $\Theta_0(dz^2/z^2)$, satisfies \[\f{2\pi^2}{\ell_j^2}\norm{\Theta_0\left(\f{dz^2}{z^2}\right)}_1\leq \f{2\pi^2}{\ell_j^2}\int_{X_j\setminus D} \abs{\f{dz^2}{z^2}}\leq  \f{2\pi^2}{\ell_j^2}\int_{X_j\setminus \widetilde{\mathcal{A}_j}}\abs{\f{dz^2}{z^2}}\leq 2\pi^2~,\]
	where the last inequality follows from the proof of the other case. Lastly, observing that, for any $i\neq j$, $\pi_i^{-1}(\mathcal{A}_j)\subseteq X_i\setminus \widetilde{\mathcal{{A}}_i}$, the $L^1$-norm on $\widetilde{\mathcal{{A}}_j}$ of each term $(d\ell_i)_X/\ell^2_i(X)$ in the last sum is bounded again as in (\ref{L1thick}). Finally, using the estimate furnished by Equation (\ref{toshow}) for the curves with indexes $j\geq k_1$, we obtain the estimate \begin{equation}\label{03}\int_{X_{\varepsilon_0}}\abs{(d\widetilde{V_R})_X}=\sum_{j=1}^k\int_{\mathcal{A}_j}\abs{(d\widetilde{V_R})_X}\leq (B+2\pi^2+2\pi^2(k_1-1))k_1+(k-k_1)\f{4\pi^3}{\varepsilon_1}.\end{equation}
	Since $k_1\leq k\leq 3g-3$, the bounds in (\ref{02}) and (\ref{03}) depend only on the genus $g$ of $X$, and the statement of the theorem follows from (\ref{01}).	
	\epf
	
	In the next theorem, we prove that the differential of the adapted renormalized volume is also bounded with respect to the infinity norm. Before stating it, we recall the definition \[\norm{q}_{\infty}=\sup_{z\in X}\norm{q}=\sup_{z\in X}\abs{\f{q}{\rho_X}}~,\]
	for $q\in Q(X)$ and where $\rho_X$ stands for the hyperbolic metric conformal to $X$. 
	This result is of particular interest, as it implies that the Weil-Petersson gradient of the adapted renormalized volume \[(\nabla\widetilde{V_R})_X=\f{\overline{dV_R}}{\rho_X}\] has bounded $L^{\infty}$ norm, i.e., bounded \textit{Teichm\"uller norm}. This last is a Finsler metric, as it is not induced by a scalar product. In turn, up to a factor $1/2$, the Teichm\"uller metric induces the \textit{Teichm\"uller distance} $d_{\text{Teich}}$ on $\mathcal{T}(S)$ (see \cite[Theorem $6.6.5$]{Hubbard2016}), with respect to which Teichm\"uller space is complete (see, for example, \cite[Proposition $11.17$]{FM2011}). This open the possibility of studying the flow lines of minus the Weil-Petersson gradient of the adapted renormalized volume.
	
	\bthm\label{boundedGradientinfty} Let $M$ be a convex co-compact hyperbolic manifold.
	At the points where it exists, the differential of the adapted renormalized volume is bounded in the $L^{\infty}$ norm by a constant that depends only on the topology of the boundary $\partial\bar M$.
	\ethm 
	
	\bpf
	As in the previous proof, if $X$ belongs to the $\varepsilon_1$-compressible thick part $\mathcal{T}_{\varepsilon_1}^{comp}(\partial\bar M)$, then, the statement follows from the already used bound on the infinity norm of the Schwarzian derivative given in Lemma $5.1$ in \cite{KM}: \[\norm{(d\widetilde{V_R})_X}_{\infty}\leq \f32 \coth^2(\varepsilon_1/4)~,\] together with a uniform bound on $\norm{\pi^3(d\ell_\gamma)_X/\ell^2_{\gamma}(X)}_{\infty}$, for $\gamma$ a simple closed compressible curve in $X$, which we show later on in this proof.
	Otherwise, following the same strategy and notations as before, we start by bounding the $L^{\infty}$ norm on the complement $X_{\varepsilon_1}^c$ of the union of the compressible  thin tubes $\mathcal{A}_j$, indexed by $j=1,\dots, k_1\leq k\leq 3g-3$. Here, the norm of the Schwarzian is again bounded by $3\coth^2(\varepsilon_1/4)/2$. Next, we have to bound the $L^{\infty}$ norm of $\pi^3d\ell_j/\ell_j^2$, on the lift $\widetilde{X_{\varepsilon_0}^c}$ of $X_{\varepsilon_0}^c$ to a fundamental domain $D\subseteq X_j$, with $X_j=\mathbb{H}^2/\langle z\rightarrow e^{\ell_j}z\rangle$ the covering associated to the core of the thin tube $\mathcal{A}_j$. In this setting (see Equation \eqref{thetaseries}) \[\pi^3\f{d\ell_j}{\ell_j^2}=2\pi^2 \sum_{\beta\in\Gamma/\langle\phi_j\rangle} \beta^*\left(\f{dz^2}{z^2}\right)~,\]
	with $\ell_j$ the hyperbolic lengths of the simple closed curves $\gamma_j$ in $X$ of length less than $\varepsilon_1$, and $\phi_j$ the corresponding hyperbolic isometries. Before going on with the proof, we remark that the bound on the infinity norm on the thick part of $X$ follows directly from Theorem \ref{boundedGradient}, as in general, if $q\in Q(X)$, then on $X_{thick}$ \[\sup_{z\in X_{thick}}\norm{q}\leq C\norm{q}_{1}~,\]
	for some constant $C$ depending just on the topology of the surface $X$. Nevertheless, we provide a more explicit argument, as it will be needed for the second part of the proof, where we handle the estimate on the thin tubes.
	
	For any fixed $j\leq k$, consider the lift $\widetilde{\gamma_j}=\{0,\infty\}$ of the core of $\mathcal{A}_j$, and choose $D$ to be the fundamental domain in $X_j$ intersecting $\widetilde{\gamma_j}$. Denoting by $\widetilde{\mathcal{A}_j}$ the lift of $\mathcal{A}_j$ contained in $D$ \[\sup_{z\in D\setminus \widetilde{X}_{\varepsilon_1}^c}\f{\pi^3}{\ell_j^2}\norm{d\ell_j}\leq \f{2\pi^2}{\ell_j^2}\sum_{\beta\in\Gamma/\langle\phi_j\rangle}\sup_{z\in\beta(D\setminus \widetilde{A}_j)}\norm{\f{dz^2}{z^2}}=\f{2\pi^2}{\ell_j^2}\sum_{\beta\in\Gamma/\langle\phi_j\rangle}\sup_{z\in\beta(D\setminus \widetilde{A}_j)}\abs{\f{y^2}{z^2}}~,\]
	where $z=x+iy=\rho e^{i\theta}$, so that $y$ is the imaginary part of $z$ and $\theta$ is its angular coordinate, and then 
	\[\sup_{z\in\beta(D\setminus \widetilde{A}_j)}\abs{\f{y^2}{z^2}}=\sup_{z\in\beta(D\setminus \widetilde{A}_j)}\abs{\sin^2(\theta)}~.\]
	In the case $\beta=[\phi_j]\in \Gamma/\langle\phi_j\rangle$, we compute 
	\[\sup_{z\in D\setminus \widetilde{\mathcal{A}_j}}\f{2\pi^2}{\ell_j^2} \abs{\f{y^2}{z^2}}=\sup_{z\in D\setminus \widetilde{\mathcal{A}_j}}\f{2\pi^2}{\ell_j^2} \abs{\f{\rho^2\sin^2(\theta)}{\rho^2}}\leq\f{2\pi^2}{\ell_j^2}\left(\f{\sinh(\ell_j/2)}{\cosh(\ell_j/2)}\right)^2\leq \pi^2~,\] where the first inequality follows from \eqref{sintheta}\[\sin(\theta)\leq \f{1}{\cosh(L_j)}=\left(\cosh\left(\argsinh\left(1/\sinh(\ell_j/2)\right)\right)\right)^{-1}= \left(\left(\f{1}{\sinh(\ell_j/2)}\right)^2+1\right)^{-1}~,\]
	with $L_j$ the width of the thin tube $\mathcal{A}_j$ (see Definition \ref{margulitubedefi}, and \cite[Lemma $5.2.7$]{Mar16}). For all the other elements $\beta\in\Gamma/\langle\phi_j\rangle$, the Prime Geodesic Theorem \cite{Huber1959}, implies that, asymptotically, for a given natural number $N$, there are approximately $e^{N}/N$ translates of the fundamental domain $D$ at hyperbolic distance $\asymp L_j+N$ from $\widetilde{\gamma_{j}}$. Therefore \[\f{2\pi^2}{\ell_j^2}\sum_{\substack{\beta\in\Gamma/\langle\phi_j\rangle \\ \beta\neq [\phi_j]}}\sup_{z\in\beta(D\setminus \widetilde{A}_j)}\abs{\sin^2(\theta)}\lesssim \f{2\pi^2}{\ell_j^2}\sum_{N}\abs{\f{e^{N}}{N\cosh^2(L_j+N)}}\lesssim \f{2\pi^2e^{-2L_j}}{\ell_j^2}\sum_{N}\abs{\f{e^{-N}}{N}}~,\]
	and, since the series converges and $e^{-2L_j}/\ell_j^2$ is bounded as before, we get \[\sup_{z\in D\setminus \widetilde{X}_{\varepsilon}^c}\f{\pi^3}{\ell_j^2}\norm{d\ell_j}\leq A~,\]
	for some constant $A>0$.
	Therefore, we proved \[\sup_{z\in{X_{\varepsilon_1}^c}}\norm{(d\widetilde{V_R})_X}\leq 3\coth^2(\varepsilon_0)/2+(3g-3)A~.\]
	Note that a similar argument works to bound the infinity norm of $\pi^3d\ell_j/(\ell_j^2)$ for $j\geq k_1$ on the whole surface: this time $\ell_j$ is bounded from below while $\L_j$ from above. In particular, this concludes the proof in the case $X$ belongs to the $\varepsilon_1$-compressible thick part of Teichm\"uller space. 
	
	Let us now consider the union $X_{\varepsilon_1}$ of the compressible thin tubes. As already shown in the proof of Theorem \ref{boundedGradient}, on any thin tube $\mathcal{A}_j\subseteq X_{\varepsilon_1}$, the diverging term of the Schwarzian derivative simplifies with the one of the holomorphic quadratic differential $\pi^3d\ell_j/\ell_j^2$, so that (see Equation \eqref{dVrthin}), on $\widetilde{\mathcal{A}_j}$:
	\[(d\widetilde{V_R})_X=\f{1}{2}\f{dz^2}{z^2}+q_{\ell_j}(z)dz^2-\f{2\pi^2}{\ell_j^2}\sum_{\substack{\beta\in\Gamma/\langle\phi_j\rangle \\ \beta\neq [\phi_j]}} \beta^*\left(\f{dz^2}{z^2}\right)-\sum_{\substack{i=1 \\ i\neq j}}^{k}\f{\pi^3}{\ell_i^2}(d\ell_i)_X~.\]
	A direct and easy computation shows that the infinity norm of the first two terms is uniformly bounded. Furthermore, each differential appearing in the last sum can be estimated as in the previous case, thanks to the containment $\widetilde{\mathcal{A}_j}\subseteq X_i\setminus \widetilde{\mathcal{A}_i}$. The same holds for the remaining terms in the $\Theta$-series of $d\ell_j$, since for any $\beta\neq[\phi_j]$, the translates satisfy $\beta(\widetilde{\mathcal{A}_j})\subseteq X_j\setminus \widetilde{\mathcal{A}_j}$.
	
	\epf
	
	As a corollary of the previous two theorems, we obtain also the bound on the $L^2$ norm, that is, on the \textit{Weil-Petersson norm} (see Section \ref{norms}).
	
	\bcor\label{boundedGradientWP}
	Let $M$ be a convex co-compact hyperbolic manifold.
	At the points where it exists, the Weil-Petersson gradient of the adapted renormalized volume is bounded in the Weil-Petersson norm by a constant that depends only on the topology of the boundary $\partial\bar M$.
	\ecor 
	
	\bpf
	The Weil-Petersson norm on harmonic Beltrami differentials coincides with the $L^2$-norm on dual holomoprhic quadratic differentials. That is, for any $q\in Q(X)$ \[\norm{q}_2^2= \int_X \f{q\bar q}{\rho_X}=\left\langle q, \f{\bar q}{\rho_X}\right\rangle=\norm{\f{\bar q}{\rho_X}}_2^2~,\]
	where $\rho_X$ denotes the conformal hyperbolic metric on $X$, and $\langle \cdot,\cdot\rangle$ the Weil-Petersson pairing (see Section \ref{norms}). Moreover, the $L^2$-norm is bounded by the $L^{\infty}$ norm and the $L^1$-norm as \begin{equation}\label{12infty}\norm{q}_2^2\leq \norm{q}_{\infty}\norm{q}_1~.\end{equation}
	Therefore, the statement follows directly from Theorems \ref{boundedGradient} and \ref{boundedGradientinfty}.
	\epf

	Others corollaries of Theorem \ref{boundedGradientinfty} are the following. The first is the Lipschitz continuity of the adapted renormalized volume. 
	
	\bcor\label{WPdistance}
	Let M be a convex co-compact hyperbolic manifold. For any couple of points $X_0, X_1\in \mathcal{T}(\partial\bar M)$, the difference of the corresponding adapted renormalized volumes is bounded by the Weil-Petersson distance between $X_0$ and $X_1$ as \[\abs{\widetilde{V_R}(X_1)-\widetilde{V_R}(X_0)}\leq \sqrt{2\pi\abs{\chi(\partial\bar M)}}C(\partial\bar M)d_{WP}(X_0, X_1)~,\]
	where $C(\partial\bar M)$ is a constant depending only on the topology of $\partial\bar M$.
	\ecor
	
	\bpf
	Let us consider the length minimizing Weil-Petersson geodesic $\alpha(s)$ connecting $X_0$ to $X_1$, and let $\mu_s$ be the tangent harmonic Beltrami differential given by its derivative at time $s$. In particular $\norm{\mu_s}_2=1$ for any $s\in[0, d_{WP}(X_0, X_1)]$. Then, by integrating on $\alpha(s)$ \[\abs{\widetilde{V_R}(X_1)-\widetilde{V_R}(X_0)}\leq \int_0^{d_{WP}(X_0, X_1)}\abs{ \langle (d\widetilde{V_R})_{\alpha(s)}, \mu_s\rangle} ds\leq \int_0^{d_{WP}(X_0, X_1)} \norm{d\widetilde{V_R}}_{\infty}\int \abs{\mu_s}da_{\alpha(s)}ds~.\]
	By Cauchy-Schwartz inequality, denoting by $da_{\alpha(s)}$ the hyperbolic area form on $\alpha(s)$, that we thought as $\rho_{\alpha(s)}^{1/2}\rho_{\alpha(s)}^{1/2}\abs{dz^2}$ \[\int \abs{\mu_s}da_{\alpha(s)}\leq \norm{\mu_s}_2 \sqrt{2\pi\chi(\alpha(s))}=\sqrt{2\pi\chi(\partial\bar M)}~.\]
	 Theorem \ref{boundedGradientinfty} now furnishes the constant $C(\partial\bar M )$ bounding $\norm{d\widetilde{V_R}}_{\infty}
	 $, and the statement follows putting together the inequalities \[\abs{\widetilde{V_R}(X_1)-\widetilde{V_R}(X_0)}\leq \sqrt{2\pi\abs{\chi(\partial\bar M)}}C(\partial\bar M)d_{WP}(X_0, X_1)~.\]
	
	\epf
	
	\bcor
	Let $M$ be a convex co-compact hyperbolic manifold, and let $\mu$ denote the infinitesimal earthquake or the infinitesimal grafting associated to a compressible simple closed curve $\gamma\subseteq X$, with $X\in\mathcal{T}(\partial\bar M)$. Then \[\abs{(d\widetilde{V_R})_X(\mu)}\leq C(\partial\bar M)\ell~,\]
	where $C(\partial\bar M)$ is a constant depending only on the genus of $X$, and $\ell$ is the length of $\gamma$ with respect to the hyperbolic metric in the conformal class of $X$.
	\ecor
	
	\bpf
	By Theorem \ref{boundedGradient}, there exists $C(\partial \bar M)$ such that, at any $X$ \[\norm{d\widetilde{V_R}}_{\infty}\leq C(\partial \bar M)~.\]
	The statement now follows by applying Lemma $5.1$ in \cite{CGS2024}: \[\abs{(d\widetilde{V_R})_X(\mu)}\leq \int_0^{\ell} \abs{(d\widetilde{V_R})_X(\dot\gamma, \dot\gamma)}dt\leq \norm{d\widetilde{V_R}}_{\infty}\int_0^{\ell}\norm{\dot\gamma}^2dt~,\] for $\gamma$ a unit length parameterization of the curve where we are grafting, and $\norm{\dot\gamma}$ the hyperbolic norm of its derivative. The proof for the infinitesimal earthquake is the same, up to substitute a $\dot\gamma$ with $i\dot\gamma$.
	\epf
	
	\bthm\label{Vcdistance}
	Let $M$ be a convex co-compact hyperbolic manifold. Then, there exists a constant $A(\partial\bar M)$ depending only on the topology of $\partial\bar M$, such that, for every $X\in\mathcal{T}(\partial\bar M)$ \[\abs{\widetilde{V_R}(X)-V_C(X)}\leq A(\partial\bar M)~,\]
	where $V_C(\cdot)=\text{Vol}(C(M(\cdot)))$ denotes the convex core volume function. 
	\ethm
	
	\bpf
	By Theorem $A.6$ in \cite{averages}, there exists a constant $A_0(\partial\bar M)$ depending only on the topology of $\partial\bar M$ such that \[\abs{V_R(X)-W(C(X))}\leq A_0(\partial\bar M)~,\]
	where $W(C(\cdot))$ denotes the $W$-volume of the convex core of $M(X)$ (see \eqref{WCC}). Therefore, for any $X\in \mathcal{T}(\partial\bar M)$ \[\abs{\widetilde{V_R}(X)-V_C(X)}\leq \abs{V_R(X)-W(C(X))}+\abs{ \max_{m\in\Gamma^{comp}(\partial\bar M)}L(X, m)-\f14 L(\beta_X)}~,\]
	where $L(\beta_X)$ denotes the length of the bending lamination of the convex core of $M(X)$. It then remains to uniformly bound the last term in the inequality above.   For any fixed $\delta>0$, both the functions $\max_{m\in\Gamma^{comp}(\partial\bar M)}L(X, m)$ and $1/4L(\beta_X)$ are bounded by a uniform constant on the $\delta$-compressible thick part of the Teichm\"uller space. Let us fix $\delta=\varepsilon_1$, with $\varepsilon_1$ as in Lemma \ref{var1}. Then, we have to bound \[\abs{\sum_{\substack{\gamma \text{\ compressible\ } \\ \ell_{\gamma}(X)<\varepsilon_1}}  \f{\pi^3}{\ell_{\gamma}(X)}-\f14 L(\beta_X)}~,\] when $X$ lies in the $\varepsilon_1$-compressible thin part.
	First, let us consider the counter image $\pi^{-1}(\widetilde{X_{\varepsilon_1}^c})\subseteq C(M(X))$ of the (lift of) the complement of the union of the compressible $\varepsilon_1$-thin tubes, through the normal projection $\pi\colon C(M(X))\rightarrow \Omega(\Gamma)$, with $M(X)=\mathbb{H}^3/\Gamma$. The length of the intersection between the bending lamination $\beta_X$ and the projection of the $\varepsilon_1$-compressible thick part of $X$ is bounded thanks to Theorem $2.16$ in \cite{BBB2018} as \[L\left(\beta_X\cap \pi^{-1}(\widetilde{X_{\varepsilon_1}^c})\right)\leq 6\pi\chi(\partial\bar M)\coth^2(\varepsilon_1/4)~.\] Let us now look at the union of the thin tubes $\mathcal{A}(\gamma)$, for each compressible simple closed curve $\gamma$ in $X$ of length less than $\varepsilon_1$. We denote by $T(\gamma)$ their counter image on the convex core through $\pi$.
	By Lemma $A.16$ in \cite{averages}, there exists a universal constant $A_1$ such that, for any compressible simple closed curve $\gamma$ of length less than $\varepsilon_1$, if $b$ is the length of a component of $\beta_X\cap T(\gamma)$, then all the other components of $\beta_X\cap T(\gamma)$ have length in $[b-A_1, b+A_1]$, and \[\abs{\text{Mod}(\gamma)-\f{b}{2\pi}}\leq A_1~,\]
	where the $\text{Mod}(\cdot)$ stands for the modulus of $\gamma$ in the conformal structure of $X$.
	Moreover, for any $\gamma$ as above of length $\ell\leq\varepsilon_1$, by the work of Maskit in \cite{Maskit1985}, there exists a constant $A_2$ depending only on $\varepsilon_1$ such that \[\abs{\text{Mod}(\gamma)-\f{\pi}{\ell}}\leq A_2~.\]
	Therefore \[\abs{\f{\pi^3}{\ell}-\f{b\pi}{2}}\leq \pi^2(A_1+A_2)~.\]
	To estimate the length of the measure bending lamination in the thin tubes, it remains to estimate the intersection number between $\beta_X$ and $\gamma$, that is, the total bending angle along $\gamma$. Consider a totally geodesic disk $D$ with boundary in $T(\gamma)$ intersecting $\gamma$ and lying in a hyperplane orthogonal to a component of $\beta_X$. This last condition implies that the length of $\partial D$ is comparable to the one of $\gamma$ in $\partial C(X)$, which we denote with $\ell_{\gamma}(\partial C(X))$. Then, a standard argument exploiting Gauss-bonnet Theorem and the Isoperimetric Inequality on the disk bounded by $\gamma$ and on $D$ shows that (see \cite{BO}) \[2\pi\leq i(\gamma, \beta_X)\leq 2\pi+2\ell_{\gamma}(\partial C(X))~.\]
	Therefore \[\f{\pi (b-A_1)}{2}\leq \f{1}{4}L(\beta_X\cap T(\gamma))\leq \f{\pi (b+A_1)}{2}+\f{(b+A_1)\ell_{\gamma}(\partial C(X))}{2}~.\]
	The boundary of the convex core is equipped with an induced hyperbolic metric \cite[Chapter 2.1]{CEM2006}. Thanks to the explicit geometry of hyperbolic thin tubes (see Section \ref{margulistube}), one can easily verify that the quantity $b\ell_{\gamma}(\partial C(X))$ is uniformly bounded, so that there exists a constant $A_3$ depending only on $\varepsilon_1$ such that \[\abs{ \f{1}{4}L(\beta_X\cap T(\gamma))-\f{b\pi}{2}}\leq A_3~.\]
	 Merging the inequalities, for any simple closed compressible geodesic $\gamma$ in $X$ of length $\ell\leq\varepsilon_1$ \[\abs{\f{\pi^3}{\ell}-\f14 L(\beta_X\cap T(\gamma))}\leq \pi^2(A_1+A_2)+A_3~.\] Since there are at most $3g-3$ such curves, for $g$ the sum of the geni of the boundary components of $\partial\bar M$, and $L(\beta_X)$ is uniformly bounded outside the union the thin tubes $T(\gamma)$, the statement follows.
	\epf
	
	\bcor\label{Vcestimate}
 Let $M$ be a convex co-compact genus $g$ handlebody, and let $X\in\mathcal{T}(\partial\bar M)$. Let $\Gamma_0^{comp}(\partial \bar M)$ denote the set of compressible multicurves on $\partial \bar M$ that decompose $M$ into a union of punctured tori. Then, for any $\ell>0$ small enough, the following inequality holds \[V_C(M(X))\leq K(g)d_{WP}(X, H_{\ell})+B(g, \ell)~,\] where $K$ and $B$ are constants depending, respectively, linearly on $g$ and on $(g, \ell, g\ell)$, and $H_{\ell}$ is the following subset of the Teichm\"uller space \[ H_{\ell}=\{X\in \mathcal{T}(\partial\bar M)\ |\  \exists m\in\Gamma_0^{comp}(\partial\bar M) \text{ s.t. } \ell_{\gamma}(X)\leq \ell\ \forall \gamma\in m\}~.\] 
	\ecor
	
	\bpf
 Let $X_0\in \bar{H_{\ell}}^{\scriptscriptstyle{WP}}$ be a point realizing the Weil-Petersson distance between $X$ and $X_0$. Then, by Corollary \ref{WPdistance} \[\abs{\widetilde{V_R}(X)-\widetilde{V_R}(X_0)}\leq K(g)d_{WP}(X, X_0)~,\] with \[K(g)=\sqrt{2\pi(2g-2)}C(g)~,\] where $C(g)$ is furnished by Theorem \ref{boundedGradientinfty}. Moreover, by the inequality above and by Theorem \ref{Vcdistance} \[\abs{V_C(X)-V_C(X_0)}\leq 2A(g)+K(g)d_{WP}(X, X_0)~,\] for $A(g)=A(S_g)$, with $S_g$ the closed genus $g$ surface. Now, Vargas-Pallete proved in \cite[Equation $21$]{VP2019} that for any $X_0\in H_{\ell}$ the convex core volume of the corresponding handlebody is bounded as \[V_C(X_0)\leq (g-1)\ell_{Th}~,\] where $\ell_{Th}$ denotes a bound on the lengths for the Thurston metric of the curves $\gamma$ of length at most $\ell$ in the hyperbolic metric at infinity. Since the Thurston metric \cite[Section $4.3$]{dumas-survey} is obtained from the induced one on the convex core by grafting on the bending lamination, and we saw in the proof of Theorem \ref{Vcdistance} that the total intersection angle is bounded by $2\pi+2\ell_{\gamma}(\partial C(X_0))\leq 2\pi+2\ell$, we obtain \[V_C(X_0)\leq (g-1)(2\pi+2\ell)~,\]
	concluding the proof by imposing \[B(g, \ell)=2A(g)+(g-1)(2\pi+2\ell)~.\]

	\epf
	
	Fixed a pants decomposition $P$ of $\partial\bar M$ and a small $\delta>0$, a $\delta$-\textit{Bers region} is the subset $\mathcal{B}_{\delta}$ of $\mathcal{T}(\partial \bar M) $ such that \[\mathcal{B}_{\delta}=\{X\in \mathcal{T}(\partial \bar M)\ |\ \ell_{\gamma}(X)\leq\delta \ \ \forall \gamma\in P \}~.\]
	In the proof of the next theorem we use \textit{Fenchel-Nielsen coordiantes} for the Teichn\"uller space (see, for example \cite[Sec 6.2]{Bu1992}, \cite[Sec 10.6]{FM2011}, or \cite[Section $7.3$]{Mar16}).
	
	\bthm\label{BersRegion}
	Let $P$ be a compressible pants decomposition of the boundary $\partial\bar M$ of a convex co-compact hyperbolic $3$-manifold. For any $\varepsilon\leq\varepsilon_1$ there exists a $\delta>0$ such that for any $X_0$ and $X_1$ in the $\delta$-Bers region of $P$ \[\abs{\widetilde{V_R}(X_1)-\widetilde{V_R}(X_0)}\leq \varepsilon~.\]
	
	\ethm
	
	\bpf
	Let $X_0$ and $X_1$ be two points of $\mathcal{B}_{\delta}$. We can modify the Fenchel-Nielsen coordinates with respect to $P$ to move from $X_0$ to $X_1$ with a complex earthquake path on $P$, i.e. first we modify the twist parameter by earth-quaking on each component of $P$, and then the lengths by grafting (see \cite{McMcomplexearth}).
	% Thanks to Lemma \ref{intbypart2} (Lemma $2.2$ in \cite{CG}), one can see how the pairing between a holomorphic quadratic differential and an infinitesimal earthquake $\nu$ or the an infinitesimal grafting $\mu$ on a simple closed curve $\gamma$ reduces to an integral on $\gamma$:  
	% Thanks to Theorem \ref{boundedGradientinfty}, by integration, a complex earthquake on $\gamma$ gives a bound on the difference of renormalized volumes linear in $\delta$ (here, integrating along the grafting path we have to take care of the variation of the length of $\gamma$, which still is bounded linearly in $\delta$, see for example Lemma \ref{graftlength} in Chapter \ref{3}).
	Let us consider $\gamma\in P$, we start by changing the twist parameter by a $t$-earthquake, reaching the surface $X_t$ whose Fenchel-Nielsen coordinates differ from the one of $X_0$ just by the twist of $\gamma$, which coincides with the one of $X_1$. In this way the standard renormalized volume changes, by \cite[Theorem $1.3$]{CG}, as	\[\abs{V_R(X_t)-V_R(X_0)}\leq 
	F(\delta)t~,\]
	with $F(\delta)$ such that \[|F(\delta)|\leq C\f{e^{-\pi^2/\delta}}{\delta}\leq C~,\] for some explicit constant $C>0$. Since, by uniformization theorem of the deformation space of convex co-compact structures, the renormalized volume descends to a function on the quotient of $\mathcal{T}(\partial\bar M)$ by the group generated by Dehn twists on compressible simple closed geodesic, we can assume $t\leq\delta$ (see Remark \ref{welldef}). Moreover, since the length of any short compressible simple closed curve remains unchanged \[\abs{\widetilde{V_R}(X_t)-\widetilde{V_R}(X_0)}=\abs{V_R(X_t)-V_R(X_0)}~.\]
	We repeat this on any $\gamma\in P$, obtaining the surface $X_{1/2}$.
	Finally, when changing the length of any $\gamma_i\in P$, for $i=1,\dots,k$, by a parameter $s_i$ grafting on it, By \cite[Theorem $1.4$]{CG}, the difference of the adapted renormalized volumes is bounded as \[\abs{\widetilde{V_R}(X_1)-\widetilde{V_R}(X_{1/2})}=\f\pi4\sum_{i=1}^{k}\left[\abs{(\ell_{\gamma_i}(X_1)-\ell_{\gamma_i}(X_{1/2}))}+O\left(e^{-\pi s_i/(2\delta)}s_i^3\right)\right]\leq \f{k\pi\delta}{4}+\left(\f{6\delta}{e\pi}\right)^3~.\]
	We proved that the total variation of the adapted renormalized volume is bounded by a linear function of $\delta$, therefore the statement follows.
	\epf
	
	Note that the hypothesis of the theorem above imply that $M$ is an handlebody.

	\subsection{Renormalized volume of a long tube}\label{VRtube}
	In this section, we are going to study how the renormalized volume diverges on long tubes. This will be useful in the next section, where we study the sequence of adapted renormalized volumes under the pinching of a compressible (multi)curve, in order to prove its continuous extension to the boundary of the Teichm\"uller space. We point out that \cite[Theorem $1.4$]{CG} already shows the asymptotic behaviour of the renormalized volume under the pinching of a compressible curve. Moreover, by looking at the proof, we see that the divergent contribution comes from the long tubes, as it is precisely there that the differential of the renormalized volume cannot be uniformly bounded. However, here we will need a more explicit control on the divergence, which also takes into account \emph{where} the tube has been truncated. More precisely, we aim to express its $W$-volume, as a domain with boundary equipped with a hyperbolic metric, in terms of the length of the boundary curves, other than the length of its core curve.
	
	We start by studying a class of examples in which the developing map of the tube and the metric are explicit. In the subsection that follows the next, we will generalize the result, analyzing the lower order terms that arise.
	\subsubsection{Toy model}\label{Toy}
	Let us consider the map $f_{\ell}\colon \mathbb{H}^2\rightarrow \mathbb{C}$ from the half-space model of the $2$-hyperbolic space, with coordinate $z=\rho e^{i\theta}$, where $\rho>0$ and $\theta \in [0, \pi]$, defined by \[f_{\ell}(z)=z^{\f{2\pi i}{\ell}}=e^{-\f{2\pi\theta}{\ell}+i\f{2\pi}{\ell}\log(\rho)}~.\] The restriction $f_{\ell,\varepsilon}$ of $f_{\ell}$ to the neighborhood of the vertical axis given by the cone of angle $\pi-2\arcsin(\ell/\varepsilon)$ is the developing map of the symmetric complex projective annulus \begin{equation}\label{Aell}\mathcal{A}_{\ell}(\varepsilon)=\left\{\rho e^{i\theta}\in\mathbb{C}\ \ |\ \ e^{-\f{2\pi^2}{\ell}+\f{2\pi}{\ell}\arcsin{\left(\f{\ell}{\varepsilon}\right)}} \leq \rho \leq e^{-\f{2\pi}{\ell}\arcsin{\left(\f{\ell}{\varepsilon}\right)}}\right\}~.\end{equation} By \textit{symmetric} here we mean that $\mathcal{A}_{\ell}(\varepsilon)$ has concentric round boundary components. Note that, equipping $\mathcal{A}_{\ell}(\varepsilon)$ with the restriction of the flat metric $h_{Eu}= \f 1{r^2}|dz|^2$ on the infinite tube $\mathbb{C}\setminus \{0\}$, one obtains a truncated euclidean tube with core length $2\pi$ and modulus such that the pair $(\mathcal{A}_{\ell}(\varepsilon), h_{Eu})$ is conformally equivalent to a truncated hyperbolic tube with core length $\ell$ and boundaries of length $\varepsilon$. Moreover, up to M\"obius transformations, any symmetric complex projective tube is realized as some $\mathcal{A}_{\ell}(\varepsilon)$ for suitable $\ell>0$ and $\varepsilon>0$, and any tube is conformal to a symmetric one (as any modulus is achievable, and this is a complete conformal invariant).
	
	The hyperbolic metric $\dfrac{1}{t^2}(dx^2+dt^2)$ on $\mathbb{H}^2=\mathbb{R}\times \mathbb{R}^{+}$, with coordinates $(x,t)$, can be pushed forward via $f_{\ell}$, obtaining the following metric tensor on $\mathcal{A}_{\ell}(\varepsilon)$:  
	
	\begin{equation}\label{Iell}\hat{I_\ell}(z)=\dfrac{\ell^2}{4\pi^2\rho^2\sin^2\left(\f{\ell}{2\pi}\log(\rho)\right)}\abs{dz}^2\end{equation}
	at $z=\rho e^{i\theta}$.
	
	\brem
	For any fixed $0 < \rho_1 < \rho_2 < 1$, the metric $\hat{I_\ell}$ converges uniformly to $\hat{I_0}$ (as in (\ref{I})) on the compact $\mathbb{D}_{\rho_1}^{\rho_2}$ (as in (\ref{D12})) as $\ell\rightarrow 0$.
	\erem

	\brem\label{Sfl}
	The Schwarzian derivative of $f_{\ell}$ is:
	\begin{align*}
		q=\mathcal S(f_{\ell})&=\biggl(\biggl(\dfrac{(f_{\ell})^{''}}{(f_{\ell})^{'}}\biggr)^{'}-\dfrac{1}{2}\biggl(\dfrac{(f_{\ell}){''}}{(f_{\ell})^{'}}\biggr)^{2}\biggr)dz^2 \\
		&=\left(\left(\left(\f{2\pi i}{\ell}-1\right)\f1z\right)^{'}-\f12  \left(\left(\f{2\pi i}{\ell}-1\right)\f1z\right)^{2}\right)dz^2 \\
		&= \f{1}{2z^2}\left(1+\f{4\pi^2}{\ell^2}\right)dz^2 ~.
	\end{align*}
	\erem

	\blem\label{double}
	Let $(\mathcal{C}_{\ell}(\varepsilon), \hat{I_\ell})$ be half of the symmetric tube $(\mathcal{A}_{\ell}(\varepsilon), \hat{I_\ell})$, with boundaries given by the core $\gamma$ of $\mathcal{A}_{\ell}(\varepsilon)$ and one of its boundary components. Then \[W(\mathcal{A}_{\ell}(\varepsilon), \hat{I_\ell})=2W(\mathcal{C}_{\ell}(\varepsilon), \hat{I_\ell})+2b(\gamma, \hat{I_\ell} )~,\]
	where $b(\gamma, \cdot)$ denotes the integral of the mean curvature terms arising from the boundary component $\gamma$:  \[b(\gamma, \hat{I_{\ell}})=\f12\int_{C(\gamma)}Hda_{C(\gamma)}+\f{\pi}{8}\ell(\partial_1C(\gamma))~,\]
	where $C(\gamma)$ is the associated caterpillar region, and $\ell(\partial_1C(\gamma))$ is the hyperbolic length of the lower boundary component of $C(\gamma)$.
	\elem
	\begin{proof}
		First, note that the core of $\mathcal{A}_{\ell}(\varepsilon)$, being the image under $f_{\ell}$ of the vertical geodesic between $0$ and $\infty$ of $\mathbb{H}^2$, corresponds to the circle of radius $e^{-\pi^2/\ell}$ in $\mathbb{C}$. We consider the hyperplane of boundary at infinity coinciding with the core, and the inversion $r(z)=(\bar z)^{-1}e^{-2\pi^2/\ell}$ along it, which is a M\"obius anti-transformation and thus an orientation reversing isometry of $\mathbb{H}^3$. In particular, r reverses the induced orientation on $\gamma$, so that \[W(\mathcal{C}_{\ell}(\varepsilon), \hat{I_\ell})=W\left(r(\mathcal{C}_{\ell}(\varepsilon)), r^*\left(\hat{I_\ell}\right)\right)-2b(\gamma, \hat{I_{\ell}})~.\] The thesis now follows, by additivity of the $W$-volume (see Lemma \ref{add}), just by observing that \[(\mathcal{A}_{\ell}(\varepsilon), \hat{I_\ell})=(\mathcal{C}_{\ell}(\varepsilon), \hat{I_\ell}) \cup (r((\mathcal{C}_{\ell}(\varepsilon)), r^*(\hat{I_\ell}))~.\]
		
	\end{proof}
	
	Let us denote by $N_{\ell}(\varepsilon)$ the compact subset of $\mathbb{H}^3$ bounded by the Epstein surface associated to the domain with boundary $(\mathcal{A}_{\ell}(\varepsilon),\hat{I_\ell})$ capped with the hyperplanes whose boundary union coincides with $\partial\mathcal{A}_{\ell}(\varepsilon)$ (as in Section \ref{BoundaryEps}).
	
	\begin{figure}
		\begin{center}
			\includegraphics[width=.8\linewidth]{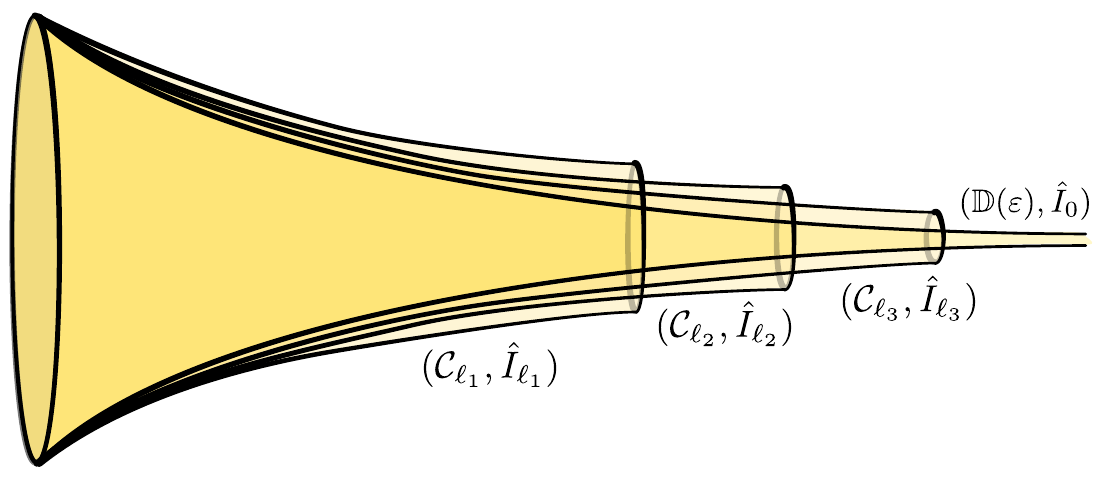}
			\medskip
			\caption{Approximating a cusp with long truncated half-tubes of core lengths $\ell_1 > \ell_2 > \ell_3$. }\label{coni}
		\end{center}
	\end{figure}

	\bprop\label{Wtube}
	The $W$-volume of $N_{\ell}(\varepsilon)$, for any fixed $\varepsilon>\ell$ small enough, has the following behavior at $\ell\sim 0$: \[W(N_{\ell}(\varepsilon))=W(\mathcal{A}_{\ell}(\varepsilon),\hat{I_\ell})= -\f{\pi^3}{\ell}+\f{2\pi^2}{\varepsilon}+2b(\varepsilon)+ O(\ell, \varepsilon)~,\] with \[b(\varepsilon)=\f{\pi^3}{4\varepsilon}+\f{\pi^2}{\varepsilon}+O(\varepsilon)~,\] and where $O(\ell, \varepsilon)$ is a function of $\ell$ and $\varepsilon$ such that the limit $\lim_{\ell\rightarrow 0}\f{O(\ell, \varepsilon)}{\ell}$ is finite.
	\eprop

	\begin{proof}
		We are going to study the $W$-volume of half of the tube  $(\mathcal{C}_{\ell}(\varepsilon), \hat{I_\ell})$, and then use Lemma \ref{double}.
		First observe that, without loss of generality, \[\mathcal{C}_{\ell}(\varepsilon)= \mathbb{D}_{\rho_1(\ell, \varepsilon)}^{\rho_2(\ell, \varepsilon)}\ , \quad \rho_1(\ell, \varepsilon)=e^{-\f{\pi^2}{\ell}} \quad \text{\ and\ } \quad \rho_2(\ell, \varepsilon)=e^{-\f{2\pi}{\ell}\arcsin{\left(\f{\ell}{\varepsilon}\right)}} ~.\]
		Then, thanks to Proposition \ref{CuspConv}: \[W(\mathcal{C}_{\ell}(\varepsilon), \hat{I_0})= -\f{\pi}{2}\log\left(\f{\rho_2}{\rho_1}\right)-b(\rho_1)+b(\rho_2)=-\f{\pi^3}{2\ell}+\f{\pi^2}{\varepsilon}-b(\rho_1)+b(\rho_2)~,\] where $b(\rho_i)$ is the boundary term given by the integral of the mean curvature for the caterpillar region, and \[b(\rho_i)=-\f{\pi^2}{8}\log(\rho_i)-\f{\pi}{2}\log(\rho_i)+O(1/\log(\rho_i))~.\]
		In particular \begin{equation}\label{brho2}
			b(\rho_2)=b(\varepsilon)=\f{\pi^3}{4\varepsilon}+\f{\pi^2}{\varepsilon}+O(\varepsilon)~.
		\end{equation}
		We now apply the Polyakov formula for the $W$-volume associated to domain with round boundaries (see Theorem \ref{Poly}), in order to estimate the difference between $W(\mathcal{
			C}_{\ell}(\varepsilon),\hat{I_\ell})$ and $W(\mathcal{C}_{\ell}(\varepsilon), \hat{I_0})$ (see Figure \ref{coni}). To this end, let us denote by $w_{\ell}$ the smooth function on $\mathcal{C}_{\ell}(\varepsilon)$ such that \[\hat{I_\ell}=e^{2w_{\ell}}\hat{I_0}~.\] We need to bound $w_{\ell}$ and its derivatives: denoting $\mathcal{C}_{\ell}(\varepsilon)$ simply by $\mathcal{C}$, Theorem \ref{Poly} states: 
		
		\[W(\mathcal{C}_{\ell}(\varepsilon), e^{w_{\ell}}\hat{I_0})-W(\mathcal{C}_{\ell}(\varepsilon), \hat{I_0}) =\] \[=-\f{1}{4}\int_\mathcal{C}\left(\abs{\nabla_{\hat{I_0}}w_{\ell}}^2+K(\hat{I_0})w_{\ell}\right)da(\hat{I_0})-\f12 \int_{\partial\mathcal{C}} k(\hat{I_0})w_{\ell}ds(\hat{I_0})\]	
		where $\nabla_{\hat{I_0}}$, $K(\hat{I_0})$, $k(\hat{I_0})$, $da(\hat{I_0})$ and $ds(\hat{I_0})$ denote, respectively, the gradient, scalar curvature, geodesic curvature of the boundary, area form, and length form, all with respect to the metric tensor $\hat{I_0}$. In particular, being $\hat{I_0}$ a hyperbolic cusp, it has scalar curvature $K(\hat{I_0})=-1$. Moreover, the round circles of the boundary of $\mathcal{C}$ have $\hat{I_0}$-geodesic curvature $k(\hat{I_0})=1$, as they are horocycles. To lighten up the notation, we define the quantity \[a(\ell)=\f{\ell}{2\pi}~,\] and we observe that $w_{\ell}$ can be made explicit, at $z=\rho e^{i\theta}$, as \[w_{\ell}(z)=w_{\ell}(\rho)=\log\left(\f{a(\ell)\log(\rho)}{\sin\left(a(\ell)\log(\rho)\right)}\right)~.\]
		Note that $w_{\ell}(\rho_2)=O(\ell^2)$, and $w_{\ell}(\rho_1)=\log(\pi/2)$ is bounded for any $\ell$. 
		Fixed $\varepsilon>0$ small enough, we can expand $w_{\ell}(\rho)$ and its derivative $(w_{\ell}(\rho))'$ near $\ell\sim 0$, getting \begin{align*}& w_{\ell}(\rho) \sim a^2(\ell)\log^2(\rho)/6~ \\& (w_{\ell}(\rho))'\sim \f{a^2(\ell)\log(\rho)}{3\rho} \end{align*}
		First, we study the boundary term. We parameterize the two connected components of $\mathcal{C}$ as $\alpha_i(s)=\rho_ie^{si/\rho_i}$, with $s\in[0, 2\pi\rho_i]$, for $i=1,2$. Then $||\dot{\alpha_i}(s)||_{\hat{I_0}}=1/(\rho_i\log(\rho_i))$ and \[-\f12 \int_{\partial\mathcal{C}} k(\hat{I_0})w_{\ell}ds(\hat{I_0})=\f12 \int_{0}^{2\pi\rho_1}\f{w_{\ell}(\rho_1)}{\rho_1\log(\rho_1)}ds-\f12 \int_{0}^{2\pi\rho_2}\f{w_{\ell}(\rho_2)}{\rho_2\log(\rho_2)}ds=\pi\left(\f{w_{\ell}(\rho_1)}{\log(\rho_1)}-\f{w_{\ell}(\rho_2)}{\log(\rho_2)}\right)\]
		but now, since $w(\rho_2)=O(\ell^2)$ and $w_{\ell}(\rho_1)/\log(\rho_1)=\ell\log(\pi/2)/\pi^2=O(\ell)$, this is $O(\ell)$. 
		
		%The second boundary term is similar, but this time we have to take care of the $\hat{I_0}$-unit normal derivative of $w(\rho)$: \[\partial_n w_{\ell}(\rho)= \rho\log(\rho)\f{\partial w_{\ell}}{\partial\rho}~,\] 
		%and  \[\f{\partial w_{\ell}}{\partial\rho}\sim a^3(\ell)\log^2(\rho)/\rho\]
		
		%therefore \[\f38\int_{\partial\mathcal{A}}\partial_nw_{\ell}ds(\hat{I_0})=2\pi\left(\rho_2\f{\partial w_{\ell}}{\partial\rho}-\rho_1\f{\partial w_{\ell}}{\partial\rho}\right)\] is again of the type $O(\ell, \varepsilon)$.
		
		Let us move on to analyzing the integrals over $\mathcal{C}$. We have \[\f{1}{4}\int_\mathcal{A} K(\hat{I_0})w_{\ell}da(\hat{I_0})=-\f{\pi}{2}\int_{\rho_1}^{\rho_2}\f{w_{\ell}(\rho)}{\rho\log^2(\rho)}d\rho \sim -\f{\pi a^2(\ell)}{12}\int_{\rho_1}^{\rho_2} \f{1}{\rho} d\rho~,\] and \[-\f{\pi a^2(\ell)}{12}\int_{\rho_1}^{\rho_2} \f{1}{\rho} d\rho=-\f{\pi a^2(\ell)}{12}\log\left(\f{\rho_2}{\rho_1}\right)\sim\f{\pi}{12}\left(\f{\ell}{2\pi}\right)^2\left(-\f{2\pi}{\varepsilon}+\f{\pi^2}{\ell}\right)\]
		which is of order $O(\ell)$ for any fixed $\varepsilon>\ell$.
		Concerning the term in the gradient, since $w_{\ell}$ depends just on the radial coordinate \[\abs{\nabla_{\hat{I_0}}w_{\ell}}^2da(\hat{I_0})=(\rho\log(\rho))^2 \left(\f{\partial w_{\ell}}{\partial\rho}\right)^2\f{d\rho d\theta}{\rho\log^2(\rho)}\sim \left(\rho\log(\rho)\f{a^2(\ell)\log(\rho)}{3\rho}\right)^2\f{d\rho d\theta}{\rho\log^2(\rho)}~,\] 
		so that \[\int_{\mathcal{C}}\abs{\nabla_{\hat{I_0}}w_{\ell}}^2da(\hat{I_0})\sim \f{2\pi a^4(\ell)}{9}\int_{\rho_1}^{\rho_2}\f{\log^2(\rho)}{\rho}d\rho=\f{2\pi a^4(\ell)}{27}\left(-\f{2\pi}{\varepsilon}+\f{\pi^2}{\ell}\right)^3~,\]
		which again is of order $O(\ell)$ for any fixed $\varepsilon>\ell$.
		Putting everything together, we conclude: \[W(\mathcal{C}_{\ell}(\varepsilon),\hat{I_\ell})=-\f{\pi^3}{2\ell}+\f{\pi^2}{\varepsilon}-b(\rho_1)+b(\rho_2)+O(\varepsilon, \ell)~.\]
		Now, by Lemma \ref{double} \[W(\mathcal{A}_{\ell}(\varepsilon), \hat{I_\ell})=2W(\mathcal{C}_{\ell}(\varepsilon), \hat{I_\ell})+2b(\alpha_1, \hat{I_\ell} )~,\] where
		\[b(\alpha_1, \hat{I_\ell} )=b(\rho_1)\]
		as \[\alpha_1(s)=\rho_1e^{si/\rho_1}\]
		for $s\in[0, 2\pi\rho_1]$ and $\rho_1=e^{-\pi^2/\ell}$. %At the end of the proof of Proposition \ref{CuspConv}, we already calculated the integral of the mean curvature term coming from the boundary: \[b(\rho_1)= - \f{\pi}{4}\log(\rho_1)+\f{\pi}{12\log(\rho_1)}=\f{\pi^3}{4\ell}+O(\ell)~.\]
		Therefore, thanks also to Equation \eqref{brho2}, the statement follows \[W(\mathcal{A}_{\ell}(\varepsilon), \hat{I_\ell})=-\f{\pi^3}{\ell}+\f{2\pi^2}{\varepsilon}+2b(\varepsilon)+O(\varepsilon, \ell)~.\]
	\end{proof}

	\subsubsection{General case}
	Let us now consider an hyperbolic domain $\Omega\subseteq\mathbb{C}$ containing a short geodesic $\gamma$ of length $\ell\leq\varepsilon_0$ with respect to the hyperbolic conformal metric $h$ on $\Omega$. For any $\varepsilon>0$ small enough, let then $\mathcal{A}_{h}(\varepsilon)$ be a tubular neighborhood of $\gamma$ obtained by cutting $\Omega$ with two hyperplanes whose boundaries at infinity have length $\varepsilon$ with respect to $h$. 
	\bthm\label{Wtubegeneral}
	The $W$-volume associated to $(\mathcal{A}_{h}(\varepsilon), h)$, with $\mathcal{A}_{h}(\varepsilon)$ the tubular neighborhood of a short simple closed curve in a domain $\Omega\subseteq\mathbb{C}$ of length $\ell$ with respect to the conformal hyperbolic metric $h$ of $\Omega$, with round boundary components of length a small enough $\varepsilon$,  has the following behaviour: \[W(\mathcal{A}_{h}(\varepsilon), h)= -\f{\pi^3}{\ell}+\f{2\pi^2}{\varepsilon}+2b(\varepsilon)+ O(\ell, \varepsilon)~~,\]
	where $O(\ell, \varepsilon)$ is a function of $\ell$ and $\varepsilon$ such that the limit $\lim_{\ell\rightarrow 0}\f{O(\ell, \varepsilon)}{\ell}$ si finite.
	\ethm
	
	\begin{proof}
		Thanks to \cite[Theorem $1.1$]{CG}, for any $\varepsilon>0$ small enough, we can assume $\mathcal{A}_h(\varepsilon)$ to be contained in a neighborhood of $\gamma$ whose developing map $f$ is defined on some neighborhood of the geodesic through $0$ and $\infty$ in $\mathbb{H}^2$, and such that \[\mathcal S(f)=\f{1}{2z^2}\left(1+\f{4\pi^2}{\ell^2}\right)dz^2+O\left(\f{e^{-\pi^2/(2\ell)}}{\ell^2}\right)dz^2~.\] 
		Note that, by Remark \ref{Sfl}, the first term in $\mathcal S(f)$ coincides with the Schwarzian derivative of the developing map of the symmetric tube $\mathcal{A}_{\ell}(\varepsilon)$.
		We can then integrate the Schwarzians $\mathcal S(f_{\ell})$ and $\mathcal S(f)$ (see \cite[Section 3.2]{dumas-survey}) to obtain, up to M\"obius transformation, the developing maps of $\mathcal{A}_h(\varepsilon)$ and $\mathcal{A}_{\ell}(\varepsilon)$, respectively. In general, every holomorphic quadratic differential $\psi(z)dz^2$ can be realized as the Schwarzian derivative of a locally injective holomorphic function equal, up to M\"obius maps, to the quotient of a base of the two-dimensional space of solutions of the following ODE (see \cite[Section $6.3$]{Hubbard2016})\[f(z)''+\f12 \psi(z)f(z)=0~.\] Since the solutions of the just above ODE depend continuously on $\psi(z)$, the supremum norm of the difference of the developing maps of $\mathcal{A}_h(\varepsilon)$ and $\mathcal{A}_{\ell}(\varepsilon)$ is bounded by a factor of order $O\left(\f{e^{-\pi^2/(2\ell)}}{\ell^2}\right)$. By continuity of the $W$-volume, the statement follows now from Proposition \ref{Wtube}.
	\end{proof}

	\subsection{Convergence of the adapted renormalized volume}
	
	We consider $M=\mathbb{H}^3/\Gamma$, a convex co-compact hyperbolic $3$-manifold with compressible boundary, and we denote by $h\in [\partial_{\infty}M]$ the unique hyperbolic representative in its conformal boundary at infinity. Let also $\gamma\subseteq \partial_{\infty}M=\Omega(\Gamma)/\Gamma $ be a compressible simple closed curve of length $\ell$ with respect to $h$, and consider its lift to $\Omega(\Gamma)$, which we still denote by $\gamma$. For any small enough $\varepsilon>0$, we furthermore consider the annulus $\mathcal{A}_h(\varepsilon)\subseteq\Omega(\Gamma)$ containing $\gamma$ obtained by cutting the boundary at infinity of $M$ with the boundaries of the two hyperplanes $H_{i}(\varepsilon)$ of length $\varepsilon$ with respect to $h$. Note that, if $h$ coincides with $\hat{I_\ell}$, then also $\mathcal{A}_h( \varepsilon)= \mathcal{A}_{\ell}( \varepsilon)$, with $\hat{I_\ell}$ and $\mathcal{A}_{\ell}( \varepsilon)$ as in Section \ref{Toy}, defined, respectively, by Equations \eqref{Iell} and \eqref{Aell}. As a direct consequence of \cite[Theorem $1.1$]{CG}, when $\ell$ is sufficiently small, the two tubes are almost projectively equivalent, as seen in the proof of Theorem \ref{Wtubegeneral}. We denote by $N_h(\varepsilon)$ the compact region in $M$ bounded by the Epstein surface associated to the domain with boundary $(\mathcal{A}_h( \varepsilon), h)$ and the two hyperplanes $H_i({\varepsilon})$, as described in Section \ref{BoundaryEps}.

	The goal of this section is to show that the function $\widetilde{V_R}$ extends continuously to the strata in the boundary of the Weil-Petersson completion $\overline{\mathcal{T}(\partial\bar{M})}^{\scriptscriptstyle{WP}}$ corresponding to compressible multicurves (see Section \ref{boundarywpclosure}), and that the limit is equal to the sum of the adapted renormalized volumes of the limit convex co-compact manifolds pointed at the pinched components of the multicurve of the corresponding stratum, defined in \ref{AVRpointed}. Recall that a pointed convex co-compact manifold is a pair $(M,P)$, where $P\subseteq \partial_{\infty}M$ is a finite set of points, and that its renormalized volume was defined in \ref{cuspedVR} and Remark \ref{P}.
	
	To this end, we first need to recall the following result on the geometric convergence of convex co-compact hyperbolic manifolds under the pinching of compressible (multi)curves (see Section \ref{Geomconv}).

	\bthm\label{GromovHausd}(\cite[Section $A.10$]{averages})
	Let $M_n = (M, g_n)$ be a sequence of convex co-compact hyperbolic $3$-manifolds, and let $m$ be a compressible multicurve in $\partial\bar{M}$ such that the conformal boundaries $\partial_{\infty}M_n$ converge to a conformal structure $X_{m}$ in the stratum of $\overline{\mathcal{T}(\partial\bar{M})}^{\scriptscriptstyle{WP}}$ corresponding to $m$. Let also $D(m)$ be a union of disks compressing $m$. For each $n$, fix a point $y_i(n)$ in the thick part of every of the $i$ connected components of $ C(M_n)\setminus D(m)$. Then $(M_n, y_i(n))$ converges in the Gromov-Hausdorff topology to a pointed complete hyperbolic manifold $(M_i, P_i)$, whose boundary at infinity is the union of the connected components of $X_m$ facing the one of $y_i(n)$ in $ C(M_n)\setminus D(m)$. Moreover, each limit $M_i$ is either convex co-compact, a hyperbolic solid torus, or $\mathbb{H}^3$.
	\ethm
	
	We remark that, since a hyperbolic manifold is irreducible, the disks compressing a fixed simple closed curve on the boundary are all isotopic.
	
	In the notation of the theorem above, we denote by $h_{P_i}$ the complete hyperbolic metric in the
	conformal class of the connected component(s) of $X_m$ corresponding to $(M_i, P_i)$ that has a cusp at each point $p\in P_i$. We also remark that the manifold $M_i$ is diffeomorphic to the $i$-th connected component of $M\setminus D(m)$.
	
	\brem
	We remark that even in the cases where $M_i$ is diffeomorphic to a solid torus or a $3$-ball, the Euler characteristic of $\partial\bar M_i\setminus P_i$ is negative, and therefore the hyperbolic metric $h_{P_i}$ still exists. The definition of renormalized volume of a pointed convex co-compact manifold extends to these cases.
	\erem
	
	%\brem\label{Remadmissible}
	%From now on, we will always assume that all pointed Gromov-Hausdorff limits are convex co-compact. However, we remark that even in the cases where $M_i$ is diffeomorphic to a solid torus or a $3$-ball, the Euler characteristic of $\partial\bar M_i\setminus P_i$ is negative, and therefore the hyperbolic metric $h_{P_i}$ still exists. The definition of renormalized volume of a pointed convex co-compact manifold seems to be extendable to these cases, as well as the result that follows.
	%\erem
	
%	\bdefi\label{admissible}
%	Let $m\subseteq S$ be a multicurve on a closed surface. We say that $m$ is \textit{admissible} if every component of its complement $S\setminus m$ is a surface of genus $g\geq 2$.
%	\edefi
	
	We now define the adapted version for the renormalized volume of pointed convex co-compact manifolds $(M,P)$ (as in Definition \ref{cuspedVR}) as a function on the strata of the boundary of the Teichm\"uller space. We first fix some notation. Given $M$ a convex cocompact manifold diffeomorphic to the interior of $N$ a tame manifold, and $m\subseteq\partial N$ a compressible multicurve, we denote by $N_i$, with $i=i,\cdots, k$, the connected components of $N\setminus D(m)$. Moreover, recall that a point in the stratum $\mathcal{T}(\partial\bar M\setminus m)$ can be seen as the union $(X,P)$ of pointed Riemann surfaces $(X_i, P_i)$, such that $\sum_i \abs{P_i}=2\abs{m}$, and each of them is homeomorphic to $\partial N_i$. We also denote by $M_i(X_i)$ the convex cocompact manifold diffeomorphic to the interior of $N_i$ and with conformal boundary $X_i$.
	
	%Before that, given the pair $(X,P)$, with $X$ a Riemann surface structure on a closed surface $S$, and $P\subseteq X$ a finite set of points, we define the \textit{Teichm\"uller space pointed at $P$}, and we denote it by $\mathcal{T}_P(S)$, as the image of the natural inclusion \[i_P\colon \mathcal{T}(S)\hookrightarrow \mathcal{T}(S\setminus P)~,\]
	%	which assign to each Riemann surface in $\mathcal{T}(S)$ the unique conformal hyperbolic metric with a cusp singularity at every point of $P$.

	\bdefi\label{AVRpointed}
	Let $M$ be a convex co-compact manifold, let $m\subseteq\partial\bar M$ be a compressible multicurve, and let $\mathcal{T}(\partial\bar M\setminus m)$ be the corresponding stratum in $\overline{\mathcal{T}(\partial\bar{M})}^{wp}$. In the notation above, we define the \textit{adapted renormalized volume on the stratum of $m$} as the function \[\widetilde{V_R}\colon\mathcal{T}(\partial\bar M\setminus m)\rightarrow \mathbb{R}\] such that \[\widetilde{V_R}(X, P)\eqdef\sum_{i=1}^k \widetilde{V_R}(M_i(X_i), P_i)\eqdef\sum_{i=1}^k \left[V_R(M_i(X_i), P_i)+\max_{m_i\in\Gamma^{comp}(\partial\bar M_i)} L(X_i, m_i)\right] \]
 where $m_i$ runs in the set of compressible multicurves in $\partial\bar{M_i}$, and 
 \[ L(X_i, m_i)=\pi^3\sum_{\substack{\gamma\in m_i }} \f{1}{\ell_{\gamma}(X_i)} \]
for $\ell_{\gamma}(\cdot)$ the length function of $\gamma$ on $\mathcal{T}(\partial M_i)$.
	\edefi
	
	A simple closed curve on a surface is called \textit{separating} if its complement is disconnected. The two cusps that arise when a curve is pinched, depending on whether this is separating or not, lie in either two different or the same connected component of the limit Riemann surface $X_m\in\mathcal{T}(\partial\bar M\setminus m)$. In the case of a multicurve one has to be more careful, as a union of non-separating curve can still have a disconnected complement. Moreover, the number of connected components of the complement of the multicurve in the boundary surface may differ from the one of the connected components of $M\setminus D(m)$, for $D(m)$ a union of disks compressing $m$ in $M$. In any case $P_i$ corresponds to the subset of the $2|m|$ cusp singularities in the components of $X_m$ in the boundary of $M_i$. To avoid the technicalities in the notation and indexing due to this phenomenon, we first state and prove the theorem for the pinching of a single separating compressible curve, whose compressing disk also disconnect the $3$-manifold.  
	
	\bthm\label{conv}
	Let $M_t=(M, g_t)$ be a path of convex co-compact hyperbolic $3$-manifolds with connected boundary obtained, via uniformization theorem, by the pinching of a separating compressible simple closed curve $\gamma$ in the conformal boundary at infinity, with compressing disk $D(\gamma)$ that separates $M$. Let also $(M_1, g_1, p_1)$ and $(M_2, g_2, p_2)$ be the two pointed convex cocompact limits of $(M_t, y_{i,t})$ in the Gromov-Hausdorff topology, with $y_{i}(t)$, for $i=1,2$, lying in the thick part of the two different connected components of $ C(M_t)\setminus D(\gamma)$. Then \[\lim_{t\rightarrow \infty} \widetilde{V_R}(M_t)= \widetilde{V_R}(M_1, g_1, p_1)+\widetilde{V_R}(M_2, g_2, p_2)~.\]
	\ethm
	
	%Before proceeding with the proof, we remark that the continuous extension fails precisely at the intersection between the closure of the codimension-one locus in which $\widetilde{V_R}$ is just lower semi-continuous (see Remark \ref{discontinuity}) and the boundary of the Teichm\"uller space.
	
	\begin{proof}
		Since the convex core of any $M_t=\mathbb{H}^3/\Gamma_t$ contains a hyperbolic geodesic $\mathcal{L}_t$ going through the (unique up to isotopy) disk $D(\gamma)$ that compresses $\gamma$, up to composition with a path of  M\"obius transformations $\phi_t$, denoting by $\Omega_t$ the domain of discontinuity of $M_t$, we can assume $\Omega_t\subseteq \mathbb{CP}^1-\{0,\infty\}$ and $\mathcal{L}_t$ to be the geodesic through $0$ and $\infty$. Let us denote by $\ell(t)$ the length of $\gamma$ with respect to the hyperbolic representative $h_t\in[\partial_{\infty}M_t]$ in the conformal class of the boundary at infinity of $M_t$. By hypothesis $\ell(t)\rightarrow 0$.
		
		For any $\varepsilon>\ell(t)>0$ small enough, we also denote by $\mathcal{A}_{h_t}(\varepsilon)\subseteq \Omega_t$ the union of annuli around the lifts of $\gamma$, each of them obtained by cutting the boundary at infinity with two hyperplanes with boundary in $(\Omega_t, h_t)$ of length $\varepsilon$ with respect to $h_t$. With an abuse of notation, we denote by $\mathcal{A}_{h_t}(\varepsilon)$ also the image via the projection to $\partial_{\infty}M_t=\Omega_t/\Gamma_t$. 
		Since by hypothesis $\partial\bar{M}_t$ is connected and $\gamma$ is separating, $\partial_{\infty}M_t-\mathcal{A}_{h_t}(\varepsilon)$ has two connected components which we denote by $(\partial_{\infty}M_t)_i(\varepsilon)$, for $i=1,2$. Let also $\Gamma_{t,i}< \mathbb{P}SL(2,\mathbb{C})$ denote the fundamental groups of the two connected manifolds obtained by cutting along the compressible disk $D(\gamma)$, whose free product gives $\Gamma_t$.
		
		We recall that the conformal boundary at infinity $\partial_{\infty}M_t$ is naturally equipped with a complex projective structure, and so also $(\partial_{\infty}M_t)_i(\varepsilon)$ are equipped with the complex projective structure induced by the restriction. Analogously, we also denote by $\Omega_{i,t}(\varepsilon)$ the components of $\Omega_t- \mathcal{A}_{h_t}(\varepsilon)$ that cover $(\partial_{\infty}M_t)_i(\varepsilon)$. In this way, the Epstein surface associated to $((\partial_{\infty}M_t)_i(\varepsilon), h_t)$, is the quotient of the one associated to $(\Omega_{i,t}(\varepsilon), h_t)$ by the action of the group of M\"obius transformation $\Gamma_{t,i}$.
		Let us fix now a small $\varepsilon>0$. By additivity of the $W$-volume (see Lemma \ref{add}), we can split the renormalized volume of $M_t$ as \[V_R(M_t)=W((\partial_{\infty}M_t)_1(\varepsilon), h_t)+W((\partial_{\infty}M_t)_2(\varepsilon), h_t)+W(\mathcal{A}_{h_t}(\varepsilon), h_t)~.\]
		Then, the adapted renormalized volume of $M_t$ splits as \[\widetilde{V_R}(M_t)=W((\partial_{\infty}M_t)_1(\varepsilon), h_t)+W((\partial_{\infty}M_t)_2(\varepsilon), h_t)+W(\mathcal{A}_{h_t}(\varepsilon), h_t)+\f{\pi^3}{\ell(t)}+L(\partial_{\infty}M_t, m_t\setminus\{\gamma\})~,\]
		where $m_t\in\Gamma^{comp}(\partial \bar M)$ is a multicurve realizing the maximum of the function \[L(\partial_{\infty}M_t, \cdot)=\pi^3\sum_{\substack{\gamma\in \cdot }} \f{1}{\ell_{\gamma}(\partial_{\infty}M_t)}~,\] and
		where every time we consider $h_t$ restricted to the domain we are referring to.
		
		 We denote by $h_{p_1}$ and $h_{p_2}$ the hyperbolic representatives in the conformal boundary at infinity of the convex co-compact manifolds $(M_1, g_1)$ and $(M_2, g_2)$ with a cusp singularity, respectively, at $p_1$ and $p_2$. We also denote by $\Omega_i(\varepsilon)$, for $i=1,2$, the complement in the discontinuity domain of $M_i=\mathbb{H}^3/\Gamma_i$ of $\Gamma_i\cdot D_{p_i}(\varepsilon)$, where $D_{p_i}(\varepsilon)$ is the disk centered at $p_i$ with boundary of length $\varepsilon$ with respect to $h_{p_i}$. Now, observe that the Gromov-Hausdorff convergence implies the one of the domains $\Omega_{i,t}(\varepsilon)\rightarrow \Omega_i(\varepsilon)$.
		% by uniqueness of the hyperbolic representative, also $h_t\rightarrow h_i$ on $\Omega_{i,t}(\varepsilon)$, for $i=1,2$.
		Therefore, by continuity of the $W$-volume \[\lim_{t\rightarrow\infty}W((\partial_{\infty}M_t)_1(\varepsilon), h_t)= W(\partial_{\infty}M_1(\varepsilon), h_{p_1})\] and equivalently
		\[\lim_{t\rightarrow\infty}W((\partial_{\infty}M_t)_2(\varepsilon), h_t)= W(\partial_{\infty}M_2(\varepsilon), h_{p_2})\]
		where again all the metrics are considered restricted to the domains indicated.
		Moreover, since $\ell(t)\rightarrow 0$, by Theorem \ref{Wtubegeneral} \[\lim_{t\rightarrow \infty}W(\mathcal{A}_{h_t}(\varepsilon), h_t)+\f{\pi^3}{\ell(t)}=\f{2\pi^2}{\varepsilon}+2b(\varepsilon)~.\]
		Therefore, for any $\varepsilon>0$ \[\lim_{t\rightarrow \infty}\widetilde{V_R}(M_t)= W(\partial_{\infty}M_1(\varepsilon), h_{p_1})+\f{\pi^2}{\varepsilon}+b(\varepsilon)+ W(\partial_{\infty}M_2(\varepsilon), h_{p_2})+\f{\pi^2}{\varepsilon}+b(\varepsilon)+\lim_{t\rightarrow\infty}\L(\partial_{\infty}M_t, m_t\setminus\{\gamma\})~.\]
		Since $\partial_{\infty}M_t$ converges to the union of $(\partial_{\infty}M_1, p_1)$ and $(\partial_{\infty}M_2, p_2)$, taking the limit for $\varepsilon\rightarrow 0$ of the above equality, by Definition \ref{cuspedVR} and Definition \ref{AVRpointed} \[\lim_{t\rightarrow \infty} \widetilde{V_R}(M_t)= \widetilde{V_R}(M_1, g_1, p_1)+\widetilde{V_R}(M_2, g_2, p_2)~.\]
		
	\end{proof}
	
	The proof of Theorem \ref{conv} similarly applies by iteration to the case of pinching a compressible multicurve, using the definition of renormalized volume for pointed convex co-compact manifolds of remark \ref{P}.

	\bthm\label{AVRContinuity}
	Let $M_t=(M, g_t)$ be a path of convex co-compact hyperbolic $3$-manifolds obtained by pinching a compressible multicurve $m$ in the conformal boundary at infinity. Let $D(m)$ be a union of disks compressing $m$, and let $(M_i, g_i, P_i)$, for $i=1,\dots,k$, be the pointed convex co-compact limits of $(M_t, y_{i}(t))$ in the Gromov-Hausdorff topology, with $ y_{i}(t)$ in the thick part of the $i$-th connected component of $C(M_t)\setminus D(m)$. Then \[\lim_{t\rightarrow\infty} \widetilde{V_R}(M_t)=\sum_{i=1}^k \widetilde{V_R}(M_i, g_i, P_i)~.\]
	\ethm

	%\thispagestyle{empty}
%	\cleardoublepage
%	\phantomsection
%	\addcontentsline{toc}{chapter}{Bibliography}
%	\markboth{References}{References}
%	\bibliographystyle{alpha}
%	\bibliography{mybibVV}{}
	
	\thispagestyle{empty}
	
{	
		\markboth{References}{References}
		\bibliographystyle{alpha}
		\bibliography{mybibVV}{}
		
}

	\end{document}

%% file: packages.tex
\usepackage{pdfsync}
\usepackage{latexsym,enumitem}
\usepackage{amssymb}
\usepackage[cp850]{inputenc}
\usepackage{epsfig}
\usepackage{psfrag}
\usepackage{amsthm}
\usepackage{amscd}
\usepackage{amsmath}
\usepackage{amsfonts}
\usepackage{graphics,caption}
\usepackage[all]{xy}
\usepackage{etoolbox}
\patchcmd{\quote}{\rightmargin}{\leftmargin 2em \rightmargin}{}{}
 \usepackage[sort,numbers]{natbib} 
 \usepackage{mathtools}
\usepackage[normalem]{ulem}
\usepackage{array}

\usepackage[english]{babel}
\usepackage{comment}
\usepackage{color}
\usepackage[colorlinks]{hyperref}

\DeclareMathOperator{\tr}{tr} %stampa tr a testo
\let\phi\varphi

\DeclareMathOperator{\eqdef}{\coloneqq} %definition
\newcommand{\f}[2]{\frac{#1}{#2}} %Fractions

\let\epsilon\varepsilon
\let\subset\subseteq
\newcommand{\be}{\begin{equation*}}
 \newcommand{\ee}{\end{equation*}}
\newcommand{\bpf}{\begin{dimo}}
\newcommand{\epf}{\end{dimo}}
\newcommand{\bdefi}{\begin{defin}}
\newcommand{\edefi}{\end{defin}}
\newcommand{\bthm}{\begin{thm}}
\newcommand{\ethm}{\end{thm}}
\newcommand{\blem}{\begin{lem}}
\newcommand{\elem}{\end{lem}}
\newcommand{\bcor}{\begin{cor}}
\newcommand{\ecor}{\end{cor}}

\newcommand{\bprop}{\begin{prop}}
\newcommand{\eprop}{\end{prop}}
\newcommand{\bese}{\begin{ese}}
\newcommand{\eese}{\end{ese}}
\newcommand{\brem}{\begin{rem}}
\newcommand{\erem}{\end{rem}}
\newcommand{\bpfc}{\begin{dimoclaim}}
\newcommand{\epfc}{\end{dimoclaim}}

\newcommand{\abs}[1]{\left\lvert#1\right\rvert}						%mod
\newcommand{\norm}[1]{\left\lVert#1\right\rVert}					%norm
\newcommand{\quotient}[2]{\left.\raisebox{.1em}{$#1\!$}\middle/\raisebox{-.1em}{$#2$}\right.}

\DeclareMathOperator{\im}{Im} %image
\DeclareMathOperator{\id}{id} %identity map

\newenvironment{quot}
{
	\vspace{-0.2cm}
	\vspace{0.2cm}
}

\theoremstyle{definition}
\newtheorem{d1}{Definition}[section] %fa partire tutto da [section]

\newenvironment{defin}
{
	\begin{quot}
		\begin{d1}
		}
		{\end{d1}
	\end{quot}

}

\theoremstyle{definition}
\newtheorem{r1}[d1]{Remark}%[section]

\newenvironment{rem}
{
	\begin{quot}
		\begin{r1}
		}
		{\end{r1}
	\end{quot}
}

\theoremstyle{definition}
\newtheorem{e1}[d1]{Exercise}%[section]

\theoremstyle{definition}
\newtheorem{ese1}[d1]{Example}

\newenvironment{ese}
{
	\begin{quot}
		\begin{ese1}
	}
	{	
		\end{ese1}
	\end{quot}
}

\theoremstyle{definition}
\newtheorem{f1}[d1]{Formulation}

\theoremstyle{definition}
\newtheorem{f2}[d1]{Fact}

\theoremstyle{definition}

\theoremstyle{definition}

\theoremstyle{definition}
\newtheorem{t1}[d1]{Theorem}%[section]

\newenvironment{thm}
{
	\begin{quot}
		\begin{t1}}
		{\end{t1}
	\end{quot}
}

\theoremstyle{definition}
\newtheorem*{T1*}{Theorem}%[section]

\newenvironment{teor*}
{
	\begin{quot}
		\begin{T1*}}
		{\end{T1*}
	\end{quot}
}

\newenvironment{dimo}
{\begin{proof}[Proof]
	}
	{\end{proof}}

\newenvironment{dimoclaim}{\emph{Proof of Claim:}\;}{\hfill$\square$}

	\theoremstyle{definition}
	\newtheorem{l1}[d1]{Lemma}%[section]
	
	\newenvironment{lem}
	{
		\begin{quot}
			\begin{l1}}
			{\end{l1}
		\end{quot}
	}
	\theoremstyle{definition}
	\newtheorem{p1}[d1]{Proposition}%[section]
	
	\newenvironment{prop}
	{
		\begin{quot}
			\begin{p1}}
			{\end{p1}
		\end{quot}
	}
	
	\theoremstyle{definition}
	\newtheorem{c1}[d1]{Corollary}%[section]
	
	\newenvironment{cor}
	{
		\begin{quot}
			\begin{c1}}
			{\end{c1}
		\end{quot}
	}

		\renewenvironment{abstract}
	{\list{}{\rightmargin\leftmargin}%
		\item[\textbf{Abstract:}]\relax}
	{\endlist}

\newtheorem{innercustomthm}{Theorem}
\newenvironment{customthm}[1]
  {\renewcommand\theinnercustomthm{#1}\innercustomthm}
  {\endinnercustomthm}

 \newtheorem*{Theorem*}{Theorem}
 \newtheorem*{Proposition*}{Proposition}
 \newtheorem*{Lemma*}{Lemma}

\newtheorem{innercustomcor}{Corollary}
\newenvironment{customcor}[1]
  {\renewcommand\theinnercustomcor{#1}\innercustomcor}
  {\endinnercustomcor}